\documentclass[aop, preprint, a4paper]{imsart}

\usepackage[T1]{fontenc}	
\usepackage{totpages}
\def\CurrentOption{} %
\usepackage{graphicx}
\usepackage{amsmath, amsthm, amssymb}
\usepackage{upgreek}
\usepackage{mathrsfs}
\usepackage{pdfsync}
\usepackage[usenames,dvipsnames]{xcolor}
\usepackage[inline,shortlabels]{enumitem}
\usepackage{verbatim}								%
\usepackage{dsfont}									%
\usepackage{bbm}
\usepackage{stmaryrd}
\usepackage{multirow}

\RequirePackage{silence}
\usepackage{cite}

\usepackage{microtype}

\usepackage[font=small]{caption}

\usepackage{mathtools}
\usepackage{xstring}

\usepackage[smalltableaux,aligntableaux=center]{ytableau}

\usepackage{tikz}
\usepackage{tikz-cd}
\usetikzlibrary{calc, positioning, arrows, patterns}
\tikzcdset{arrow style=tikz}
\tikzset{>=latex}

\definecolor{grey1}{RGB}{71,71,71}

\RequirePackage{silence}
\usepackage[pdftex,bookmarks,
						colorlinks,				%
						breaklinks,unicode,
						citecolor=OliveGreen,
						linkcolor=Maroon
						]{hyperref} %

\usepackage{xparse}

\startlocaldefs

\ExplSyntaxOn
\NewDocumentCommand{\makeabbrev}{mmm}
 {
  \yoruk_makeabbrev:nnn { #1 } { #2 } { #3 }
 }

\cs_new_protected:Npn \yoruk_makeabbrev:nnn #1 #2 #3
 {
  \clist_map_inline:nn { #3 }
   {
    \cs_new_protected:cpn { #2 } { #1 { ##1 } }
   }
 }
 \ExplSyntaxOff

\makeabbrev{\textbf}{tbf#1}{a,b,c,d,e,f,g,h,i,j,k,l,m,n,o,p,q,r,s,t,u,v,w,x,y,z,A,B,C,D,E,F,G,H,I,J,K,L,M,N,O,P,Q,R,S,T,U,V,W,X,Y,Z}

\makeabbrev{\textbf}{bf#1}{a,b,c,d,e,f,g,h,i,j,k,l,m,n,o,p,q,r,s,t,u,v,w,x,y,z,A,B,C,D,E,F,G,H,I,J,K,L,M,N,O,P,Q,R,S,T,U,V,W,X,Y,Z}

\makeabbrev{\textsf}{tsf#1}{a,b,c,d,e,f,g,h,i,j,k,l,m,n,o,p,q,r,s,t,u,v,w,x,y,z,A,B,C,D,E,F,G,H,I,J,K,L,M,N,O,P,Q,R,S,T,U,V,W,X,Y,Z}

\makeabbrev{\mathsf}{mss#1}{a,b,c,d,e,f,g,h,i,j,k,l,m,n,o,p,q,r,s,t,u,v,w,x,y,z,A,B,C,D,E,F,G,H,I,J,K,L,M,N,O,P,Q,R,S,T,U,V,W,X,Y,Z}

\makeabbrev{\mathfrak}{mf#1}{a,b,c,d,e,f,g,h,i,j,k,l,m,n,o,p,q,r,s,t,u,v,w,x,y,z,A,B,C,D,E,F,G,H,I,J,K,L,M,N,O,P,Q,R,S,T,U,V,W,X,Y,Z,
sl,gl
}

\makeabbrev{\mathrm}{mrm#1}{a,b,c,d,e,f,g,h,i,j,k,l,m,n,o,p,q,r,s,t,u,v,w,x,y,z,A,B,C,D,E,F,G,H,I,J,K,L,M,N,O,P,Q,R,S,T,U,V,W,X,Y,Z}

\makeabbrev{\mathbf}{mbf#1}{a,b,c,d,e,f,g,h,i,j,k,l,m,n,o,p,q,r,s,t,u,v,w,x,y,z,A,B,C,D,E,F,G,H,I,J,K,L,M,N,O,P,Q,R,S,T,U,V,W,X,Y,Z}

\makeabbrev{\mathcal}{mc#1}{A,B,C,D,E,F,G,H,I,J,K,L,M,N,O,P,Q,R,S,T,U,V,W,X,Y,Z}

\makeabbrev{\mathbb}{mbb#1}{A,B,C,D,E,F,G,H,I,J,K,L,M,N,O,P,Q,R,S,T,U,V,W,X,Y,Z}

\makeabbrev{\mathscr}{ms#1}{A,B,C,D,E,F,G,H,I,J,K,L,M,N,O,P,Q,R,S,T,U,V,W,X,Y,Z}

\makeabbrev{\mathrm}{#1}{
Id,id,ran,rk,diag,stab,ann,conv,pr,ev,tr,End,Hom,sgn,im,op,can,fin,ext,red,Leb,lex,rsl,Aut,Inn,
rot,usc,lsc,Lip,lip,bSymLip,osc,AC,loc,z,all,
Opt,Adm,Cpl,Geo,GeoOpt,GeoAdm,GeoCpl,reg,res,graph,
bd,co,Ric,Exp,dExp,dist,seg,Seg,cut,fcut,Cut,SDiff,Iso,Isom,diam,cl,Homeo,Diff,Der,vol,dvol,inj,relint, Graph, sub,
var,law,Var,Poi,Gam,pa,iso,fs,inv,pqi,mix,erg,Cat,so,sort,
TestF,
}

\makeabbrev{\mathsf}{#1}{CD,BE,MCP,Ent,wMTW,MTW,Ch,RCD,EVI,Rad,dRad,SL,cSL,dSL,ScL,Irr,SC,wFe,VA,Dec}

\makeabbrev{\mathsc}{#1}{mmaf,cg}

\newcommand{\T}{\tau} %
\newcommand{\GP}[1]{\mcG_{#1}}

\newcommand{\DF}[1]{\mcD_{#1}}
\newcommand{\LP}{\mcL}

\newcommand{\Mbp}{\msM_b^+}

\newcommand{\imu}{\mrmi}

\newcommand{\proc}{\bullet}

\newcommand{\eps}{\varepsilon}

\newcommand{\defeq}{\eqqcolon}
\renewcommand{\complement}{\mathrm{c}}

\newcommand{\mathsc}[1]{\text{\textsc{#1}}}
\newcommand{\emparg}{{\,\cdot\,}}

\newcommand{\as}[1]{\quad #1\text{-a.s.}}

\newcommand{\Dir}[1]{D_{#1}}
\newcommand{\PD}[1]{\Pi_{#1}}
\newcommand{\Esf}[2]{E^{#1}_{#2}}
\newcommand{\pEsf}[2]{\mathrm{\msE}^{#1}_{#2}}
\newcommand{\rEsf}[2]{M^{\mathrm{Ewens}}_{#2;#1}}%

\newcommand{\probP}{\mrmP}
\newcommand{\ExpectP}{\mrmE}
\newcommand{\CondExpectP}[2]{\mrmE\left[#1\, \middle| \, #2\right]}
\newcommand{\probQ}{\mrmQ}

\newcommand{\BlockSupp}[1]{\Sigma(#1)}

\DeclareMathOperator{\eqdef}{\coloneqq}

\renewcommand{\deg}{\mathrm{deg}}

\newcommand{\supp}[1]{\mathsf{supp}#1}

\let\epsilon\varepsilon

\newcommand{\longrar}{\longrightarrow}

\newcommand{\nlim}{\lim_{n}}								%

\newcommand{\diff}{\mathop{}\!\mathrm{d}}						%

\newcommand{\tabs}[1]{\big\lvert#1\big\rvert}	
\newcommand{\abs}[1]{\left\lvert#1\right\rvert}						%

\newcommand{\floor}[1]{\left\lfloor#1\right\rfloor}						%

\newcommand{\norm}[1]{\left\lVert#1\right\rVert}					%
\newcommand{\set}[1]{\left\{#1\right\}}							%
\newcommand{\tset}[1]{\big\{#1\big\}}							%
\newcommand{\ttset}[1]{\{#1\}}									%
\newcommand{\paren}[1]{\left(#1\right)}							%
\newcommand{\tparen}[1]{\big({#1}\big)}
\newcommand{\ttparen}[1]{({#1})}
\newcommand{\bracket}[1]{\left[#1\right]}							%
\newcommand{\tbraket}[1]{\big[#1\big]}

\newcommand{\class}[2][]{\llbracket#2\rrbracket_{#1}}						%

\newcommand{\sym}[1]{{\scriptscriptstyle{(#1)}}}
\newcommand{\tgeq}{{\scriptscriptstyle{\geq}}}

\newcommand{\tbracket}[1]{\big[#1\big]}							%
\newcommand{\seq}[1]{\paren{#1}}								%
\newcommand{\tseq}[1]{{\big(#1\big)}}
\newcommand{\ttseq}[1]{(#1)}
\newcommand{\Cb}{\mcC_b}									%

\newcommand{\pfwd}{\sharp}
\newcommand{\mpfwd}{*}
\DeclareMathOperator*{\esssup}{esssup}

\DeclareMathOperator{\car}{\mathbf 1}

\newcommand{\CSimpl}{\Sigma}

\DeclareMathOperator{\emp}{\varnothing}
\newcommand{\N}{{\mathbb N}}
\newcommand{\R}{{\mathbb R}}

\newcommand{\DRG}[1]{\mathrm{DRG}_{#1}}
\newcommand{\CPS}[1]{\mathrm{PPS}_{#1}}

\newcommand{\restr}{\big\lvert}
\newcommand{\mrestr}[1]{\!\!\downharpoonright_{#1}}

\newcommand{\comma}{\,\mathrm{,}\;\,}

\newcommand{\fstop}{\,\mathrm{.}}

\DeclareMathOperator{\zero}{{\mathbf 0}}
\newcommand{\uno}{{\mathbf 1}}

\newcommand{\Bary}[1]{\mbfB_{#1}}

\newcommand{\Beta}{\mathrm{B}}

\newcommand{\boldalpha}{\boldsymbol\alpha}

\newcommand{\boldtheta}{\boldsymbol\theta}
\newcommand{\boldnu}{{\boldsymbol\nu}}
\newcommand{\boldlambda}{\boldsymbol\lambda}

\newcommand{\boldell}{\boldsymbol\ell}

\newcommand{\length}[1]{\abs{#1}}
\newcommand{\card}[1]{\#{#1}}

\newcommand{\coloring}{\mssC}
\newcommand{\Pmom}[1]{P_{#1}}
\newcommand{\DustV}{\mbfD}
\newcommand{\Dust}{D}
\newcommand{\BulkDustMap}[1]{\mbfS(#1)}
\newcommand{\Bulk}{\rho}

\newcommand{\had}{\diamond}
\newcommand{\Poch}[2]{\left\langle#1\right\rangle_{#2}}
\newcommand{\AntiPoch}[2]{#1^{\underline{#2}}}
\newcommand{\shape}{\mathsf{shape}}
\newcommand{\struct}{\mathsf{struct}}
\newcommand{\colore}{\mathsf{color}}

\newcommand{\col}[1]{\langle\!\!\langle#1\rangle\!\!\rangle}

\newcommand{\uni}[1]{\mathsf{dust}(#1)}

\newcommand{\size}[1]{\mathsf{s}(#1)}

\newcommand{\transpose}{\top}
\newcommand{\Multi}[1]{M}
\newcommand{\MultiSecond}{M}
\newcommand{\momNA}[1]{m_{#1}}
\newcommand{\mom}[2]{\momNA{#2}(#1)}
\newcommand{\momentumNA}[1]{\kappa_{#1}}
\newcommand{\momentum}[2]{\momentumNA{#2}(#1)}

\newcommand{\shortlex}{\mathrm{sl}}

\newcommand{\Part}[2]{\mcA^\sym{#1}_{#2}}
\newcommand{\PartM}[2]{\mcM^\sym{#1}_{#2}}
\newcommand{\KSimplex}[1]{\boldsymbol\nabla^\sym{#1}_{\!0}}
\newcommand{\OSimplex}[1]{\boldsymbol\nabla^\sym{#1}}

\newcommand{\IsoSimplex}{J}

\newcommand{\BackwardK}{\mssS}
\newcommand{\ForwardK}{\mssT}

\newcommand{\PotK}{\mssF}
\newcommand{\MartinK}{\mssK}

\let\temp\phi
\let\phi\varphi
\let\varphi\temp

\numberwithin{equation}{section}
\theoremstyle{plain}
\newtheorem{thm}{Theorem}[section]
\newtheorem*{thm*}{Theorem}
\newtheorem*{mthm*}{Main Theorem}
\newtheorem{thmIntro}{Theorem}
\newtheorem{propIntro}[thmIntro]{Proposition}
\newtheorem{prop}[thm]{Proposition}%
\newtheorem*{prop*}{Proposition}%
\newtheorem{lem}[thm]{Lemma}%
\newtheorem{cor}[thm]{Corollary}%

\theoremstyle{definition}
\newtheorem{defs}[thm]{Definition}%
\newtheorem*{defs*}{Definition}%

\newtheorem{rem}[thm]{Remark}%
\newtheorem{ese}[thm]{Example}%
\newtheorem*{ass*}{Assumption}%

\AtBeginDocument{%
  \setlength{\abovedisplayskip}{4pt plus 2pt minus 1pt}%
  \setlength{\belowdisplayskip}{4pt plus 2pt minus 1pt}%
  \setlength{\abovedisplayshortskip}{3pt plus 2pt}%
  \setlength{\belowdisplayshortskip}{3pt plus 2pt minus 1pt}
  }

\renewcommand{\paragraph}[1]{\smallskip\noindent\emph{#1}.\quad}

\allowdisplaybreaks

\endlocaldefs

\begin{document}

\begin{frontmatter}
\title{Polychromatic Partitions, Kingman Theory,\\and Ewens Sampling Formula}

\runtitle{Polychromatic Partitions, Kingman Theory, and ESF}
    
\begin{aug}
\author[A]{\fnms{Lorenzo}~\snm{Dello Schiavo}\ead[label=e1]{delloschiavo@mat.uniroma2.it}\orcid{0000-0002-9881-6870}}
\thanks{L.D.S.\ gratefully acknowledges funding by Austrian Science Fund (FWF) project \href{https://doi.org/10.55776/ESP208}{{10.55776/ESP208}}; and from PRIN Department of Excellence MatMod@TOV (CUP: E83C23000330006).}
\and
\author[B]{\fnms{Filippo}~\snm{Quattrocchi}\ead[label=e2]{filippo.quattrocchi@universite-paris-saclay.fr}\orcid{0009-0000-9773-1931}}
\thanks{F.Q.\ gratefully acknowledges support by the Austrian Science Fund (FWF) project \href{https://doi.org/10.55776/F65}{10.55776/F65} and from the Agence nationale de la recherche (ANR), through the PEPR PDE-AI project (ANR-23-PEIA-0004).}
\thanks{The authors are grateful to Professor~Nathana\"el Berestycki and Dr.~Alice Callegaro for several helpful suggestions, to Nicola Battisti and Dr.~Elizabeth Hollwey for discussions on DNA methylation, to Vladimir Badalyan for pointing out reference~\cite{FraLijPruReb25}.}
\address[A]{Dipartimento di Matematica, Università degli Studi di Roma ``Tor Vergata''%
\printead[presep={ ,\ }]{e1}}
\address[B]{Laboratoire de Mathématiques d'Orsay, Université Paris-Saclay\printead[presep={ ,\ }]{e2}}
\end{aug}

\begin{abstract}
A polychromatic partition is a partition of a set of colored (marked) elements in which each block is recorded only by the number of elements of each color in that block.

We introduce a notion of polychromatic partition structure, i.e.\ a family of random polychromatic partitions consistent under deletion of elements uniformly at random, simultaneously extending (standard) partition structures, multipartition structures, and partition structures of partially exchangeable type.
As an example thereof, we define a polychromatic analogue of the celebrated Ewens Sampling Formula, proving polychromatic versions of Kingman's characterization and of Hoppe's urn model.

As our main result, we give a complete characterization of polychromatic partition structures, proving a triple affine homeomorphism of Bauer simplices among: polychromatic partition structures, harmonic densities for the polychromatic Hoppe urn, and probability measures on a polychromatic Kingman simplex.
Additionally, we prove that this new simplex is homeomorphic to the Martin boundary of the polychromatic Hoppe urn inside its Martin compactification.

We discuss various applications including closed non-recursive non-iterative expressions for multivariate moments of several random measures in terms of cycle index polynomials whose monomials are indexed by polychromatic partitions.%
\end{abstract}



\end{frontmatter}

\setcounter{tocdepth}{3}

\tableofcontents

\section{Introduction}\label{s:Intro}
Kingman’s partition structures classify random partitions consistent under deletion of a uniformly chosen element. In many sampling problems, however, a block also carries an internal composition vector recording observable marks or colors.
This leads to partitions whose blocks have full color-count vectors.
The purpose of this paper is to develop the corresponding Kingman theory.

\paragraph{Integer partitions, set partitions, permutations}
Let~$[n]\eqdef\set{1,2,\dotsc, n}$. %
Partitioning~$[n]$ into disjoint blocks induces an (\emph{integer}) \emph{partition of~$n$}, recording the size of the blocks.
Throughout, we encode this partition by the vector $\boldlambda=\seq{\lambda_1,\dotsc,\lambda_n}$, where~$\lambda_i$ denotes the number of blocks of size~$i$.  
Thus,~$\sum_i i\lambda_i=n$, and we write $\boldlambda\vdash n$.

Integer partitions also classify \emph{permutations} of~$[n]$ according to their cycle structure.
The number of permutations sharing the same cycle structure is the multinomial coefficient~$\MultiSecond(\boldlambda)$ in~\eqref{eq:MultiSecond} which will play a central role throughout this work.

\subsection{Ewens Sampling Formula and Kingman's partition structures}\label{ss:ESFIntro}
Ewens~\cite{Ewe72} (cf.\ also Karlin--McGregor~\cite{KarMcG72}) computed the law of the allelic partition of a sample of~$n$ genes drawn from a large haploid population evolving under selectively neutral, infinitely-many-alleles mutation with scaled mutation rate~$\theta>0$, viz.
\begin{equation}
\label{eq:ESFIntro}
\Esf{n}{\theta}[\boldlambda]=\frac{n!}{\Poch{\theta}{n}} \prod_{i=1}^{n}\frac{\theta^{\lambda_i}}{i^{\lambda_i}\,\lambda_i!}
= \frac{\theta^{\lambda_1+\cdots+\lambda_n}}{\Poch{\theta}{n}} \MultiSecond(\boldlambda) \comma
    \qquad \boldlambda\vdash n \comma
\end{equation}
where $\Poch{\theta}{n}\eqdef \theta(\theta+1)\cdots(\theta+n-1)$ is the Pochhammer symbol.

The \emph{Ewens Sampling Formula}~\eqref{eq:ESFIntro} (in short: \emph{ESF}) also arises from the Wright--Fisher, Moran, and Cannings models, as explained \emph{a posteriori} by Kingman's coalescent theory~\cite{Kin82}, %
and it is ubiquitous in the combinatorics of random permutations and of logarithmic combinatorial structures; see e.g.\ the monograph~\cite{ArrBarTav03}. 
At~$\theta=1$, it gives the frequency~$\MultiSecond(\boldlambda)/n!$ of permutations of~$[n]$ with a given cycle structure~$\boldlambda$; for general~$\theta$ it is the law of the cycle structure of a $\theta$-biased permutation.

\paragraph{Partition structures}
In~\cite{Kin77,Kin78,Kin78b}, Kingman developed a structure theory around~\eqref{eq:ESFIntro}.
A \emph{partition structure} is a sequence~$\seq{\probP_n}_n$, with~$\probP_n$ a probability measure on the partitions of~$n$, such that, whenever~$n$ objects are partitioned according to $\probP_n$ and one object is subsequently deleted independently and uniformly at random, the partition of the remaining~$n-1$ objects is distributed according to~$\probP_{n-1}$. 

The ESF is the \emph{only} partition structure satisfying a recurrence relation of the form~\eqref{eq:KingmanDefinitionA},~\cite[\S3]{Kin78b}.
Furthermore, every partition structure is represented by one and only one Borel probability measure~$\mu$ on the ordered extended simplex~$\nabla_0$, as a mixture of extremal `paintbox' structures~\cite{Kin78}, the ESF corresponding to the \emph{Poisson--Dirichlet distribution} $\PD{\theta}$ of~\cite{Kin75}.

\subsection{Polychromatic partitions}
Suppose that each element of~$[n]$ is assigned one out~of~$q$ different colors.
Our main purpose is to describe (random) \emph{$q$-chromatic partitions of~$n$}; e.g.,
\[
\tset{\set{{\color{gray}{\bullet}}}, \set{\circ}, \set{{\color{gray}{\bullet}}\circ{\color{black}{\bullet}}}, \set{{\color{black}{\bullet}}\,\,{\color{black}{\bullet}}\,\circ\circ\,\,{\color{gray}{\bullet}}}, \set{{\color{gray}{\bullet}}\circ{\color{black}{\bullet}}\circ{\color{black}{\bullet}}}, \set{\circ\circ{\color{black}{\bullet}}\circ\circ\,\circ}}
\]
is a $3$-chromatic partition of~$21$ consisting of one block of type~$\tiny\ydiagram[*(gray)]{1}_1$, one~block of type~$\tiny\ydiagram[*(white)]{1}_1$, one~block of type~$\tiny\ydiagram[*(black)]{1}_1~\tiny\ydiagram[*(gray)]{1}_1~\tiny\ydiagram[*(white)]{1}_1$, \emph{two} blocks of type~$\tiny\ydiagram[*(black)]{1}_2~\tiny\ydiagram[*(gray)]{1}_1~\tiny\ydiagram[*(white)]{1}_2$, and one block of type~$\tiny\ydiagram[*(black)]{1}_1~\tiny\ydiagram[*(white)]{1}_5$.
Here, for colors~$\tiny\ydiagram[*(black)]{1}$, $\tiny\ydiagram[*(gray)]{1}$, and~$\tiny\ydiagram[*(white)]{1}$, the phrasing `\emph{type~$\tiny\ydiagram[*(black)]{1}_a~\tiny\ydiagram[*(gray)]{1}_b~\tiny\ydiagram[*(white)]{1}_c$}' indicates a block with $a$~elements of color~$\tiny\ydiagram[*(black)]{1}$, $b$~elements of color~$\tiny\ydiagram[*(gray)]{1}$, and $c$~elements of color~$\tiny\ydiagram[*(white)]{1}$.
Let us note explicitly that a polychromatic partition accounts only for the \emph{number} of elements of each color in each block, but not for their label as elements of~$[n]$: the blocks~$\set{{\color{black}{\bullet}}\,{\color{black}{\bullet}}\circ\circ\,{\color{gray}{\bullet}}}$ and~$\set{{\color{gray}{\bullet}}\circ{\color{black}{\bullet}}\circ{\color{black}{\bullet}}}$ have the same type and are thus indistinguishable. 
(That is, we distinguish neither the elements' labels nor the order in which the colors appear.)

It is often convenient to keep track of the total number of elements of~$[n]$ painted in each color.
To this end, let~$\mbfn\eqdef\seq{n_1,\dotsc, n_q}$ be non-negative integers and set~$\length{\mbfn}\eqdef n_1+\cdots + n_q$.
We call \emph{$\mbfn$-chromatic partition} any $q$-chromatic partition of~$\length{\mbfn}$ totaling $n_j$ elements of color~$j$ across all blocks, for each~$j\in [q]$; see Fig.~\ref{fig:LongList}.%

\paragraph{Multisets}
The following discussion rests on the encoding of polychromatic partitions by multisets.
Set $\N^q_*\eqdef \N_0^q\setminus\set{\zero}$.
A $q$-chromatic partition is a (finite finitely supported) $\N^q_*$-multiset, that is, a map $A\colon\N^q_*\to\N_0$, the value $A(\mbfa)$ denoting the number of blocks whose color count is $\mbfa$.
We write: $\car_{\mbfa}$ for the multiset consisting of the single element $\mbfa$, so that $A-\car_{\mbfa}$ is~$A$ with one block of type $\mbfa$ removed; $\supp{A}$ for the support of $A$; $\size{A}\eqdef \sum_{\mbfa} A(\mbfa)$ for the total number of blocks; $\col{A}\eqdef \sum_{\mbfa}A(\mbfa)\,\mbfa$ for the vector of total color counts; $A\vDash\mbfn$ when $\col{A}=\mbfn$; $\shape(A)$ for the integer partition of~$\length{\mbfn}$ obtained by forgetting the coloring; and $\Part{q}{n}$ for the set of all $q$-chromatic partitions of~$n$.
Finally, the \emph{multinomial coefficient} of $A\vDash\mbfn$ is, cf.~Definition~\ref{d:MultinomialCoeff},
\begin{align}\label{eq:MultiIntro}
\Multi{\mbfn}(A)\eqdef&\ n_1!\cdots n_q! \prod_{\mbfa\in \supp{A}} \frac{\binom{\length{\mbfa}}{\mbfa}^{A(\mbfa)}}{\length{\mbfa}^{A(\mbfa)} A(\mbfa)!}\comma \qquad \binom{\length{\mbfa}}{\mbfa} \eqdef \binom{a_1+\cdots + a_q}{a_1,\dotsc, a_q} \fstop
\end{align}

\subsection{Main results}
We briefly illustrate our main results, focusing on probabilistic ones.
For combinatorial results, precise statements and references will be presented in the next sections, once the necessary notation and definitions have been introduced.

\subsubsection{Probability}
We present polychromatic analogues of~\cite{Kin77,Kin78,Kin78b,Hop84,LzDS19a} and others, as well as extensions to polychromatic partitions of recent results for multipartitions~\cite{Str24a,Str26}.
In particular, we define a \emph{polychromatic ESF} (see below), a \emph{polychromatic Hoppe Chain} (PHC, see~\S\ref{ss:Hoppe}), and polychromatic versions of Kingman's Characterization and Representation Theorems, completing the latter with a full characterization of the PHC \emph{Martin boundary}. %

\paragraph{Consistency}
We introduce a notion of \emph{consistency} for random polychromatic partitions (Dfn.~\ref{d:Consistency}), simultaneously extending:
\begin{itemize}
\item Kingman's \emph{partition structures}~\cite{Kin78,Kin78b};
\item Strahov's \emph{multipartition structures}~\cite{Str24a,Str26} (cf.~\eqref{eq:Multipartition});
\item the polychromatic partitions \emph{of partially exchangeable type} (Dfn.~\ref{d:PEPS}) arising from the \emph{partially exchangeable partition probability functions} in~\cite{FraLijPruReb25}.
\end{itemize}
Precisely, we call \emph{$q$-chromatic partition structure} ($\CPS{q}$) any sequence $\probP_\proc\eqdef\seq{\probP_n}_n$ with $\probP_n$ a probability on~$\Part{q}{n}$ and satisfying $\probP_m=\sigma^\sym{q}_{m,n}\probP_n$ for all $m\le n$, where $\sigma^\sym{q}$ is the projective system~\eqref{eq:Consistency} induced by deletion of one element chosen uniformly at random.

This notion of consistency is stable under the identification of colors (Prop.~\ref{p:ConsistencyDegeneracy}).

\paragraph{Representation}
We prove a Representation Theorem for $\CPS{q}$'s, analogous to~\cite[Theorem, p.~376]{Kin78}, as we now explain.
For each~$n$, let~$\Part{q}{n}$ be the set of $q$-chromatic partitions of~$n$, see Definition~\ref{d:MultiYoung}.
Further let~$\KSimplex{q}\subset ([0,1]^q)^{\N_0}$ be the \emph{$q$-chromatic Kingman simplex}~\eqref{eq:KingmanSimplex} and, for each $n$ and each~$A\in\Part{q}{n}$, let~$\psi_A\colon \KSimplex{q}\to [0,1]$ be the polynomial~\eqref{eq:DefPsi}.

\begin{thmIntro}[Polychromatic Kingman Representation, see Thm.~\ref{t:KingmanRepresentation}]\label{t:KingmanIntro}
For every Borel probability measure~$\mu$ on~$\KSimplex{q}$ and every~$n \in \N_0$, the assignment
\begin{equation}\label{eq:p:KingmanRepresentationConverseIntro:0}
\probP_n\colon A\longmapsto \int_{\KSimplex{q}} \psi_A(\mbfx) \diff\mu(\mbfx) \comma \qquad A\in\Part{q}{n}\comma
\end{equation}
defines a probability measure~$\probP_n$ on~$\Part{q}{n}$, and the family~$\seq{\probP_n}_n$ is a $q$-chromatic partition structure.
Conversely, every $q$-chromatic partition structure is represented by~\eqref{eq:p:KingmanRepresentationConverseIntro:0} for one and only one Borel probability measure~$\mu$ on~$\KSimplex{q}$.
\end{thmIntro}

We will detail the proof strategy in~\S\ref{sss:KingmanRepresentationDiscussion} below.
We further show that the representation is functorial in the sense of Proposition~\ref{p:Functoriality}.
This property is a genuine feature of the polychromatic picture, void in the monochromatic case.
We collect several examples of $\CPS{q}$'s and their corresponding~$\mu$'s in \S\S\ref{sss:Strahov} and~\ref{sss:FurtherExamples}.

\paragraph{Doob--Martin boundary and ergodic theory}
Our proof of the polychromatic Representation Theorem is different from known proofs in the monochromatic case. 
Firstly, at fixed coloring, polychromatic partition structures are not partially exchangeable sequences in the standard sense, so that Aldous' proof~\cite{Ald85} does not apply;
{when the coloring is itself random, the exchangeability of polychromatic partitions on the extended space~$\N_1\times [q]$ only yields the existence of a de~Finetti measure on some abstract indexing space.}
Secondly, the multivariate nature of the block labeling requires, upon deletion of one element, a reordering of the blocks, so that the proportion of elements of a given color in %
the $k$ largest blocks need not be a reverse sub-martingale, which breaks Kingman's original argument~\cite{Kin78,Kin78b}.
Thirdly, for~$q\geq 2$, the reordering map onto~$\KSimplex{q}$ is \emph{not} continuous for the natural (\emph{shortlex}, see Dfn.~\eqref{eq:ShortlexOrder}) order on~$([0,1]^q)^{\N_0}$.
We overcome this obstacle by representing points in~$\KSimplex{q}$ as purely atomic probabilities on~$\Delta^{q-1}$, see~\S\ref{sss:MartinCompactness}.

We thus prove the Theorem via potential theory for the PHC, completely characterizing its Martin compactification and boundary, and the class of its harmonic functions, obtaining as well the following equivalence of Bauer simplices:

\begin{thmIntro}[Martin correspondence, Cor.~\ref{c:MartinCorrespondence}]\label{t:MartinIntro}
Endowed with an appropriate compact Polish topology~$\T_*$, the polychromatic Kingman simplex~$\KSimplex{q}$ is (homeomorphic to) the Martin boundary of the PHC.
In particular, there is a triple continuous affine one-to-one correspondence among
\begin{equation*}
\set{\parbox{0.23\textwidth}{\centering polychromatic\\ partition structures}} \longleftrightarrow \set{\parbox{0.22\textwidth}{\centering Borel probability\\measures on~$\KSimplex{q}$}} \longleftrightarrow \set{\parbox{0.18\textwidth}{\centering PHC-harmonic\\densities}} \fstop
\end{equation*}
Under this correspondence, the extremal $\CPS{q}$'s~$\psi_\proc(\mbfx)$ are identified with the pure states~$\delta_\mbfx$ on~$\KSimplex{q}$ and with minimal PHC-harmonic densities.

Furthermore, for a~$\CPS{q}$~$\probP_\proc$ the following are equivalent:
\begin{itemize}
\item $\probP_\proc$ is extremal, i.e.~$\probP_\proc=\psi_\proc(\mbfx)$ for some~$\mbfx\in\KSimplex{q}$;
\item the tail $\sigma$-algebra generated by~$\probP_\proc$ is trivial;
\item empirical block-color frequencies converge $\probP_\proc$-a.s.\ to a deterministic~$\mbfx\in\KSimplex{q}$.
\end{itemize}
\end{thmIntro}
(See~\S\ref{sss:MartinCompactness} for the definition and properties of~$\T_*$ on~$\KSimplex{q}$, \S\ref{sss:PHCRevisited} for results on the PHC, \S\ref{sss:MartinConvergence} for the convergence of Martin kernels, and~\S\ref{sss:Minimality} for the minimality assertion.
The characterization of tail-triviality is obtained as a byproduct of the proof.)

Write~$\mbfx\eqdef \seq{\mbfx_0,\mbfx_1,\mbfx_2,\dotsc}\in\KSimplex{q}$ with~$\mbfx_i\in [0,1]^q$ for all~$i$, and~$\mbfe_j$ for the $j^\text{th}$ vector of the canonical basis of~$\R^q$.
Summing the entries of~$\mbfx$ returns a vector~$\DustV(\mbfx)\eqdef \mbfx_0+\mbfx_1+\cdots$, the \emph{total frequency of each color} in the population described by~$\mbfx$.
Thus, a law~$\DustV_\pfwd\mu$ of color frequencies is assigned to every polychromatic partition structure.

Let~$\Delta^{q-1}$ be the standard $q$-simplex.
We have the following specialization to~$\CPS{q}$'s of special type:

\begin{propIntro}
Let~$\probP_\bullet$ be a~$\CPS{q}$ represented by~$\mu$ in the sense of Theorem~\ref{t:KingmanIntro}. Then,
\begin{itemize}
\item $\probP_\bullet$ is a standard partition structure if and only if~$\mbfx_i\propto \mbfe_j$ for every~$i\geq 0$ $\mu$-a.s., for some fixed~$j\in [q]$ (trivial);
\item $\probP_\bullet$ is a multipartition structure in the sense of~\cite{Str24a} if and only if~$\mbfx_i\propto \mbfe_{j_i}$ for some~$j_i=j_i(\mbfx)\in [q]$ depending on~$i$ and on~$\mbfx$, for every~$i\geq 1$, $\mu$-a.s.\ (Prop.~\ref{p:RepresentationStrahov});
\item $\probP_\bullet$ is of partially exchangeable type if and only if~$\DustV_\pfwd \mu= \delta_\mbfp$ for some~$\mbfp\in\Delta^{q-1}$, {where~$\DustV\colon \KSimplex{q}\to \Delta^{q-1}$ is the map in~\eqref{eq:IsoSimplex}} (Prop.~\ref{p:DiracCriterion}) .
\end{itemize}
\end{propIntro}

\paragraph{The polychromatic Ewens Sampling Formula}
Fix $\theta>0$ and a Borel probability measure~$\nu$ on the simplex~$\Delta^{q-1}$, and write $\mom{\mbfa}{\nu}\eqdef \int_{\Delta^{q-1}}u_1^{a_1}\cdots u_q^{a_q}\diff\nu(\mbfu)$
for the moments of $\nu$.
We introduce the \emph{polychromatic ESF} (in short: \emph{pESF}) with parameters~$\theta$ and~$\nu$, viz.
\begin{equation}\label{eq:pESFIntro}
\pEsf{n}{\theta,\nu}[A] \eqdef \frac{n!\,\theta^{\size{A}}}{\Poch{\theta}{n}}\, \bracket{\prod_{\mbfa\in\supp A} \mom{\mbfa}{\nu}^{A(\mbfa)}}\, \frac{\Multi{\col{A}}(A)}{\col{A}!}
\comma \qquad A\in \Part{q}{n} \comma
\end{equation}
which is a probability distribution on $\Part{q}{n}$ (Prop.~\ref{p:ESFProbab}) and defines a $q$-chromatic partition structure (Prop.~\ref{p:RepresentationPESF}).
The parameter~$\theta$ has the same interpretation as in the ESF, while the parameter measure~$\nu$ determines the frequency of colors:
\begin{itemize}
\item when~$\nu=\delta_{\mbfe_j}$ for some~$j\in [q]$, the pESF reduces to the standard Ewens Sampling Formula, as is readily seen comparing~\eqref{eq:pESFIntro} and~\eqref{eq:MultiIntro} with~\eqref{eq:ESFIntro}; also cf.\ Remark~\ref{r:UnidimensionalSecond};
\item when~$\nu=\delta_\mbfp$ for some~$\mbfp\in\Delta^{q-1}$, then~$\mom{\mbfa}{\nu}=p_1^{a_1}\cdots p_q^{a_q}$ and each individual acquires marker~$j$ independently with probability~$p_j$. In this case (and only in this case) the pESF is \emph{of partially exchangeable type}; see Corollary~\ref{c:PESFNotPE};
\item if~$\nu$ is supported only on the vertices of the simplex~$\Delta^{q-1}$, the pESF reduces to Strahov's \emph{refined Ewens Sampling Formula}~\cite{Str26} (rESF); see Proposition~\ref{p:StrahovReduction}.
\end{itemize}

The representing measure for the pESF (isomorphic to the law of a Dirichlet--Ferguson process with intensity~$\theta\nu$) is discussed in Proposition~\ref{p:RepresentationPESF};
its asymptotics as~$\theta\to0,\infty$ in Corollary~\ref{c:Asymptotics}

\paragraph{Characterization} 
We prove the polychromatic analogue of Kingman's characterization in~\cite[\S3, p.~5]{Kin78b}. 
Here we state a simplified version; see Theorem~\ref{t:CharacterizationNew} for a more general statement under \emph{necessary and sufficient assumptions} on~$\probP_\proc$.
The special case of the rESF is noted in Proposition~\ref{p:CharacterizationRESF}.
\begin{thmIntro}[\ref{t:CharacterizationNew}]
Let~$\probP_\proc$ be a $q$-chromatic partition structure satisfying
\begin{enumerate}[$(a)$, wide]
\item $\probP_n[A]>0$ for every~$n\in\N_1$ and every~$A\in\Part{q}{n}$; 
\item for every~$n\in\N_1$, every~$A\in\Part{q}{n}$, and every~$\mbfr\in \N^q_*$ with~$\car_\mbfr \leq A$, it holds that
\begin{equation}\label{eq:PESFRecursionIntro}
\length{\mbfr}\, A(\mbfr) \, \probP_n[A] = b(n,\mbfr)\, \probP_{n-\length{\mbfr}}[A-\car_\mbfr] \comma 
\end{equation}
for some constant~$b(n,\mbfr)$ independent of~$A$. 
\end{enumerate}
Then,~$\probP_n=\pEsf{n}{\theta,\nu}$ for a unique~$\theta>0$ and a unique Borel probability measure~$\nu$ on~$\Delta^{q-1}$.
\end{thmIntro}

Let us note that the coefficients~$b(n,\mbfr)$ above are entirely unspecified; the pESF satisfies~\eqref{eq:PESFRecursionIntro} with the explicit choice~\eqref{eq:KingmanDefinitionB} for~$b(n,\mbfr)$, see Proposition~\ref{p:pESFRecurrence}.
As an effect of the coloring, the proof is different from Kingman's and {additionally requires a characterization of harmonic functions for the $\mbfp$-biased random walk on the $q$-dimensional directed Pascal graph~$\N_0^q$. %
This fact is known for the Pascal triangle (i.e.~$q=2$, see e.g.~\cite[first Ex.\ on p.~5679]{GnePit06} and references therein) and for the uniform case~$\mbfp=(\frac{1}{q},\dotsc,\frac{1}{q})$ (see~\cite[Ex.~6]{GerGruHag16}).}

\paragraph{Related work}
In the recent work~\cite{FraLijPruReb25}, the authors develop a general theory of \emph{partially exchangeable partition probability functions} (pEPPF, previously considered in specific instances in~\cite{LeiLij11,LijNipPru14}), and use it to characterize the random partition induced by a partially exchangeable array.
(We refer to~\cite{FraLijPruReb25} for a comparison between this notion and the one of partially exchangeable random partition by Pitman~\cite{Pit95}.)
Our Representation Theorem~\ref{t:KingmanIntro} may be used to infer a {polychromatic Kingman paintbox} (cf.~\cite[\S1.1]{Ber09}) as well as the characterization of partially exchangeable partition structures in~\cite[Thm.~5]{FraLijPruReb25} (see~\S\ref{sss:PEPS}), thus providing a direct construction of pEPPF's via $\CPS{q}$'s. The lack of such construction, in comparison with the exchangeable case, is noted in~\cite[p.~17]{FraLijPruReb25}.
We note that Theorem~\ref{t:KingmanIntro} and~\cite[Thm.~5]{FraLijPruReb25} are not equivalent, for a structural reason which we can already state here.
We show (Prop.~\ref{p:DiracCriterion}) that a polychromatic partition structure comes from a partially exchangeable array precisely when the law $\DustV_\pfwd\mu$ of block-color frequencies is a Dirac mass~$\delta_\mbfp$, that is, precisely when the color frequencies are \emph{deterministic}.
Partially exchangeable partition structures are thus (indexed by) elements of the fibers of the single map~$\DustV$, and a general polychromatic partition structure is a mixture of them over~$\mbfp\in\Delta^{q-1}$.

The inclusion is strict already for the objects motivating this work.
Under $\pEsf{n}{\theta,\nu}$, each block has its own color distribution drawn according to~$\nu$, so the color frequencies of the whole population are random unless~$\nu$ is a Dirac mass.
In particular, Strahov's rESF, obtained for~$\nu$ supported on the vertices of~$\Delta^{q-1}$, has Dirichlet distributed color frequencies and lies outside the scope of~\cite{FraLijPruReb25}; see Cor.~\ref{c:PESFNotPE}.

Furthermore, the shortlex ordering in~\eqref{eq:KingmanSimplex} fixes the labeling of the atoms in~\eqref{eq:IsoSimplex}, thereby removing the ambiguity left by their unlabeled weights and ensuring uniqueness of the representing measure~$\mu$.

Finally, whereas the affine bijections in Theorem~\ref{t:MartinIntro} and the
identification of their extreme points are measure-theoretic (and are deduced
here from Theorem~\ref{t:KingmanIntro}, cf.\ Remark~\ref{r:LogicalDependencies}),
all topological statements --- the homeomorphism with the Martin boundary, the
Martin compactification, the continuity of the correspondence, and the Bauer
simplex structure --- rely on the compact topology~$\T_*$
of~\S\ref{sss:MartinCompactness} and are thus not accessible by purely
measure-theoretic techniques.

\paragraph{Monochromatic case}
When~$q=1$, %
the topology~$\T_*$ reduces to the topology on~$\nabla_0$ discussed by Kingman in~\cite[p.~379]{Kin78}.
The compactness of~$\nabla_0$ in this topology is stated in passing in~\cite[p.~379]{Kin78}, while the continuity of Kingman's correspondence is shown in~\cite[\S5]{Kin78}.
It was shown by Kerov, Okounkov, and Olshanski in~\cite[\S4]{KerOkoOls98} that the boundary of the {Kingman graph} coincides with~$\nabla_0$ and can be described as the Martin boundary induced by the Martin kernel of the \emph{dimension function}~\cite[Eqn.~(2.2)]{KerOkoOls98}.
Our proof recovers {at once} all these results.

\subsubsection{Combinatorics} 
We compute the generating function~$Z^\sym{q}_\mbfn$ for the frequency of $\mbfn$-chromatic set partitions of~$[\length{\mbfn}]$ %
as a generalized \emph{De~Bruijn's pattern inventory}. %
By efficiently encoding polychromatic partitions as \emph{multisets}, we further provide a robust and versatile framework for counting {and ordering} similarly polychromatic objects, including polychromatic permutations, Young diagrams, and necklaces.
We show that the generating function of~$\pEsf{n}{\theta,\nu}[A]$ {is expressed in terms of~$Z^\sym{q}_\mbfn$} for~$A\vDash\mbfn$, which allows us to easily compute the pESF in the extremal cases~$\theta\in\set{0,\infty}$.

\subsubsection{Applications}
As a theoretical application of polychromatic partitions, \S\ref{s:ApplMoments} provides an effective way to compute all multivariate moments of Dirichlet, Dirichlet--Ferguson, and Gamma random measures in terms of the cycle index polynomials~{$Z^\sym{q}_\mbfn$.}

In Appendix~\ref{ss:PopGen}, we show how the pESF can be used to describe allele frequencies in a population when also accounting for epigenetic modifications such as DNA methylation.
Indeed, we have the following correspondence between the objects above and the sampling of a population subject to \emph{in situ} epigenetic modifications:
\begin{itemize}
\item allelic type $\rightsquigarrow$ block of the polychromatic partition;
\item epigenetic marker $\rightsquigarrow$ color of an element;
\item allele-marker group $\rightsquigarrow$ entry $a_j$ of the color count $\mbfa$ of a block;
\item removal of one individual chosen uniformly at random from the sample $\rightsquigarrow$ the kernels~\eqref{eq:Kernels} and the consistency relation~\eqref{eq:Consistency};
\item inability to resolve (all or some) markers~$\rightsquigarrow$ the degeneracy maps of~\S\ref{sss:Degeneracy}.
\end{itemize}

\section{Terminology and notation}
We use the term `chromatic partition' (and, in general, the term `chromatic', rather than `colored') to distinguish the objects we consider from `colored partitions' (also: `prefabs', see~\cite[\S3.14, p.~92]{Wil94}; `multipartitions', see~\cite{And08}; `multiple partitions', see~\cite{Str26,Str24a}), i.e.\ partitions in which each \emph{block} is assigned one of a fixed set of colors.
Whereas, in the probabilistic literature, the terminology of `colored partitions' appears to be more common than that of `multipartitions', in order to further avoid confusion, from here onwards we shall only use the term `multipartition' to refer to these objects.

Multipartitions form a special (quite small) subset of polychromatic partitions, see Remark~\ref{r:ColoredVsPolychromatic} below.
Also, polychromatic partitions cannot be translated into semistandard (nor standard) Young tableaux; see Remark~\ref{r:Tableaux}.

\subsection{Notation}
The following shorthand notation greatly facilitates the presentation of concise proofs.
Let~$\N_1\eqdef \set{1,2,\dotsc}$,~$\N_0\eqdef \N_1\cup \set{0}$, and~$\overline\N_1\eqdef \N_1 \cup \set{\infty}$,~$\overline\N_0\eqdef \overline\N_1\cup\set{0}$.
Let~$\Gamma$ be \emph{Euler's Gamma function}, $\Poch{\alpha}{k}\eqdef \Gamma(\alpha+k)/\Gamma(\alpha)$ be the \emph{Pochhammer symbol} (raising factorial) and~$\AntiPoch{\alpha}{k}\eqdef \Poch{\alpha-k+1}{k}$ the \emph{falling factorial} of~$\alpha>0$.
Finally, let~$\Beta(\mbfx)\eqdef \Gamma(x_1)\cdots\Gamma(x_k)/\Gamma(x_1+\cdots+x_k)$ be the \emph{multivariate Beta function}.
Conventionally, throughout this work,~$0^0\eqdef 1$.

\paragraph{Permutations}
For~$n\in\N_1$ let~$\mfS_n$ be the symmetric group of degree~$n$, naturally acting on~$[n]$ by permutation of its elements.
Usually, we represent a permutation~$\pi$ with $r$~cycles in its cycle notation, viz.
\begin{equation}\label{eq:StdRepPi}
\pi=\seq{y_{1,1}\ \ y_{1,2}\ \cdots)(y_{2,1}\ \ y_{2,2}\ \cdots)\cdots (y_{r,1}\ \ y_{r,2}\ \cdots} \fstop
\end{equation}
Contrary to the standard convention, 
\begin{center}
$r=r_\pi$ is the total number of cycles of~$\pi$, \textbf{including fixed points}.
\end{center}
For a permutation~$\pi\in\mfS_n$, we let~$\boldlambda_\pi\vdash n$ count the number of cycles in~$\pi$ of given length:~$\lambda_{\pi,1}$ is the number of fixed points of~$\pi$, $\lambda_{\pi,2}$ the number of $2$-cycles in~$\pi$, and so on.

\paragraph{Partitions}
$\boldlambda\in \N_0^n$ is called a (\emph{integer}) \emph{partition} of~$n$ if~$\sum_i i\lambda_i=n$, in which case we write~$\boldlambda\vdash n$.
We call \emph{multinomial coefficient~$\MultiSecond(\boldlambda)$ of~$\boldlambda\vdash n$} the number of permutations in~$\mfS_n$ with $\lambda_i$ $i$-cycles for each~$i$, viz.\ (See~\cite[Vol.~1, Prop.~1.3.2]{Sta01}.)%
\begin{align}\label{eq:MultiSecond}
\MultiSecond(\boldlambda)\eqdef n! \prod_{i=1}^n\frac{1}{i^{\lambda_i} \lambda_i!}\fstop
\end{align}

\paragraph{Multisets}
Given a set~$S$, an \emph{$S$-multiset} is any map~$m\colon S\to \overline\N_0$. 
We denote by~$\supp\,{m}$ the \emph{support} of~$m$.
The \emph{size}~$\size{m}$ of~$m$ is the sum of all its values, viz.
\[
\size{m}\eqdef \sum_{s\in S}m(s) \fstop
\]

Given a map~$f\colon S\to T$, the \emph{push-forward via~$f$} of an $S$-multiset~$m$ is the $T$-multiset 
\begin{equation}\label{eq:MultiPush}
f_\mpfwd m \eqdef \sum_{s\in\supp\,{m}} m(s) \car_{f(s)} \fstop
\end{equation}

\paragraph{Vectors} 
Whenever no confusion may arise, we do not distinguish between row and column vectors.
Let~$\mbfe_i%
$ be the $i^\textrm{th}$ vector of the canonical basis of~$\R^n$ ($n$~is apparent from context and therefore omitted from the notation), and set~$\zero\eqdef\seq{0, \dotsc, 0}$ and~$\uno\eqdef \seq{1,\dotsc, 1}$.
For vectors~$\mbfx,\mbfy\in \R^n$ and~$m\in\N_0$, write%
\begin{align*}
\mbfx\cdot\mbfy\eqdef&\ x_1y_1+\cdots + x_ny_n\comma & \mbfx\had\mbfy\eqdef& \seq{x_1y_1,\dotsc, x_ny_n}\comma
\\
\mbfx^{\had m}\eqdef& \seq{x_1^m,\dotsc, x_n^m} \comma
& \length{\mbfx}\eqdef&\ \abs{x_1}+\cdots + \abs{x_n}\fstop %
\end{align*}
For vectors~$\mbfx,\mbfy\in\R^n$, we write
\[
\mbfx\leq_\had\mbfy \quad\iff\quad x_i\leq y_i \ \text{ for each }\ i\in [n] \fstop
\]

\paragraph{Matrices}
For a matrix~$\mbfM\eqdef \tbraket{m_{i,j}}_{i\in[a], j\in [b]}\in \R^{a\times b}$ ($a$ rows,~$b$ columns) set
\begin{align*}
\mbfM_i\eqdef&\ \mbfe_i^\transpose\, \mbfM \comma & \mbfM^j\eqdef&\ \mbfM\, \mbfe_j \comma
\end{align*}
In words:~$\mbfM_i$ is the $i^\textrm{th}$ row of~$\mbfM$ and~$\mbfM^j$ is the $j^\textrm{th}$ column of~$\mbfM$. %
For matrices (or vectors)~$\mbfS,\mbfM\in\R^{a\times b}$ set~$\mbfS^\mbfM\eqdef\prod_{i,j}^{a,b} s_{i,j}^{m_{i,j}}$ and~$\mbfM!\eqdef\prod_{i,j}^{a,b} m_{i,j}!$.

\paragraph{Simplex}
We denote by~$\Delta^{k-1}$ the \emph{standard simplex}
\begin{equation*}%
\Delta^{k-1}\eqdef \tset{\mbfx\in \R^k : \mbfx\geq_\had \zero\comma\length{\mbfx}=1} \fstop
\end{equation*}

\section{Counting pattern inventories}\label{s:Patterns}
Throughout,~$n\in \N_0$, $q,r\in\N_1$, and~$\mbfn\eqdef\seq{n_1,\dotsc,n_q}\in\N_0^q$.
In the next definition, let~$S$ be any non-empty set.

\begin{defs}\label{d:Category}
A \emph{$q$-coloring} of~$S$ is any function~$\mssC\colon S\to [q]$, assigning to each element of~$S$ a color in~$[q]$.
An \emph{$\mbfn$-coloring} is any $q$-coloring~$\mssC$ of~$[\length{\mbfn}]$ with~$\card{\mssC^{-1}(j)}=n_j$ for~$j\in [q]$.
\end{defs}

If not otherwise stated, everywhere in the following we fix an $\mbfn$-coloring~$\mssC_\mbfn$.
Assertions which depend on~$\mssC_\mbfn$ do so only via the choice of~$\mbfn$, which is always clear from the notation.

\subsection{Polychromatic partitions}
Set
$\N^q_*\eqdef \N_0^q\setminus \set{\zero}$. %
For a finite finitely supported $\N^q_*$-multiset~$A$, we set
\begin{equation}
\col{A}\eqdef \sum_{\mbfa\in \supp{A}} A(\mbfa)\, \mbfa  \quad \in \N^q_* 
\fstop
\end{equation}

\begin{defs}[Polychromatic partitions]\label{d:MultiYoung}
An $\N_*^q$-multi\-set~$A$ is called a \emph{$q$-chromatic partition of~$\mbfn$} if~$\col{A}= \mbfn$, in which case we write~$A\vDash\mbfn$.
We denote by~$\mcA_\mbfn$ the family of all $q$-chromatic partitions of~$\mbfn$, and we set
\begin{align}\label{eq:MultiSets}
\Part{q}{0}\eqdef\set{0}\comma \quad \Part{q}{n}\eqdef \bigcup_{\mbfn\in \N^q_*: \length{\mbfn}=n} \mcA_\mbfn \comma \quad n\in \N_1\comma \qquad \text{and} \qquad \Part{q}{\cup}\eqdef \bigcup_{n\geq 0} \Part{q}{n} \fstop
\end{align}
\end{defs}

Note that~$\Part{q}{\cup}$ is the family of all finitely supported finite $\N^q_*$-multisets.
Also, by definition, a $q$-chromatic partition~$A\in\Part{q}{n}$ is a partition of elements colored with \emph{up to} (as opposed to: \emph{exactly})~$q$ colors.

Chromatic partitions of~$n$ have been considered in the combinatorics literature with various names, most frequently \emph{partitions of~$n$ objects of~$q$ colors}.
For~$q=2$, they appear to have been first considered by %
P.A.~MacMahon in~\cite[\S1.3, p.~483]{MMM1890}.
The fundamental representation by $\N^q_*$-multisets, however, appears to be novel.

Notably, the generating function of~$\card{\Part{q}{n}}$ is known, see OEIS~\cite{A217093},
\[
\card{\Part{q}{n}} = [x^n] \paren{\prod_{i=1}^\infty (1-x^i)^{-\binom{i+q-1}{q-1}}} \fstop
\]

\begin{defs}[Shape]
The \emph{shape} of a $q$-chromatic partition~$A\vDash\mbfn$ is the partition%
\[
\shape(A)\vdash \length{\mbfn} \comma \qquad \shape(A)_i\eqdef \sum_{\mbfa\in \supp{A}\, :\,  \length{\mbfa}=i} A(\mbfa) \comma \qquad i\in [\length{\mbfn}]\fstop
\]
\end{defs}

Note that~$\Part{q}{n}$ is the family of all $\N^q_*$-multisets~$A$ with~$\shape(A)\vdash n$.

\smallskip

In the graphical representation of a $q$-chromatic partition described in the introduction, each vector $\mbfa=\seq{a_1,a_2,\dotsc, a_q}\in\N^q_*$ is but the \emph{color count} of the type~$\tiny\ydiagram[*(black)]{1}_{a_1}~\tiny\ydiagram[*(gray)]{1}_{a_2} \cdots~\tiny\ydiagram[*(white)]{1}_{a_q}$, and~$A(\mbfa)$ is the number of blocks of type (indexed by)~$\mbfa$.
(Since every block in a polychromatic partition contains elements of at least one color, $\mbfa$~cannot be the $\zero$~vector, thus~$\mbfa\in\N^q_*$ rather than~$\N_0^q$.)
The~$\shape$ of~$A\vDash\mbfn$ is the (standard) partition of~$\length{\mbfn}$ obtained by forgetting the coloring.

\begin{rem}\label{r:ColoredVsPolychromatic}
A \emph{$q$-multipartition} of~$n$ is a partition of~$n$ for which each block is assigned one of a fixed set of~$q$ colors.
Multipartitions are of interest in number theory (see e.g.~\cite{BecCaiCheDilGonSu23}), and in combinatorics (see e.g., the survey~\cite{And08}).
In our language, $q$-multipartitions form a subset of all $q$-chromatic partitions.
Precisely, a $q$-chromatic partition~$A$ of~$n$ is a $q$-multipartition of~$n$ if and only if each~$\mbfa\in\supp A$ has exactly one non-zero entry.
$q$-Chromatic partitions have a much higher complexity than $q$-multipartitions, already when~$q=2$, cf.\ Table~\ref{tab:Uno}.
\end{rem}

\begin{table}[tb!]
\begin{tabular}{r | c | c | c | c c c c c c c c c c c}
\multirow{2}{*}{$n$} & \multirow{2}{*}{$P_2(n)$} & \multirow{2}{*}{$\card{\Part{2}{n}}$} & $p(n)$ & \multicolumn{11}{c}{$\card{\mcA_{(n-k,k)}}$}
\\
& & & $k=0$ & $1$ & 2 & 3 & 4 & 5 & 6 & 7 & 8 & 9 & 10 & 11
\\
\hline
1& 2 & 2 & 1 & 1
\\
2& 5 & 6 & 2 & 2 & 2
\\
3& 10 & 14 & 3 & 4 & 4 & 3
\\
4& 20 & 33 & 5 & 7 & 9 & 7 & 5
\\
5& 36 & 70 & 7 & 12 & 16 & 16 & 12 & 7
\\
6& 65 & 149 & 11 & 19 & 29 & 31 & 29 & 19 & 11
\\
7& 110 & 298 & 15 & 30 & 47 & 57 & 57 & 47 & 30 & 15
\\
8& 185 & 591 & 22& 45& 77& 97& 109& 97& 77& 45& 22
\\
9& 300 & 1132 & 30 & 67& 118& 162& 189& 189& 162& 118& 67& 30
\\
10& 481 & 2139 & 42& 97& 181& 257& 323& 339& 323& 257& 181& 97& 42
\\
11& 752 & 3948 & 56 & 139 & 267& 401& 522& 589& 589& 522& 401& 267& 139& 56
\\
\end{tabular}
\medskip
\caption{$P_2(n)$ is the number of (integer) multipartitions of~$n$ into blocks of two types, OEIS~\cite{A000712}.
$p(n)$ is the number of (integer) partitions of~$n$, OEIS~\cite{A000041}.
$\card{\Part{2}{n}}$ is the number of all 2-chromatic partitions of~$n$, OEIS~\cite{A005380} obtained by summing over~$k=0,\dotsc, n$ (i.e., over rows) the number~$\card{\mcA_{(n-k,k)}}$ of partitions with coloring~$\mbfn=(n-k,k)$. Note that this includes partitions of only one color for either color, namely partitions in~$\mcA_{(n,0)}$ and~$\mcA_{(0,n)}$.
}
\label{tab:Uno}
\end{table}

\subsubsection{The Young diagram of a polychromatic partition}\label{sss:ColorCountDiagrams}
Polychromatic partitions have a natural total order, the \emph{reverse shortlex} order which we now describe.
We denote by~$<_\lex$ the standard \emph{lexicographic order} on~$\N^q_*$, induced by the order of~$[q]$.
Recall that the \emph{shortlex order} on~$[0,\infty)^q$ is the total order defined as
\begin{equation}\label{eq:ShortlexOrder}
\mbfs <_{\shortlex} \mbft \qquad \iff \qquad (\length{\mbfs}, s_1,\dotsc, s_q) <_\lex (\length{\mbft}, t_1,\dotsc, t_q) %
\fstop
\end{equation}
The shortlex order on~$\N^q_*$ induces {both a standard representation of~$A \in \Part{q}{n}$ and a total order on~$\Part{q}{n}$ itself as follows.

\begin{defs}[Enumeration]\label{d:Enumeration}
	The \emph{enumeration} of the blocks in~$A$ is the tuple $\seq{\mbfa_h}_{h\in[\size{A}]}$ of all~$\mbfa\in\supp A$, each listed with multiplicity~$A(\mbfa)$, presented decreasingly in the shortlex order.
\end{defs}}

\begin{defs}[Reverse shortlex order]\label{d:OrderMultiYoung}
For polychromatic partitions $A,B\in\Part{q}{n}$ we say that~$A$ precedes~$B$ in the \emph{reverse shortlex order} (write~$A <_{\rsl} B$) if there exists~$\mbfb\in \N^q_*$ such that~$A(\mbfb) < B(\mbfb)$, and~$A(\mbfa)=B(\mbfa)$ for all $\mbfa>_{\shortlex} \mbfb$.
\end{defs}

{To each $q$-chromatic partition~$A \in \Part{n}{q}$ with enumeration~$\seq{\mbfa_h}_{h\in[\size{A}]}$, we may associate a colored Young diagram matching its color count as follows: the~$h^\text{th}$ of its~$\size{A}$ rows has (in order)~
$a_{h,1}$ boxes of color~$1$, $a_{h,2}$ boxes of color~$2$, etc., see Figure~\ref{fig:LongList}.}

\begin{rem}[Polychromatic partitions vs.\ Young tableaux]\label{r:Tableaux}
We prefer the reverse shortlex order since we think of Young diagrams in the English notation. The (forward) shortlex order on polychromatic partitions would result in Young diagrams in the French notation.

It is apparent from Figure~\ref{fig:LongList} that replacing colors~$\tiny\ydiagram[*(black)]{1},~\tiny\ydiagram[*(gray)]{1},\dotsc,~\tiny\ydiagram[*(white)]{1}$ by $q$ consecutive numbers in increasing order does \emph{not} result in a semistandard (nor standard) Young tableau: whereas the tableaux so obtained have weakly increasing rows, labels in their columns are not ordered.
\end{rem}

\begin{figure}[tb!]
\begin{gather*}
\underset{\scriptstyle{\car_{(4,2)}}}{\tiny\ydiagram[*(gray)]{4}*[*(white)]{4+2}}
\quad
\underset{\scriptstyle{\car_{(4,1)}+\car_{(0,1)}}}{\tiny\ydiagram[*(gray)]{4}*[*(white)]{4+1,1}}
\quad
\underset{\scriptstyle{\car_{(3,2)}+\car_{(1,0)}}}{\tiny\ydiagram[*(gray)]{3,1}*[*(white)]{3+2}}
\quad
\underset{\scriptstyle{\car_{(4,0)}+\car_{(0,2)}}}{\tiny\ydiagram[*(gray)]{4}*[*(white)]{0,2}}
\quad
\underset{\scriptstyle{\car_{(3,1)}+\car_{(1,1)}}}{\tiny\ydiagram[*(gray)]{3,1}*[*(white)]{3+1,1+1}}
\quad
\underset{\scriptstyle{\car_{(2,2)}+\car_{(2,0)}}}{\tiny\ydiagram[*(gray)]{2,2}*[*(white)]{2+2}}
\quad
\underset{\scriptstyle{\car_{(4,0)}+2\cdot \car_{(0,1)}}}{\tiny\ydiagram[*(gray)]{4}*[*(white)]{0,1,1}}
\\[6pt]
\underset{\scriptstyle{\car_{(3,1)}+\car_{(1,0)}+\car_{(0,1)}}}{\tiny\ydiagram[*(gray)]{3,1}*[*(white)]{3+1,0,1}}
\quad
\underset{\scriptstyle{\car_{(2,2)}+2\cdot \car_{(1,0)}}}{\tiny\ydiagram[*(gray)]{2,1,1}*[*(white)]{2+2}}
\quad
\underset{\scriptstyle{\car_{(3,0)}+\car_{(1,2)}}}{\tiny\ydiagram[*(gray)]{3,1}*[*(white)]{0,1+2}}
\quad
\underset{\scriptstyle{2\cdot \car_{(2,1)}}}{\tiny\ydiagram[*(gray)]{2,2}*[*(white)]{2+1,2+1}}
\quad
\underset{\scriptstyle{\car_{(3,0)}+\car_{(0,2)}+\car_{(1,0)}}}{\tiny\ydiagram[*(gray)]{3,0,1}*[*(white)]{0,2}}
\\[6pt]
\underset{\scriptstyle{\car_{(3,0)}+\car_{(1,1)}+\car_{(0,1)}}}{\tiny\ydiagram[*(gray)]{3,1}*[*(white)]{0,1+1,1}}
\quad
\underset{\scriptstyle{\car_{(2,1)}+\car_{(2,0)}+\car_{(0,1)}}}{\tiny\ydiagram[*(gray)]{2,2}*[*(white)]{2+1,0,1}}
\quad
\underset{\scriptstyle{\car_{(1,2)}+\car_{(2,0)}+\car_{(1,0)}}}{\tiny\ydiagram[*(gray)]{1,2,1}*[*(white)]{1+2}}
\quad
\underset{\scriptstyle{\car_{(2,1)}+\car_{(1,1)}+\car_{(1,0)}}}{\tiny\ydiagram[*(gray)]{2,1,1}*[*(white)]{2+1,1+1}}
\\[6pt]
\underset{\scriptstyle{\car_{(3,0)}+\car_{(1,0)}+2\cdot \car_{(0,1)}}}{\tiny\ydiagram[*(gray)]{3,1}*[*(white)]{0,0,1,1}}
\quad
\underset{\scriptstyle{\car_{(2,1)}+2\cdot \car_{(1,0)}+\car_{(0,1)}}}{\tiny\ydiagram[*(gray)]{2,1,1}*[*(white)]{2+1,0,0,1}}
\quad
\underset{\scriptstyle{\car_{(1,2)}+3\cdot\car_{(1,0)}}}{\tiny\ydiagram[*(gray)]{1,1,1,1}*[*(white)]{1+2}}
\quad
\underset{\scriptstyle{2\cdot \car_{(2,0)}+\car_{(0,2)}}}{\tiny\ydiagram[*(gray)]{2,2}*[*(white)]{0,0,2}}
\quad
\underset{\scriptstyle{\car_{(2,0)}+2\cdot \car_{(1,1)}}}{\tiny\ydiagram[*(gray)]{2,1,1}*[*(white)]{0,1+1,1+1}}
\\[6pt]
\underset{\scriptstyle{2\cdot \car_{(2,0)}+2\cdot\car_{(0,1)}}}{\tiny\ydiagram[*(gray)]{2,2}*[*(white)]{0,0,1,1}}
\quad
\underset{\scriptstyle{\car_{(2,0)}+\car_{(0,2)}+2\cdot \car_{(1,0)}}}{\tiny\ydiagram[*(gray)]{2,0,1,1}*[*(white)]{0,2}}
\quad
\underset{\scriptstyle{\car_{(2,0)}+\car_{(1,1)}+\car_{(1,0)}+\car_{(0,1)}}}{\tiny\ydiagram[*(gray)]{2,1,1}*[*(white)]{0,1+1,0,1}}
\quad
\underset{\scriptstyle{2\cdot \car_{(1,1)}+2\cdot\car_{(1,0)}}}{\tiny\ydiagram[*(gray)]{1,1,1,1}*[*(white)]{1+1,1+1}}
\\[6pt]
\underset{\scriptstyle{\car_{(2,0)}+2\cdot\car_{(1,0)}+2\cdot\car_{(0,1)}}}{\tiny\ydiagram[*(gray)]{2,1,1}*[*(white)]{0,0,0,1,1}}
\quad
\underset{\scriptstyle{\car_{(0,2)}+4\cdot\car_{(1,0)}}}{\tiny\ydiagram[*(gray)]{0,1,1,1,1}*[*(white)]{2}}
\quad
\underset{\scriptstyle{\car_{(1,1)}+3\cdot\car_{(1,0)}+\car_{(0,1)}}}{\tiny\ydiagram[*(gray)]{1,1,1,1}*[*(white)]{1+1,0,0,0,1}}
\quad
\underset{\scriptstyle{4\cdot\car_{(1,0)}+2\cdot\car_{(0,1)}}}{\tiny\ydiagram[*(white)]{0,0,0,0,1,1}*[*(gray)]{1,1,1,1}}
\end{gather*}
\caption{The color-count Young diagrams displaying all block types of $(4,2)$-chromatic partitions of~$6$ described in~\S\ref{sss:ColorCountDiagrams} and the associated polychromatic partitions (multisets) in Definition~\ref{d:MultiYoung}, first grouped by shape and subsequently arranged in lexicographic order. 
(Note that this is different from the shortlex order in Definition~\ref{d:OrderMultiYoung}.)
}
\label{fig:LongList}
\end{figure}

\subsubsection{Coloring permutations}
Each permutation~$\pi\in\mfS_n$ induces a \emph{set} partition~$\Lambda_\pi$ of~$[n]$ consisting of $\pi$'s orbits.
When~$\pi$ is presented in the cycle representation~\eqref{eq:StdRepPi}, we have
\begin{equation}\label{eq:StdSetPartitionPi}
\Lambda_\pi = \tset{\set{y_{1,1}, y_{1,2}, \dotsc}, \set{y_{2,1}, y_{2,2}, \dotsc}, \dotsc, \set{y_{r,1}, y_{r,2}, \dotsc}} \fstop
\end{equation}
The corresponding (integer) partition~$\boldlambda_\pi$ of~$n$ is obtained by counting the size of blocks in~$\Lambda_\pi$, and~$\pi$ is said to have \emph{cycle structure~$\boldlambda_\pi$}.
We denote the corresponding assignment by
\[
\struct\colon \pi \longmapsto \boldlambda_\pi \fstop
\]

In a similar fashion, in the presence of a $q$-coloring~$\coloring_\mbfn$ of~$[n]$, we may naturally assign to each permutation~$\pi$ acting on~$[n]$ a $q$-chromatic partition~$A_\pi$, defined as follows.
For every~$\mbfa \in\N^q_*$, we let~$A_\pi(\mbfa)$ be the number of blocks in~$\Lambda_\pi$ with color count~$\mbfa$.
In this case we say that~\emph{$\pi$ has coloring~$A_\pi$} and we define an assignment
\begin{equation}\label{eq:Pi}
\colore=\colore_{\coloring_\mbfn} \colon \mfS_{\length{\mbfn}} \longrar \Part{q}{n} \comma \qquad \colore\colon \pi\longmapsto A_\pi \fstop
\end{equation}

Clearly, if~$\pi$ has coloring~$A_\pi$, then it has cycle structure~$\shape(A_\pi)$. Formally,
\begin{equation}\label{eq:ShapeColorStruct}
\shape\circ \colore=\struct \fstop
\end{equation}

\subsection{Pattern inventories}\label{ss:Patterns}
The \emph{cycle index polynomial}~$Z^G$ of a %
group~$G< \mfS_n$ of degree~$n$~is
\[
Z^G(\mbft)\eqdef \frac{1}{\card{G}}\sum_{\pi\in G}\mbft^{\boldlambda_\pi} \comma \qquad \mbft=\seq{t_1,\dotsc, t_n} \fstop
\]
We denote by~$Z_n\eqdef Z^{\mfS_n}$ the cycle index polynomial of~$\mfS_n$. 
It is not difficult to show that%
\begin{equation}\label{eq:CycleIndex}
Z_n(\mbft)=\frac{1}{n!} \sum_{\boldlambda\vdash n} \MultiSecond(\boldlambda)\, \mbft^{\boldlambda}\comma \qquad \mbft=\seq{t_1,\dotsc, t_n} \fstop
\end{equation}
The powers of each $\mbft$-monomial correspond to the partition~$\boldlambda_\pi$ encoding the cycle structure of~$\pi$.
That is, $Z_n$~is the generating function for the cycle structure of permutations on~$[n]$.

\paragraph{Patterns}
Similarly to the case of~$Z_n$, we aim at defining a generating function~$\msZ^\sym{q}_\mbfn$ for the coloring of permutations of~$[\length{\mbfn}]$. %
We first need to define an appropriate set of monomials, reflecting the different roles of cycle structure and coloring.
To this end, fix~$k\in\N_1$ and let~$\mbfs_1,\dotsc,\mbfs_q, \boldalpha\in\R^k$ and~$\mbfS\eqdef \seq{\mbfs_1,\dotsc,\mbfs_q}\in\R^{k\times q}$.
Finally, for every~$\mbfa\leq_\had\mbfn$ set
\begin{align}\label{eq:YDefinition}
\begin{aligned}
\omega_\mbfa[\mbfS;\boldalpha]\eqdef \tparen{\mbfs_1^{\had a_1}\had\cdots \had \mbfs_q^{\had a_q}}\cdot\boldalpha \comma 
\qquad
\Omega_\mbfn[\mbfS;\boldalpha]\eqdef \tseq{\omega_\mbfa[\mbfS;\boldalpha]}_{\mbfa\leq_\had\mbfn}
 \fstop
\end{aligned}
\end{align}

We call $\mbfn$-\emph{pattern} of a polychromatic partition~$A$ the quantity
\begin{equation}\label{eq:Pattern}
w[\mbfS;\boldalpha](A) \eqdef \prod_{\mbfa\in \supp A} \omega_\mbfa[\mbfS;\boldalpha]^{A(\mbfa)} \fstop
\end{equation}

The following definition of the polynomial~$\msZ^\sym{q}_\mbfn$ is inspired by P\'olya Enumeration Theory.
\begin{defs}%
The \emph{pattern inventory of~$\mssC_\mbfn$} is the polynomial
\begin{equation}\label{eq:PatternInventory}
\msZ^\sym{q}_\mbfn[\mbfS;\boldalpha]=\msZ^\sym{q}_\mbfn[\mbfs_1,\dotsc,\mbfs_q;\boldalpha]\eqdef \frac{1}{\length{\mbfn}!}\sum_{\pi\in\mfS_{\length{\mbfn}}} w[\mbfs_1,\dotsc,\mbfs_q;\boldalpha]\tparen{\colore(\pi)} \fstop
\end{equation}
\end{defs}

From a structural viewpoint, the expression for~$\msZ^\sym{q}_\mbfn$ in~\eqref{eq:PatternInventory} is unsatisfactory, since many monomials in the summation coincide, as it is the case when~$q=1$; cf.~\eqref{eq:CycleIndex}.
Furthermore, the same expression is also unfeasible from a computational viewpoint due to the growth of the number~$\length{\mbfn}!$ of the summands.%

In order to simplify the expression~\eqref{eq:PatternInventory} of~$\msZ^\sym{q}_\mbfn$, we need to collect all identical monomials, and thus count how many such monomials appear.
Similarly to the case of~$q=1$ (when they are indexed by partitions), in the general case~$q>1$ they are indexed by $q$-chromatic partitions.

\begin{defs}%
\label{d:MultinomialCoeff}
The \emph{multinomial coefficient} of a $q$-chromatic partition~$A$ of~$\mbfn$ is
\begin{align}\label{eq:MultinomialMultiYoung}
\Multi{\mbfn}(A)\eqdef&\ \mbfn! \prod_{\mbfa\in \supp{A}} \frac{\binom{\length{\mbfa}}{\mbfa}^{A(\mbfa)}}{\length{\mbfa}^{A(\mbfa)} A(\mbfa)!}\fstop
\end{align}
\end{defs}

\begin{lem}\label{l:MultiA}
The number of permutations of~$[\length{\mbfn}]$ with coloring~$A\vDash \mbfn$ is~$\Multi{\mbfn}(A)$.
\begin{proof}
Fix~$A\vDash\mbfn$ {with enumeration~$(\mbfa_h)_{h \in [\size{A}]}$} and note that a permutation{~$\pi$} of~$[\length{\mbfn}]$ is a set partition{~$\Lambda_\pi$}
of~$[\length{\mbfn}]$ together with a cyclic order on each of its blocks.
Thus~$\colore^{-1}(A)$ is in bijection with the set of pairs consisting of a set partition
of~$[\length{\mbfn}]$ with exactly~$A(\mbfa)$ blocks of color count~$\mbfa$ for every~$\mbfa\in\N^q_*$,
and of a cyclic order on each of its blocks.

{First, we observe that each set partition as above is induced by the fibers of exactly~$\prod_\mbfa A(\mbfa)!$ many maps~$f \colon [\length{\mbfn}] \to [\size{A}]$ such that~$\card{\tparen{f^{-1}(h) \cap \mssC_\mbfn^{-1}(j)}} = a_{h,j}$ for every~$h \in [\size{A}]$ and~$j \in [q]$. Second, all such maps~$f$ are found by distributing, for each color~$j \in [q]$, the~$n_j$ elements of~$\mssC_\mbfn^{-1}(j)$ among the ordered
blocks according to their prescribed color counts; hence, their total number is}
\[
\prod_{j=1}^q \frac{n_j!}{\prod_{h=1}^{\size{A}} a_{h,j}!}
= \mbfn! \prod_{\mbfa\in\supp A} \frac{1}{\mbfa!^{A(\mbfa)}} \fstop
\]
Finally, each block with color count~$\mbfa$ carries exactly~$\paren{\length{\mbfa}-1}!$ cyclic orders.
{Combining these three observations} and using~$\paren{\length{\mbfa}-1}!\big/\mbfa!=
\tfrac{1}{\length{\mbfa}}\binom{\length{\mbfa}}{\mbfa}$ proves the assertion.
\end{proof}
\end{lem}

\begin{rem}[$q=1$]\label{r:UnidimensionalFirst}
When~$q=1$ (hence~$\mbfn=n$), %
a $1$-chromatic partition is identical to its $\shape$.
Thus, the multinomial coefficient of~$A\in\mcA_\mbfn$ reduces to the one~\eqref{eq:MultiSecond} of the integer partition~$\boldlambda\eqdef \shape(A)\vdash n$, viz.
\begin{align*}
\Multi{\mbfn}(A)= n! \prod_{i\in \supp{A}}\frac{1}{i^{\lambda_i} \lambda_i!}\comma \qquad \lambda_i\eqdef A(i) \fstop
\end{align*}
\end{rem}

\begin{cor}%
\label{c:Multi}
For every~$\mbfn\in\N^q_*$ we have
\[
\sum_{\substack{A \vDash \mbfn\\ \shape(A)=\boldlambda}} \Multi{\mbfn}(A)= \MultiSecond(\boldlambda) \fstop
\]
\begin{proof}
In light of~\eqref{eq:ShapeColorStruct},~$\struct^{-1}(\boldlambda)$ coincides with the disjoint union over~$A \vDash \mbfn$ with~$\shape(A)=\boldlambda$ of the preimages~$\colore^{-1}(A)$.
Since the union is disjoint, Lemma~\ref{l:MultiA} and~$\card{\struct^{-1}(\boldlambda)}=\MultiSecond(\boldlambda)$ imply the conclusion.
\end{proof}
\end{cor}

We now define a polychromatic cycle index polynomial, similarly to~\eqref{eq:CycleIndex}.

\begin{defs}
The $q$-\emph{chromatic cycle index polynomial of~$\mbfn$} is the polynomial
\begin{align}\label{eq:ZDefinition}
\begin{aligned}
Z^\sym{q}_\mbfn(\mbft)\eqdef \frac{1}{\length{\mbfn}!} \sum_{A\vDash \mbfn} \Multi{\mbfn}(A)\, \prod_{\mbfa\in \supp{A}} t_\mbfa^{A(\mbfa)} \comma \qquad 
Z^\sym{q}_{\zero}(\mbft)\eqdef 1 \comma \qquad \mbft\eqdef \seq{t_\mbfa}_{\mbfa\leq_\had\mbfn}
\fstop
\end{aligned}
\end{align}
\end{defs}

\begin{ese}
The $2$-chromatic cycle index polynomial of~$(4,2)$ may be computed from Figure~\ref{fig:LongList}.
Omitting parentheses and commas in the vector notation for brevity, and listing monomials corresponding to $2$-chromatic partitions of~$(4,2)$ as in Figure~\ref{fig:LongList}, we have
\begin{align*}
Z^\sym{2}_{42}(\mbft) &= \tfrac{1}{6!}\big(
120\, t_{42}+ 48\, t_{41}t_{01}+ 96\, t_{32}t_{10}+ 6\, t_{40}t_{02}+ 48\, t_{31}t_{11}+ 36\, t_{22}t_{20}+ 6\, t_{40}t_{01}^2
\\
&\qquad + 48\, t_{31}t_{10}t_{01}+ 36\, t_{22}t_{10}^2+ 16\, t_{30}t_{12}+ 24\, t_{21}^2+ 8\, t_{30}t_{02}t_{10}
\\
&\qquad + 16\, t_{30}t_{11}t_{01}+ 24\, t_{21}t_{20}t_{01}+ 24\, t_{12}t_{20}t_{10}+ 48\, t_{21}t_{11}t_{10}
\\
&\qquad + 8\, t_{30}t_{10}t_{01}^2+ 24\, t_{21}t_{10}^2t_{01}+ 8\, t_{12}t_{10}^3+ 3\, t_{20}^2t_{02}+ 12\, t_{20}t_{11}^2
\\
&\qquad + 3\, t_{20}^2 t_{01}^2+ 6\, t_{20}t_{02}t_{10}^2+ 24\, t_{20}t_{11}t_{10}t_{01}+ 12\, t_{11}^2t_{10}^2
\\
&\qquad + 6\, t_{20}t_{10}^2 t_{01}^2 + t_{02}t_{10}^4+ 8\, t_{11}t_{10}^3t_{01}+ t_{10}^4t_{01}^2
\big) \fstop
\end{align*}
\end{ese}

\begin{prop}\label{p:Z=Z}
Choosing~$\Omega_\mbfn[\mbfS;\boldalpha]$ as in~\eqref{eq:YDefinition} we have
\begin{align}\label{eq:t:Zn:0} \begin{aligned}
\msZ^\sym{q}_\mbfn[\mbfS;\boldalpha]%
=Z^\sym{q}_\mbfn\tparen{\Omega_\mbfn[\mbfS;\boldalpha]} %
=\frac{1}{\length{\mbfn}!} \sum_{A\vDash \mbfn} \Multi{\mbfn}(A)\, \prod_{\mbfa\in \supp{A}} \tparen{ \tparen{\mbfs_1^{\had a_1}\had\cdots \had \mbfs_q^{\had a_q}}\cdot\boldalpha}^{A(\mbfa)} \fstop
\end{aligned}
\end{align}

\begin{proof}
It suffices to collect all terms in~$\msZ^\sym{q}_\mbfn$ with the same monomials. 
By Lemma~\ref{l:MultiA}, for each~$\pi\in\mfS_{\length{\mbfn}}$ there are~$\Multi{\mbfn}(\colore(\pi))$ permutations~$\pi'$ with~$\colore(\pi')=\colore(\pi)\defeq A$.
The conclusion follows using that~$w_\mbfn[\mbfS;\boldalpha](\colore(\pi))= \prod_{\mbfa\in\supp{A}}\omega_\mbfa[\mbfS;\boldalpha]^{A(\mbfa)}$.
\end{proof}
\end{prop}

\section{A polychromatic Ewens Sampling Formula}\label{s:Polychromatic}
We denote by~$r=r_\pi$ the total number of cycles (including fixed points) of a permutation~$\pi\in\mfS_n$.
Let~$\theta>0$ and recall that a probability distribution on~$\mfS_n$ is \emph{$\theta$-biased} if its value on each~$\pi$ is proportional to~$\theta^r$. %
The \emph{Ewens Sampling Formula} (ESF) with parameter~$\theta$ is the probability distribution
\[
\Esf{n}{\theta}[\boldlambda]\eqdef \frac{n!}{\Poch{\theta}{n}} \prod_{i=1}^n \frac{\theta^{\lambda_i}}{i^{\lambda_i}\lambda_i!} = \frac{\theta^{\length{\boldlambda}}}{\Poch{\theta}{n}} \MultiSecond(\boldlambda) \comma \qquad \boldlambda\eqdef\seq{\lambda_1,\dotsc, \lambda_n} \vdash n \fstop
\]
(Recall that~$\Poch{\theta}{n}\eqdef \theta(\theta+1)\cdots (\theta+n-1)$ is the Pochhammer symbol.)
It gives the probability that a $\theta$-biased permutation has cycle structure~$\boldlambda$.
In particular, the distribution~$\Esf{n}{1}$ describes the frequency of a permutation in~$\mfS_n$ with a given cycle structure.

We refer to the recent surveys~\cite{Cra16,Tav21} and references therein for a complete account of the history and importance of the ESF throughout mathematics and in theoretical biology.

\paragraph{Relations between the ESF and the Dirichlet distribution}
For~$\boldalpha\in (0,\infty)^k$, the \emph{Dirichlet distribution}~$\Dir{\boldalpha}$ is the probability measure with density 
\begin{equation*}\label{eq:DirichletDistribution}
\frac{\car_{\Delta^{k-1}}(\mbfx)}{\Beta(\boldalpha)}\, \mbfx^{\boldalpha-\uno}
\end{equation*}
w.r.t.\ the standard Lebesgue measure on the hyperplane of equation~$\length{\mbfx}=1$.
When~$k=1$ we have~$\Delta^0= \set{1}$ and~$\Dir{\boldalpha}=\delta_1$ for every~$\boldalpha=\alpha\in (0,\infty)$.
Univariate moments of the Dirichlet distribution are the generating function of the standard ESF.
Indeed, we can express the generating function of~$\Esf{n}{\theta}$ %
in terms of the cycle index polynomial~$Z_n$, viz.
\begin{equation}\label{eq:DirichletEwens}
\sum_{\boldlambda \vdash n} \Esf{n}{\theta}[\boldlambda] \, \mbft^{\boldlambda} = \frac{n!}{\Poch{\theta}{n}}  Z_n(\theta\, \mbft)\comma %
\end{equation}
which is an expression for certain moments of~$\Dir{\boldalpha}$, as shown in~\cite{LzDS19a}.
We refer to~\S\ref{s:ApplMoments} for a thorough account; here, we are interested in establishing the polychromatic analogue of~\eqref{eq:DirichletEwens}.

\subsection{Polychromatic ESF}\label{ss:PESF}
Let us now turn to a polychromatic generalization of the ESF. 
Throughout this section, let~$n\in\N_0$,~$q\in\N_1$,~$\theta>0$.
For a Borel probability~$\nu$ on~$\Delta^{q-1}$, we set
\[
\mom{\mbfa}{\nu} \eqdef \int_{\Delta^{q-1}} \mbfu^\mbfa\diff\nu(\mbfu) %
\qquad \text{and}\qquad
\momentum{\mbfa}{\theta,\nu}\eqdef \frac{\theta}{\length{\mbfa}} \binom{\length{\mbfa}}{\mbfa} \mom{\mbfa}{\nu} \comma \qquad \mbfa \in \N^q_* \fstop
\]

\begin{defs}[Polychromatic ESF]\label{d:PolyESF}
We call \emph{polychromatic ESF} (in short: \emph{pESF}) with parameters~$\theta>0$ and~$\nu\in\msP(\Delta^{q-1})$ the function
\begin{equation}\label{d:eq:PolyESF}
\pEsf{n}{\theta,\nu}[A] %
\eqdef%
\frac{\theta^{\size{A}}}{\Poch{\theta}{n}}\, \frac{n!\Multi{\col{A}}(A)}{\col{A}!} \prod_{\mbfa\in\supp A} \mom{\mbfa}{\nu}^{A(\mbfa)}
= \frac{n!}{\Poch{\theta}{n}} \prod_{\mbfa\in \N^q_*} \frac{\momentum{\mbfa}{\theta,\nu}^{A(\mbfa)}}{A(\mbfa)!} 
\comma
\qquad A\in \Part{q}{n} \fstop
\end{equation}
\end{defs}

The parameter~$\theta>0$ has the same interpretation as in the ESF, while the parameter measure~$\nu$ determines the frequency of colors.

\begin{rem}[$q=1$]\label{r:UnidimensionalSecond}
When~$q=1$ %
formula~\eqref{d:eq:PolyESF} reduces to the ESF by Remark~\ref{r:UnidimensionalFirst}.
\end{rem}

\begin{rem}[Non-uniqueness]
For fixed~$n$, the association~$%
\nu\mapsto \pEsf{n}{\theta,\nu}$ is \emph{not} injective.
Indeed,~$\pEsf{n}{\theta,\nu}$ only depends on the $\mbfa$-moments~$\mom{\mbfa}{\nu}$ of~$\nu$ of order~$\length{\mbfa}\leq n$.
However, it will follow from Theorem~\ref{t:CharacterizationNew} below that the association~$%
(\theta,\nu)\mapsto \tseq{\pEsf{n}{\theta, \nu}}_n$ is injective.
\end{rem}

\begin{prop}\label{p:ESFProbab}
$\pEsf{n}{\theta,\nu}$ is a probability distribution on~$\Part{q}{n}$ and we have 
\begin{equation}\label{eq:p:Aggregation:0.2}
\shape_\pfwd \pEsf{n}{\theta,\nu} = \Esf{n}{\theta} \fstop
\end{equation}

\begin{proof}
In order to prove~\eqref{eq:p:Aggregation:0.2} it suffices to note that, by the Multinomial Theorem,
\[
\sum_{\mbfa\in \N^q_* : \length{\mbfa}=i} \momentum{\mbfa}{\theta,\nu} = \frac{\theta}{i} \int_{\Delta^{q-1}} \length{\mbfu}^i\diff\nu(\mbfu) = \frac{\theta}{i}\comma \qquad i \in [n]\fstop
\]
Raising both sides to the power~$\lambda_i$ and applying again the Multinomial Theorem yields~\eqref{eq:p:Aggregation:0.2}.
As a consequence, the total mass of~$\pEsf{n}{\theta,\nu}$ is equal to the total mass of~$\Esf{n}{\theta}$.
\end{proof}
\end{prop}

Similarly to~\eqref{eq:DirichletEwens}, the pESF may be encoded by a generating function; in the polychromatic case the dummy variables are indexed by block~\emph{type} rather than by block~\emph{size}.
Fix~$\mbfn\in\N^q_0$, set~$n\eqdef\length{\mbfn}$, and let~$\mbft\eqdef \seq{t_\mbfa}_{\mbfa\leq_\had\mbfn}$ and~$\momNA{\nu}\eqdef\seq{\mom{\mbfa}{\nu}}_{\mbfa\leq_\had\mbfn}$, so that~$\momNA{\nu}\had\mbft=\seq{\mom{\mbfa}{\nu} t_\mbfa}_{\mbfa\leq_\had\mbfn}$.
Further set, in analogy with~$\mbft^{\boldlambda}$,
\[
\mbft^A\eqdef \prod_{\mbfa\in\supp A} t_\mbfa^{A(\mbfa)}\comma \qquad A\in\Part{q}{n}\fstop
\]
Virtually every quantity related to~$\pEsf{n}{\theta,\nu}$ may be computed from its generating function. 
{The following proposition readily follows from the definitions of~$\pEsf{n}{\theta,\nu}$ and~$Z^\sym{q}_\mbfn$.}

\begin{prop}[Generating function of the pESF]\label{p:DirichletEwensPoly}
It holds that%
\begin{equation}\label{eq:DirichletEwensPoly:1}
\sum_{A\vDash \mbfn} \pEsf{n}{\theta,\nu}[A]\, \mbft^A= \frac{n!}{\Poch{\theta}{n}} \binom{n}{\mbfn}\, Z^\sym{q}_\mbfn\tparen{\theta\, \mbfm_\nu\had\mbft} \comma \qquad \mbfn\in\N^q_0\comma \quad \length{\mbfn}=n\comma
\end{equation}
and therefore summing over the colorings~$\mbfn$ gives the generating function.

\end{prop}

By a simple computation, together with the explicit expression~\eqref{eq:ZDefinition} for~$Z^\sym{q}_\mbfn$, we have

\begin{cor}[Asymptotics]\label{c:AsymptoticsPESF}
For every~$\mbfn$, for every~$A\vDash \mbfn$,
\begin{align*}
\pEsf{n}{0,\nu}[A]\eqdef&  \lim_{\theta\downarrow 0} \pEsf{n}{\theta,\nu}[A] = \begin{cases} \binom{n}{\mbfa} \mom{\mbfa}{\nu} & \text{if } A=\car_\mbfa \comma \length{\mbfa}=n \\ 0 & \text{if } \size{A}\geq 2\end{cases} \comma
\\
\pEsf{n}{\infty,\nu}[A] \eqdef& \lim_{\theta\uparrow \infty} \pEsf{n}{\theta,\nu}[A] = \begin{cases} \binom{n}{\mbfm} \Bary{\nu}^\mbfm & \text{if } A=m_1\car_{\mbfe_1}+ \cdots + m_q\car_{\mbfe_q}\comma \length{\mbfm}=n \\ 0 & \text{otherwise} \end{cases}\comma
\end{align*}
where~$\Bary{\nu}\eqdef \seq{\mom{\mbfe_1}{\nu},\dotsc,\mom{\mbfe_q}{\nu}}\in\Delta^{q-1}$ is the barycenter of~$\nu$.
\end{cor}

Note that the assignment~$\nu\mapsto\pEsf{\bullet}{\infty,\nu}$ is \emph{not} injective, since~$\pEsf{\bullet}{\infty,\nu}$ depends on~$\nu$ only via~$\Bary{\nu}$.

We denote by~$\mathring{\Delta}^{q-1}$ the relative interior of~$\Delta^{q-1}$, and, for each~$j\in [q]$, we denote by~$\Delta^{q-1}_j\eqdef \Delta^{q-1}\cap\set{\mbfx\in \R^q: x_j=0}$ the $j^\text{th}$ face of~$\Delta^{q-1}$.

\begin{prop}[Support properties]\label{p:SupportProp}
The following assertions hold:
\begin{enumerate}[$(i)$]
\item\label{i:p:SupportProp:1} For every~$A\in\Part{q}{n}$,
\[
\pEsf{n}{\theta,\nu}[A]>0 \qquad \iff \qquad \nu\bracket{\bigcap_{j: a_j>0} \tparen{\Delta^{q-1}\setminus \Delta^{q-1}_j} }>0 \quad \text{for every~$\mbfa\in \supp A$.}
\]

\item\label{i:p:SupportProp:2} In particular, if~$n\geq q$, then $\pEsf{n}{\theta,\nu}[A]>0$ for every~$A\in\Part{q}{n}$ if and only if $\nu \mathring{\Delta}^{q-1}>0$.

\item\label{i:p:SupportProp:3} For every~$n\geq 2$, for every~$B\in\Part{q}{n-1}$, if~$\pEsf{n-1}{\theta,\nu}[B]=0$, then~$\pEsf{n}{\theta,\nu}[B+\car_{\mbfe_j}]=0$ for every~$j\in [q]$ and~$\pEsf{n}{\theta,\nu}[B+\car_{\mbfb+\mbfe_j}-\car_\mbfb]=0$ for every~$j\in [q]$ and every~$\mbfb\in\supp B$.
\end{enumerate}

\begin{proof}
It is clear that
\begin{equation}\label{eq:p:SupportProp:1}
\pEsf{n}{\theta,\nu}[A]=0\qquad \iff \qquad \mom{\mbfa}{\nu}=0 \quad \text{for some} \quad \mbfa\in\supp A\fstop
\end{equation}
For each~$\mbfa\in \N_0^q$ it is readily verified that
\begin{equation}\label{eq:p:SupportProp:2}
\mom{\mbfa}{\nu} = 0 \quad \iff \quad (\mbfu\mapsto \mbfu^\mbfa) \equiv 0 \;\;\;\; \nu\text{-a.e.} \quad \iff \quad \nu\bracket{\bigcup_{j: a_j>0}\Delta^{q-1}_j}= 1 \fstop
\end{equation}
By contrapositive, this proves~\ref{i:p:SupportProp:1}.
\ref{i:p:SupportProp:2} follows from~\ref{i:p:SupportProp:1} by choosing~$A\geq \car_\mbfa\in\Part{q}{n}$ for some~$\mbfa$ with~$a_j>0$ for all~$j\in [q]$, which exists since~$n\geq q$.
To show~\ref{i:p:SupportProp:3} choose~$B\in \Part{q}{n-1}$ with~$\pEsf{n-1}{\theta,\nu}[B]=0$, and fix~$j\in [q]$ and~$\mbfb\in\supp B$. By~\eqref{eq:p:SupportProp:1} there exists~$\mbfa\in \supp B$ such that~$\mom{\mbfa}{\nu}=0$.
Since~$\mbfa\in \supp B \subset \supp (B+\car_{\mbfe_j})$, by the converse implication in~\eqref{eq:p:SupportProp:1} we have~$\pEsf{n}{\theta,\nu}[B+\car_{\mbfe_j}]=0$.
If~$\mbfb\neq \mbfa$, then the same reasoning shows that~$\pEsf{n}{\theta,\nu}[B+\car_{\mbfb+\mbfe_j}-\car_\mbfb]=0$.
If~$\mbfb= \mbfa$, the assertion follows from the converse implication in~\eqref{eq:p:SupportProp:1} if we show that~$\mom{\mbfb+\mbfe_j}{\nu}=0$.
The conclusion follows from the chain of implications in~\eqref{eq:p:SupportProp:2}, since
\[
\nu\bracket{\bigcup_{k: b_k+\delta_{jk}>0}\Delta^{q-1}_k} \geq \nu\bracket{\bigcup_{k: a_k>0}\Delta^{q-1}_k} = 1 \fstop
\qedhere
\]
\end{proof}
\end{prop}

\subsection{Polychromatic Kingman consistency}\label{ss:Consistency} 
In~\cite{Kin78b,Kin78}, J.F.C.~Kingman introduced a celebrated notion of \emph{consistency} for stochastic processes on partitions, and showed that a sequence of random partitions~$\seq{\boldlambda_n}_n$ with~$\boldlambda_n\vdash n$ distributed according to~$\Esf{n}{\theta}$, satisfies this notion.
Precisely, if~$n$ objects are partitioned into blocks of sizes~$\boldlambda_n$, and one object is deleted uniformly at random, independently of~$\boldlambda_n$, the~$n-1$ remaining objects are partitioned into blocks of sizes (distributed as)~$\boldlambda_{n-1}$, cf.\ e.g.~\cite[p.~146]{Pit95}.

In this section, we introduce a similar consistency property for polychromatic partitions.
Firstly, we define a projective system of kernels on polychromatic partitions, induced by deletion uniformly at random of an element of~$[n]$.
Recall~\eqref{eq:MultiSets} and, for each~$A\in \Part{q}{\cup}$, set
\begin{equation}\label{eq:AmenoEj}
A_{\setminus \mbfa,j}\eqdef \begin{cases} A-\car_\mbfa & \text{if } \mbfa=\mbfe_j 
\\
A-\car_\mbfa +\car_{\mbfa-\mbfe_j} & \text{otherwise}
\end{cases}\comma \qquad j\in [q] \comma  \mbfa\in \supp{A}\text{ with } a_j \ge 1 \fstop
\end{equation}
Note that~$(\mbfa,j) \mapsto A_{\setminus \mbfa,j}$ is injective for every fixed~$A$.
Following~\cite{Kin78}, for each~$q\in\N_1$ we define a system~$\BackwardK^\sym{q}=\BackwardK^\sym{q}_{m,n}$, $n,m\in \N_0$, $m\leq n$, of probability kernels on~$\Part{q}{\cup}$.
Firstly, set
\begin{subequations}\label{eq:Kernels}
\begin{align}
\BackwardK^\sym{q}_{n,n}(A,B)\eqdef&\ \car_{A=B}\comma & A,B&\in\Part{q}{n} \comma
\\
\label{eq:Kernels:B}
\BackwardK^\sym{q}_{n-1,n}(A,B)\eqdef&\ \begin{cases} 
\frac{a_j A(\mbfa)}{n} & \text{if } B=A_{\setminus \mbfa,j} 
\\
0 & \text{otherwise}
\end{cases} \comma & A\in\Part{q}{n}\comma& B\in\Part{q}{n-1}\comma \qquad n\geq 1\fstop
\end{align}
\end{subequations}
(Recall that~$\Part{q}{0}\eqdef\set{0}$.)
Note that~$\BackwardK^\sym{q}_{n-1,n}(A,\emparg)$ is a probability on~$\Part{q}{n-1}$ for every~$A\in\Part{q}{n}$.
Secondly, let~$\BackwardK^\sym{q}=\BackwardK^\sym{q}_{m,n}$ be the unique system of kernels extending~\eqref{eq:Kernels} and satisfying the cocycle relation
\begin{equation}\label{eq:Cocycle}
\BackwardK^\sym{q}_{\ell, n}(A,C)= \sum_{B\in\Part{q}{m}} \BackwardK^\sym{q}_{\ell,m}(B,C)\, \BackwardK^\sym{q}_{m,n}(A,B)\comma \qquad 
\begin{gathered}
A\in\Part{q}{n}\comma C\in\Part{q}{\ell}\comma 
\\
\quad 0\leq \ell < m < n\fstop\quad\,
\end{gathered}
\end{equation}
Further note that~$\BackwardK^\sym{q}_{m,n}(A,\emparg)$ is a probability on~$\Part{q}{m}$ for every~$m$ and every~$A\in\Part{q}{n}$, since it is so for~$m=n-1$ as noted above, and in light of~\eqref{eq:Cocycle}.

\begin{rem}\label{r:KernelInterpretation}
Analogously to the case of usual partitions,~$\BackwardK^\sym{q}_{m,n}(A,B)$ may be interpreted as the selection of a random sampling (uniform, without replacement) of $m$~elements from~$[n]$ when~$[n]$ is partitioned according to a given $q$-chromatic partition~$A\in\Part{q}{n}$, resulting in the $q$-chromatic partition~$B\in\Part{q}{m}$.
The cocycle relation~\eqref{eq:Cocycle} is then a consequence of the consistency of random sub-sampling.
\end{rem}

Let us now turn to probability measures on~$\Part{q}{n}$.
For~$n\in\N_1$ let~$\msP(\Part{q}{n})$ be the set of all probability measures on~$\Part{q}{n}$.
Define a system~$\sigma^\sym{q}$ of maps~$\sigma^\sym{q}_{m,n}\colon \msP(\Part{q}{n})\to\msP(\Part{q}{m})$ by
\begin{align}\label{eq:Consistency}
\tparen{\sigma^\sym{q}_{m,n} \probP_n}[B] = \sum_{A\in \Part{q}{n}} \probP_n[A]\, \BackwardK^\sym{q}_{m,n}(A,B) \comma \qquad B\in\Part{q}{m} \comma \qquad \probP_n \in \msP(\Part{q}{n})\comma
\end{align}
and note that~$\sigma^\sym{q}$ satisfies the cocycle relation
\begin{equation}\label{eq:CocycleSigma}
\sigma^\sym{q}_{\ell, n}=\sigma^\sym{q}_{\ell, m}\circ\sigma^\sym{q}_{m, n} \comma \qquad \ell<m<n\fstop
\end{equation}

\begin{defs}%
\label{d:Consistency}
A sequence~$\probP_\bullet\eqdef \seq{\probP_n}_{n\in\N_0}$, with~$\probP_n\in\msP(\Part{q}{n})$, is \emph{consistent} if%
\begin{equation}\label{eq:d:Consistency:0}
\probP_m=\sigma^\sym{q}_{m,n}\probP_n\comma \qquad m\leq n\fstop
\end{equation}
We call any consistent~$\probP_\bullet$ a ($q$-)\emph{chromatic partition structure} (in short:~$\CPS{q}$).
Note that, for every such~$\probP_\bullet$, we have~$\probP_0[0] = 1$.
We write~$\CPS{q}$ also for the set of all~$\CPS{q}$'s.
\end{defs}

Let us note explicitly how~\eqref{eq:d:Consistency:0} reads for a $\CPS{q}$ when~$m=n-1$. Letting~$\mbfa\eqdef\mbfb+\mbfe_j$ and~$A=B+\car_{\mbfb+\mbfe_j}-\car_\mbfb$ in~\eqref{eq:Kernels:B}, we have
\begin{equation}\label{eq:ConsistencySpecial}
\begin{aligned}
\probP_{n-1}[B] &= \sum_{j=1}^q \probP_n[B+\car_{\mbfe_j}] \frac{B(\mbfe_j)+1}{n}
\\
&\qquad + \sum_{j=1}^q \sum_{\mbfb\in\supp B} \probP_n [B-\car_\mbfb+\car_{\mbfb+\mbfe_j}] \frac{(b_j+1) \tparen{B(\mbfb+\mbfe_j)+1}}{n} \comma
\end{aligned}
\end{equation}
to be compared with the monochromatic case ($q=1$) in~\cite[Eqn.~(2.6)]{Kin78b}.

\begin{ese}
In Proposition~\ref{p:RepresentationPESF} below we show that~$\pEsf{\bullet}{\theta,\nu}\eqdef \tseq{\pEsf{n}{\theta,\nu}}_n$ is a $\CPS{q}$.
\end{ese}

\subsection{Characterization of the pESF}\label{ss:CharacterizationESF}  
In~\cite{Kin78b} Kingman proved that the ESF is the only partition structure satisfying a recurrence relation in the form, cf.~\cite[Eqn.~(3.1)]{Kin78b},
\begin{equation}\label{eq:KingmanDefinitionA}
r\lambda_r\, \probP_n[\lambda_1, \dotsc, \lambda_r, \dotsc] = b(n,r)\, \probP_{n-r}[\lambda_1,\dotsc \lambda_{r-1}, \lambda_r-1,\lambda_{r+1},\dotsc] \comma %
\end{equation}
for~$r\in [n]$ with unspecified constants~$b(n,r)>0$.
We discuss an analogous characterization for the pESF.
Let~$q\in\N_1$, $\theta>0$, and~$\nu\in \msP(\Delta^{q-1})$. 
Omitting~$\theta$ and~$\nu$ from the notation, set
\begin{equation}\label{eq:KingmanDefinitionB}
b_*(n,\mbfr)\eqdef \binom{\length{\mbfr}}{\mbfr}\, \theta\, \mom{\mbfr}{\nu} \frac{\Poch{\theta}{n-\length{\mbfr}}}{\Poch{\theta}{n}} \frac{n!}{(n-\length{\mbfr})!} \comma\qquad n\in\N_1 \comma \mbfr \in \N^q_* \, : \, \length{\mbfr} \le n\fstop
\end{equation}

\begin{prop}%
\label{p:pESFRecurrence}
The pESF~$\probP_n \coloneqq \pEsf{n}{\theta,\nu}$ satisfies
\begin{equation}\label{eq:Kingman(3.1)b_*}
\length{\mbfr}\, A(\mbfr) \, \probP_n[A] = b_*(n,\mbfr)\, \probP_{n-\length{\mbfr}}[A-\car_\mbfr] \comma 
\end{equation}
for every~$A\in\Part{q}{n}$, every~$\mbfr\in \N^q_*$ with~$\car_\mbfr \leq A$, with~$b_*(n,\mbfr)$ as in~\eqref{eq:KingmanDefinitionB}. 

\begin{proof}
	The proof is straightforward from the definition~\eqref{d:eq:PolyESF} of~$\pEsf{n}{\theta,\nu}$.
\end{proof}

\end{prop}

\begin{rem}[$q=1$]
When~$q=1$%
~\eqref{eq:Kingman(3.1)b_*} is indeed~\cite[Eqn.~(3.1)]{Kin78b} with~$a_r\eqdef A(\mbfr)$.
\end{rem}

For the standard ESF, it is shown by Kingman in~\cite[pp.~5f.]{Kin78b} that, whenever~$\probP_\bullet$ is a partition structure satisfying~\eqref{eq:KingmanDefinitionA} %
for every~$r,n$, then~$\probP_n=\Esf{n}{\theta}$ for some~$\theta>0$.
In order to state the polychromatic analogue of this characterization, we shall need the following definitions.

\begin{defs}[Moment support]
We call \emph{moment support} of~$\nu\in\msP(\Delta^{q-1})$ the set
\begin{equation}\label{eq:SuppMomentsNu}
\Sigma_\nu\eqdef \set{\mbfa\in \N_0^q : \mom{\mbfa}{\nu}>0} \fstop
\end{equation} 
\end{defs}

The next lemma is straightforward.
\begin{lem}\label{l:HoppeTransition}
\begin{enumerate*}[$(i)$]
\item\label{i:l:HoppeTransition:1} $\Sigma_\nu$ is (downward) order-closed: if~$\mbfb\in\Sigma_\nu$ and~$\mbfa\leq_\had \mbfb$, then~$\mbfa\in \Sigma_\nu$. In particular,~$\zero\in \Sigma_\nu$;
\item\label{i:l:HoppeTransition:2} $\sum_{j=1}^q \mom{\mbfb+\mbfe_j}{\nu}=\mom{\mbfb}{\nu}$ for every~$\mbfb\in \N_0^q$.
\end{enumerate*}
\end{lem}

\begin{defs}[Block support of a~$\CPS{q}$]
The \emph{block support} of~$\probP_\bullet\in\CPS{q}$ is the set
\[
\BlockSupp{\probP_\bullet} \eqdef \set{\mbfa\in \N_0^q : \probP_{\length{\mbfa}}[\car_\mbfa]>0}\fstop
\]
\end{defs}
Note that~$\zero\in\BlockSupp{\probP_\bullet}$ since~$\probP_0[0]=1$.

\begin{lem}\label{l:BlockSupp}
For every~$\probP_\bullet\in\CPS{q}$, the following hold:
\begin{enumerate*}[$(i)$]
\item\label{i:l:BlockSupp:1} $\BlockSupp{\probP_\bullet}$ is order-closed: if~$\mbfb \in \BlockSupp{\probP_\bullet}$ and~$\mbfa\leq_\had \mbfb$ then~$\mbfa\in\BlockSupp{\probP_\bullet}$;
\item\label{i:l:BlockSupp:2} if~$\probP_n[A]>0$ then~$\supp A\subset \BlockSupp{\probP_\bullet}$.
\end{enumerate*}

\begin{proof}
Consistency~\eqref{eq:ConsistencySpecial} with~$m=n-1$ gives
\begin{equation}\label{eq:l:BlockSupp:1}
\probP_{n-1}[A_{\setminus\mbfa,j}] \geq \probP_n[A] \frac{a_j\, A(\mbfa)}{n} \comma \qquad \mbfa\in\supp A: a_j\geq 1 \fstop
\end{equation}
Deleting elements from~$[n]$ until~$\car_\mbfa$ is reached and iteratively applying the estimate~\eqref{eq:l:BlockSupp:1} shows~\ref{i:l:BlockSupp:1}.
For~\ref{i:l:BlockSupp:2} we show that~$\probP_{\length{\mbfa}}[\car_\mbfa]>0$ for every~$\mbfa\in\supp A$ for every~$A\in\Part{q}{n}$ with~$\probP_n[A]>0$.
For any such~$A$ and~$\mbfa$ write~$A=\car_\mbfa+B$ with~$B\in\Part{q}{n-\length{\mbfa}}$.
Successively deleting elements belonging to blocks represented by~$B$ leaves the probability of the resulting multiset positive by~\eqref{eq:l:BlockSupp:1}.
Deleting elements until only~$\car_\mbfa$ remains proves the assertion.
\end{proof}
\end{lem}

We prove a characterization for the pESF analogous to Kingman's characterization in the monochromatic case.

\begin{thm}[Characterization]\label{t:CharacterizationNew}
Let~$\probP_\bullet\in \CPS{q}$ satisfy
\begin{enumerate}[$({{a}}_{\arabic*})$]
\item\label{i:t:CharacterizationNew:1} $\mbfe_j\in \BlockSupp{\probP_\bullet}$ for every~$j\in [q]$;
\item\label{i:t:CharacterizationNew:2} $\probP_n[A]>0$ for every~$n\in \N_1$ and every~$A\in\Part{q}{n}$ with~$\supp A\subset \BlockSupp{\probP_\bullet}$;
\item\label{i:t:CharacterizationNew:3} there exists~$\mbfc\in\BlockSupp{\probP_\bullet}$ with~$\length{\mbfc}=2$;
\end{enumerate}
\begin{enumerate}[$(a)$]
\setcounter{enumi}{1}
\item\label{i:t:CharacterizationNew:4} for every~$n\in\N_1$, every~$A\in\Part{q}{n}$, and every~$\mbfr\in \N^q_*$ with~$\car_\mbfr \leq A$, it holds that
\begin{equation}\label{eq:Kingman(3.1)}
\length{\mbfr}\, A(\mbfr) \, \probP_n[A] = b(n,\mbfr)\, \probP_{n-\length{\mbfr}}[A-\car_\mbfr] \comma 
\end{equation}
for some constant~$b(n,\mbfr)$ independent of~$A$. 
\end{enumerate}
Then,~$\probP_n=\pEsf{n}{\theta,\nu}$ for a unique~$\theta>0$ and a unique Borel probability measure~$\nu$ on~$\Delta^{q-1}$.
Furthermore,~$\BlockSupp{\probP_\bullet}=\Sigma_\nu$ as in~\eqref{eq:SuppMomentsNu} and~$b\equiv b_*$ as in~\eqref{eq:KingmanDefinitionB}.
\end{thm}

\begin{rem}[Necessity of the assumptions]
Note that:
\begin{enumerate*}[$(1)$]
\item Assumption~\ref{i:t:CharacterizationNew:2} is the converse to the one in Lemma~\ref{l:BlockSupp}\ref{i:l:BlockSupp:2}.
\item Assumptions~\ref{i:t:CharacterizationNew:1}-\ref{i:t:CharacterizationNew:3} are satisfied in particular if~$\probP_n[A]>0$ for every~$A\in\Part{q}{n}$ for every~$n\in\N_1$, in which case Proposition~\ref{p:SupportProp}\ref{i:p:SupportProp:2} forces~$\nu\mathring{\Delta}^{q-1}>0$.
\item If Assumption~\ref{i:t:CharacterizationNew:1} is violated for some~$j\in [q]$, then~$\probP_{\length{\mbfn}}[A]=0$ whenever~$A\vDash \mbfn$ with~$n_j>0$, as a consequence of Lemma~\ref{l:BlockSupp}; that is, if~$\probP_1[\car_{\mbfe_j}]=0$ for some~$j\in [q]$, then~$\probP_\bullet$ does not charge any multiset in which the color~$j$ appears.
This shows that assumption~\ref{i:t:CharacterizationNew:1} is necessary to avoid degenerate cases that can be obtained for lower~$q$.
\item We will show in Remark~\ref{r:NecessityAssChar} that also~\ref{i:t:CharacterizationNew:2} and~\ref{i:t:CharacterizationNew:3} are necessary, and that~\ref{i:t:CharacterizationNew:4} does not detect their failure.
In particular, the validity of~\ref{i:t:CharacterizationNew:2} excludes the limiting case~$\theta\to 0$ and the validity of~\ref{i:t:CharacterizationNew:3} excludes the limiting case~$\theta\to\infty$. 
\end{enumerate*}
\end{rem}

\subsubsection{A discrete random growth model}\label{sss:DRG}
Our proof of Theorem~\ref{t:CharacterizationNew} requires the classification of harmonic functions of a Markov chain describing a random growth on~$[q]$, which we now discuss.
We mostly adhere to the terminology and notation in~\cite[\S24, pp.~256-259]{Woe00}; also see~\cite[Ch.~7]{Rev84}. %
Throughout, fix~$\mbfp\in\Delta^{q-1}$ and set, conventionally,~$0^0=1$.

Let~$\seq{j_t}_{t\geq 1}$ be i.i.d.\ random variables with categorical distribution~$\Cat_\mbfp$ on~$[q]$.
We call \emph{discrete} \emph{$\mbfp$-random growth on~$[q]$} (in short:~$\DRG{\mbfp}$) the discrete-time Markov chain~$\mbfr_\proc\eqdef\seq{\mbfr_t}_{t\geq 0}$ with state space~$\N^q_0$ defined by~$\mbfr_0\eqdef \zero$ and~$\mbfr_t\eqdef\sum_{s=1}^t \mbfe_{j_s}$.
We may and shall assume with no loss of generality, up to reducing~$q$ as necessary, that~$\mbfp\in \mathring{\Delta}^{q-1}$.
Thus, $\DRG{\mbfp}$ is the $\mbfp$-biased random walk on the $q$-dimensional directed Pascal graph~$\N_0^q$.
A function~$h\colon \N^q_0 \to\R$ is ($\DRG{\mbfp}$-)\emph{harmonic} if
\begin{equation}\label{eq:d:Harmonic}
\sum_{k=1}^q h(\mbfr +\mbfe_k)\, p_k = h(\mbfr) \comma \qquad \mbfr\in\N^q_0 \fstop
\end{equation}
Note that $\mbfp^{-\mbfr}\mom{\mbfr}{\nu}$ satisfies this relation for any~$\nu\in\msP(\Delta^{q-1})$ by Lemma~\ref{l:HoppeTransition}\ref{i:l:HoppeTransition:2}.

We denote by~$\mcH^\tgeq_\mbfp$ the space of all non-negative harmonic functions. 
An~$h_0\in\mcH^\tgeq_\mbfp$ is \emph{minimal} if~$h_0(\zero)=1$ and, whenever~$h_0\geq h$ for some~$h\in\mcH^\tgeq_\mbfp$, then~$h=\tau h_0$ for some~$\tau\in [0,1]$.

\begin{prop}\label{p:RepresentationH}
The Martin boundary of~$\DRG{\mbfp}$ is (homeomorphic to)~$\Delta^{q-1}$ and the Martin compactification is (homeomorphic to) $\N_0^q\sqcup \Delta^{q-1}$.
Furthermore, the Martin kernel (based at~$\zero\in\N_0^q$) of~$\DRG{\mbfp}$ satisfies
\begin{align} \label{eq:KernelDRGp}
\MartinK(\mbfr,\mbfs) = \begin{cases}\displaystyle \frac{(\length{\mbfs}-\length{\mbfr})!}{\length{\mbfs}!} \frac{\mbfs!}{(\mbfs-\mbfr)!} \, \mbfp^{-\mbfr} & \text{if } \mbfs,\mbfr\in \N_0^q\comma \mbfs\geq_\had \mbfr \\ \mbfp^{-\mbfr}\mbfs^\mbfr & \text{if } \mbfs\in\Delta^{q-1} %
\\ 0 & \text{otherwise} \end{cases}\comma\qquad \mbfr\in\N_0^q \comma
\end{align}
and the extended kernel~$\MartinK(\emparg,\mbfs)\colon \mbfr \mapsto \mbfp^{-\mbfr}\mbfs^\mbfr$ is minimal harmonic for every~$\mbfs\in\Delta^{q-1}$.
Finally, for each~$h\in\mcH^\tgeq_\mbfp$ there exists a unique Borel measure~$\nu^h$ on~$\Delta^{q-1}$ such that
\begin{equation}\label{eq:l:MinimalBoundary:1}
h(\mbfr)=\int_{\Delta^{q-1}} \MartinK(\mbfr,\mbfu)\diff\nu^h(\mbfu) = \mbfp^{-\mbfr} \int_{\Delta^{q-1}} \mbfu^\mbfr \diff\nu^h(\mbfu) \comma \qquad \mbfr\in\N^q_0 \fstop
\end{equation}

\end{prop}

Noting that~$h$ is $\DRG{\mbfp}$-harmonic if and only if~$h'\colon \mbfr\mapsto q^{\length{\mbfr}}\mbfp^\mbfr h(\mbfr)$ is $\DRG{q^{-1}\car}$-harmonic, Proposition~\ref{p:RepresentationH} follows from the results in~\cite{GerGruHag16}.
In particular, the above characterization of the Martin boundary together with~\eqref{eq:l:MinimalBoundary:1} are noted (in the unbiased case~$\DRG{q^{-1}\car}$) in~\cite[Ex.~6]{GerGruHag16}.
We provide an independent proof for completeness, and record a proof of minimality, which does not follow from~\cite{GerGruHag16}.

\begin{proof}[Proof of Proposition~\ref{p:RepresentationH}]
	\paragraph{Martin boundary}
	The formula
	\[
		\MartinK(\mbfr,\mbfs) = \frac{(\length{\mbfs}-\length{\mbfr})!}{\length{\mbfs}!} \frac{\mbfs!}{(\mbfs-\mbfr)!} \, \mbfp^{-\mbfr} \comma \qquad \mbfs,\mbfr\in \N_0^q\comma \mbfs\geq_\had \mbfr
	\]
	is straightforward and gives
	\begin{equation}\label{eq:p:MartinBd:1}
		\abs{\MartinK(\mbfr,\mbfs) - \length{\mbfs}^{-\length{\mbfr}} \mbfs^\mbfr \mbfp^{-\mbfr} } \leq \mbfp^{-\mbfr}\paren{1-{\AntiPoch{\length{\mbfs}}{\length{\mbfr}}}/{\length{\mbfs}^{\length{\mbfr}}}} \comma \qquad \mbfs,\mbfr\in \N_0^q\comma \mbfs\geq_\had \mbfr \fstop
	\end{equation}
	Moreover,~$\MartinK(\mbfe_j,\mbfs)= p_j^{-1} s_j/\length{\mbfs}$ for every~$j\in [q]$ and every~$\mbfs \in \N_0^q$,~$\mbfs \geq_\had \mbfe_j$.
	Thus, a sequence~$\seq{\mbfs_n}_n$ converges to the Martin boundary if and only if~$\mbfs_n\to\infty$ (i.e., $\length{\mbfs_n}\to\infty$) and there exists the vector~$\mbfs\in\Delta^{q-1}$ of asymptotic frequencies~$\mbfs \eqdef \lim_{n\to\infty} \mbfs_n/\length{\mbfs_n}$, for which we have~$%
	\nlim \MartinK(\mbfr,\mbfs_n) = \mbfp^{-\mbfr} \mbfs^\mbfr  %
	$%
	, for every~$\mbfr\in\N_0^q$.
	Since~$\mbfs$ uniquely identifies each limit~$\mbfr\mapsto \mbfp^{-\mbfr} \mbfs^\mbfr = \lim_{n\to\infty} \MartinK(\mbfr,\mbfs_n)$ along a given~$\mbfs_n\to\infty$ with frequencies~$\mbfs$, the proof of the first assertion and of~\eqref{eq:KernelDRGp} is concluded.
	
	\paragraph{The extended kernel is harmonic}
	For~$\mbfs\in\Delta^{q-1}$, let~$h_\mbfs \coloneqq \MartinK(\emparg,\mbfs)$. We have
	\[
	\sum_{j=1}^q h_\mbfs(\mbfr+\mbfe_j) p_j = \sum_{j=1}^q \mbfs^{\mbfr} \mbfp^{-\mbfr} s_j p_j^{-1} p_j = h_\mbfs(\mbfr) \sum_{j=1}^q s_j = h_\mbfs(\mbfr)\comma \qquad \mbfr\in\N^q_0 \comma
	\]
	that is,~$h_\mbfs$ is harmonic. Furthermore, it is clear that~$h_\mbfs(\zero)=1$.
	
	\paragraph{Representing measure} Let~$h\in\mcH^\tgeq_\mbfp$.
	Using that~$w\colon \mbfr \mapsto \mbfp^\mbfr h(\mbfr)$ satisfies the assumptions in~\cite[Thm~2.1]{HelSto03}, we deduce the existence of~$\nu^h$ as in~\eqref{eq:l:MinimalBoundary:1}, as the unique solution of the \emph{generalized Dale moment problem} on the simplex. (Note that~\cite{HelSto03} uses the corner simplex: the equivalent assertion we need here follows by the standard change of variables from the corner simplex to the standard simplex.)
	Since~$\Delta^{q-1}\subset \R^q$ is compact, the multi-dimensional Hausdorff moment problem on~$\Delta^{q-1}$ is well-posed, which uniquely fixes~$\nu^h$; also cf.~\cite[p.~264]{HelSto03}.
	
	\paragraph{Minimality}
	Let~$h\in\mcH^\tgeq_\mbfp$ with~$h\leq h_\mbfs$ and let~$\nu^h$ be its representing measure. If~$s_j=0$ for some~$j\in [q]$, then~$h\leq h_\mbfs$ forces~$\nu^h \ttparen{\Delta^{q-1}\setminus \Delta^{q-1}_j}=0$ and the proof reduces to that of the analogous statement on a lower-dimensional simplex.
	Thus, we may and shall assume, with no loss of generality up to a simple induction argument on the dimension of the support of~$\nu^h$, that~$\mbfs\in \mathring{\Delta}^{q-1}$.
	If~$\nu^h\ll \delta_\mbfs$ then~$h\equiv \tau h_\mbfs$ for some~$\tau\in [0,1]$ and the proof is complete.
	Thus, we argue by contradiction that~$\nu^h\not\ll \delta_\mbfs$.	
	Let~$\mbfu\in\Delta^{q-1}$ with~$\mbfu \neq \mbfs$ and note that there exists~$j\in [q]$ such that~$u_j>s_j$.
	Thus, the sets~$U_j\eqdef \ttset{\mbfu\in\Delta^{q-1}: u_j>s_j}$ satisfy~$\Delta^{q-1}\setminus\set{\mbfs} = \cup_{j\in [q]} U_j$.
	Since~$\nu^h\not\ll \delta_\mbfs$, there exists at least one~$j\in [q]$ for which~$\nu^h U_j>0$.
	Then, by definition of~$U_j$ and by Fatou's lemma,
	\[
	\liminf_{n\to\infty} \frac{h(n\mbfe_j)}{h_\mbfs(n\mbfe_j)} \geq \liminf_{n\to\infty} \frac{1}{\mbfs^{n\mbfe_j}} \int_{U_j} \mbfu^{n\mbfe_j} \diff\nu^h(\mbfu) \ge \int_{U_j} \liminf_{n\to\infty} (u_j/s_j)^n \diff\nu^h(\mbfu) = \infty \comma
	\]
	which contradicts the assumption that~$h\leq h_\mbfs$.
	(Note that~$\mbfs^{n\mbfe_j}\neq 0$ since~$\mbfs\in\mathring{\Delta}^{q-1}$.)
\end{proof}

\subsubsection{Proofs} We are now ready to prove the characterization theorem.

\begin{proof}[Proof of Theorem~\ref{t:CharacterizationNew}]
We divide the proof into several steps.

\paragraph{Characterization of~$\BlockSupp{\probP_\bullet}$}
Let~$\mbfr \in \N_*^q$ and~$n \ge \length{\mbfr}$. Choosing~$A=\car_\mbfr + (n-\length{\mbfr}) \car_{\mbfe_1}$ in~\eqref{eq:Kingman(3.1)} we have
\[
	\probP_{n-\length{\mbfr}}[(n-\length{\mbfr}) \car_{\mbfe_1}] \, b(n,\mbfr)= \length{\mbfr} \, A(\mbfr) \, \probP_{n}[A] 
\]
The probability~$\probP_{n-\length{\mbfr}}[(n-\length{\mbfr}) \car_{\mbfe_1}]$ is strictly positive by~\ref{i:t:CharacterizationNew:1}-\ref{i:t:CharacterizationNew:2}. So is~$\probP_n[A]$ when~$\mbfr \in \BlockSupp{\probP_\bullet}$. Conversely, if~$\probP_n[A] > 0$, then~$\mbfr \in \BlockSupp{\probP_\bullet}$ by Lemma~\ref{l:BlockSupp}\ref{i:l:BlockSupp:2}. Therefore,
\begin{equation}\label{eq:t:CharacterizationProofs:1}
\mbfr\in\BlockSupp{\probP_\bullet} \iff b(n,\mbfr)>0 \comma \qquad n\in\N_1 \comma \mbfr \in \N^q_* \, : \, \length{\mbfr} \le n \fstop
\end{equation}

\paragraph{Step~1: Properties of~$b$}
Fix~$n\geq 2$. Given~$\mbfr,\mbfs\in\BlockSupp{\probP_\bullet}\setminus\set{\zero}$ with~$\length{\mbfr}+\length{\mbfs}\leq n$, note that $A\eqdef\car_\mbfr+\car_\mbfs+(n-\length{\mbfr}-\length{\mbfs})\car_{\mbfe_1}$ has~$\supp A=\set{\mbfr,\mbfs,\mbfe_1}\subset\BlockSupp{\probP_\bullet}$ by~\ref{i:t:CharacterizationNew:1} and the assumption on~$\mbfr,\mbfs$.
Thus,~$\probP_n[A]>0$ by~\ref{i:t:CharacterizationNew:2}, and, analogously,~$\probP_{n-\length{\mbfr}-\length{\mbfs}}[A-\car_\mbfr-\car_\mbfs]>0$.
Apply~\eqref{eq:Kingman(3.1)} twice to~$A\in\Part{q}{n}$ to obtain
\[
0<\probP_n[A] = \frac{b(n,\mbfr)}{\length{\mbfr}\, A(\mbfr)} \probP_{n-\length{\mbfr}}[A-\car_\mbfr] = \frac{b(n,\mbfr)}{\length{\mbfr}\, A(\mbfr)}  \frac{b(n-\length{\mbfr},\mbfs)}{\length{\mbfs}\, (A-\car_\mbfr)(\mbfs)} \probP_{n-\length{\mbfr}-\length{\mbfs}}[A-\car_\mbfr-\car_\mbfs] \fstop
\]
Exchanging~$\mbfr$,~$\mbfs$, and simplifying~$\probP_{n-\length{\mbfr}-\length{\mbfs}}[A-\car_\mbfr-\car_\mbfs]>0$, we see that~$b(n,\mbfr)$ must satisfy
\begin{equation}\label{eq:t:CharacterizationProofs:2}
b(n,\mbfr)\, b(n-\length{\mbfr},\mbfs) = b(n,\mbfs)\, b(n-\length{\mbfs},\mbfr) \comma \qquad \mbfr,\mbfs\in\BlockSupp{\probP_\bullet}\setminus\set{\zero} \, : \, \length{\mbfr}+\length{\mbfs}\leq n \fstop
\end{equation}

\paragraph{Step~2: $b$ on singletons}
Since~$\probP_1[\car_{\mbfe_j}]>0$ by assumption~\ref{i:t:CharacterizationNew:1}, we have~$p_j\eqdef \probP_1[\car_{\mbfe_j}]>0$ for every~$j\in [q]$. Since~$\probP_1$ is a probability measure,~$\mbfp\in\Delta^{q-1}$.
Taking~$\mbfr\eqdef \mbfe_j$ and~$\mbfs\eqdef \mbfe_k$ in~\eqref{eq:t:CharacterizationProofs:2} shows that~$b(n,\mbfe_j)/b(n,\mbfe_k)$ is independent of~$n$, hence equal to~$p_j/p_k$.
Thus %
\[
b(n,\mbfe_j) = n\beta_n p_j \qquad \text{for some}\quad \beta_n>0 \comma \quad \beta_1\eqdef 1\comma \qquad n\in\N_1 \fstop
\]

\paragraph{Step 3: General~$b$}
For~$\mbfr\in\BlockSupp{\probP_\bullet}\setminus\set{\zero}$ and~$n\geq \length{\mbfr}+1$, choosing~$\mbfs\eqdef\mbfe_j$ in~\eqref{eq:t:CharacterizationProofs:2} yields
\begin{equation}\label{eq:t:CharacterizationProofs:3}
\frac{b(n,\mbfr)}{b(n-1,\mbfr)} = \frac{b(n,\mbfe_j)}{b(n-\length{\mbfr},\mbfe_j)}= \frac{n\beta_n}{(n-\length{\mbfr})\beta_{n-\length{\mbfr}}} \fstop
\end{equation}
Note that the ratio on the right-hand side above depends on~$\mbfr$ only via~$\length{\mbfr}$. 
Now, set
\[
g(\mbfr)\eqdef b(\length{\mbfr},\mbfr)\Big/\prod_{i=1}^{\length{\mbfr}} i\beta_i \comma \qquad \mbfr\in\N^q_*\fstop
\]
By~\eqref{eq:t:CharacterizationProofs:1} we have $g(\mbfr)>0$ if~$\mbfr\in\BlockSupp{\probP_\bullet}\setminus\set{\zero}$, and~$g(\mbfr)= 0$ otherwise.
Solving~\eqref{eq:t:CharacterizationProofs:3} for~$b(n,\mbfr)$ and iterating in the right-hand side eventually yields
\begin{equation}\label{eq:t:CharacterizationProofs:4}
b(n,\mbfr) = g(\mbfr) \prod_{i=n-\length{\mbfr}+1}^n i\beta_i {\comma\qquad n\in\N_1 \comma \mbfr \in \N^q_* \, : \, \length{\mbfr} \le n \fstop}
\end{equation}
{(This formula holds as well for~$\mbfr \notin \BlockSupp{\probP_\bullet}$ by~\eqref{eq:t:CharacterizationProofs:1}.)} In particular, choosing~$\mbfr=\mbfe_j$ and~$n=1$ above, we have~$g(\mbfe_j)=p_j$.

\paragraph{Step 4: Consistency determines~$\beta$ and~$g$}
Iteratively substituting~\eqref{eq:t:CharacterizationProofs:4} into~\eqref{eq:Kingman(3.1)} shows
\begin{equation}\label{eq:t:CharacterizationProofs:5}
\probP_n[A] = \prod_{i=1}^n i\beta_i  %
 \prod_{\mbfa\in\supp A} \frac{g(\mbfa)^{A(\mbfa)}}{\length{\mbfa}^{A(\mbfa)}A(\mbfa)!}\comma \qquad A\in\Part{q}{n} \comma \quad n \in \N_1 \fstop
\end{equation}
Substituting~\eqref{eq:t:CharacterizationProofs:5} into~\eqref{eq:ConsistencySpecial} and simplifying, %
for each~$B\in\Part{q}{n-1}$ with~$\probP_{n-1}[B]>0$,
\begin{align}\label{eq:t:Characterization:2}
\frac{1}{\beta_n} = \sum_{k=1}^q g(\mbfe_k) + \sum_{\mbfr\in\supp B} B(\mbfr) f(\mbfr) = 1 + \sum_{\mbfr\in\supp B} B(\mbfr) f(\mbfr)\comma \qquad n \geq 2\comma
\end{align}
where we used~$g(\mbfe_k)=p_k$ and~$\length{\mbfp}=1$, and we set
\[
f(\mbfr)\eqdef \sum_{k=1}^q \length{\mbfr} \frac{g(\mbfr+\mbfe_k)}{g(\mbfr)} \frac{r_k+1}{\length{\mbfr}+1} \fstop
\]

We now specialize~\eqref{eq:t:Characterization:2} to three particular choices of~$B$, admissible (i.e.\ satisfying~$\probP_{n-1}[B]>0$) in light of assumptions~\ref{i:t:CharacterizationNew:1}-\ref{i:t:CharacterizationNew:2}:
\begin{itemize}%
\item With~$B\eqdef(n-1)\car_{\mbfe_j}$ we have~$\beta_n^{-1}= 1+ (n-1)f(\mbfe_j)$. Since the left-hand side does not depend on~$j\in [q]$ we have~$f(\mbfe_j)=\gamma$ for some~$\gamma\geq 0$ and~$\beta_n^{-1}=1+(n-1)\gamma$.
\item By~\ref{i:t:CharacterizationNew:3} there exists~$\mbfc=\mbfe_j+\mbfe_k\in \BlockSupp{\probP_\bullet}$, possibly with~$j=k$.
Since~$\mbfc\in \BlockSupp{\probP_\bullet}$, we have~$g(\mbfc)>0$. 
With~$n=2$ and~$B=\car_{\mbfe_k}$ we see that~$\gamma=f(\mbfe_j)\geq \frac{g(\mbfc)}{g(\mbfe_j)}\frac{\delta_{jk}+1}{2}>0$ by definition of~$f$.
Setting~$\theta\eqdef 1/\gamma\in (0,\infty)$ then gives
\begin{equation}\label{eq:t:CharacterizationProofs:8}
\beta_n=\theta/(\theta+n-1)\comma \qquad n\in \N_1 \fstop
\end{equation}
\item With~$B\eqdef \car_\mbfr+(n-1-\length{\mbfr})\car_{\mbfe_1}$ for~$\mbfr\in\BlockSupp{\probP_\bullet}\setminus\set{\zero}$ we have~$\beta_n^{-1}= 1+ f(\mbfr)+(n-1-\length{\mbfr})\gamma$,~hence
\begin{equation}\label{eq:t:CharacterizationProofs:9}
f(\mbfr)= \length{\mbfr}\theta^{-1} \comma \qquad \mbfr\in\BlockSupp{\probP_\bullet} \setminus\set{\zero}\fstop
\end{equation}
\end{itemize}

\paragraph{Step~5: Harmonicity and identification of~$\BlockSupp{\probP_\bullet}$}
Define~$h\colon \N_0^q\to\R$ by
\[
h(\zero)\eqdef 1\comma \qquad h(\mbfr)\eqdef g(\mbfr)\theta^{\length{\mbfr}-1} \binom{\length{\mbfr}}{\mbfr}^{-1} \mbfp^{-\mbfr} \comma \qquad \mbfr\in \N^q_* \fstop
\]
Then, using that~$\binom{\length{\mbfr}+1}{\mbfr+\mbfe_k}\frac{r_k+1}{\length{\mbfr}+1}=\binom{\length{\mbfr}}{\mbfr}$, the identity~\eqref{eq:t:CharacterizationProofs:9} becomes
\begin{equation}\label{eq:t:CharacterizationProofs:11}
\sum_{k=1}^q h(\mbfr+\mbfe_k) p_k = h(\mbfr)\comma \qquad \mbfr\in \N_0^q\fstop
\end{equation}
(For~$\mbfr=\zero$ it reads~$\sum_{k=1}^q p_k=1$ and for~$\mbfr\not\in\BlockSupp{\probP_\bullet}$ both sides vanish identically because~$\BlockSupp{\probP_\bullet}$ is order-closed by Lemma~\ref{l:BlockSupp}\ref{i:l:BlockSupp:1}.)
Equation~\eqref{eq:t:CharacterizationProofs:11} is precisely~\eqref{eq:d:Harmonic}, which shows that~$h$ is $\DRG{\mbfp}$-harmonic.
Proposition~\ref{p:RepresentationH} proves the existence of a unique Borel measure~$\nu$ on~$\Delta^{q-1}$ such that
$%
h(\mbfr) = \mbfp^{-\mbfr} \mom{\mbfr}{\nu} %
$ %
for every~$\mbfr\in\N^q_0$,
which identifies~$g$ as
\begin{equation}\label{eq:t:CharacterizationProofs:13}
g(\mbfr)= \theta^{1-\length{\mbfr}}\binom{\length{\mbfr}}{\mbfr} \mom{\mbfr}{\nu}\comma \qquad \mbfr\in \N^q_*\comma
\end{equation}
positive if and only if~$\mbfr\in \BlockSupp{\probP_\bullet}=\Sigma_\nu$. (Note that~$\zero$, which is not covered by the above discussion, is an element of both.)

\paragraph{Step~6: Identifications}
$\nu$ is a probability measure since~$h(\zero)=1$.
Substituting~\eqref{eq:t:CharacterizationProofs:13} and~\eqref{eq:t:CharacterizationProofs:8} into~\eqref{eq:t:CharacterizationProofs:4} and~\eqref{eq:t:CharacterizationProofs:5} proves that~$b\equiv b_*$ %
and~$\probP_n=\pEsf{n}{\theta,\nu}$.

\paragraph{Step~7: Uniqueness ($\theta$ and~$\nu$ are uniquely determined)}
We have $\sum_{\mbfa:\length{\mbfa}=2}\probP_2[\car_\mbfa]=\frac{1}{\theta+1}$, which determines~$\theta$.
Making~$b$ explicit in~\eqref{eq:Kingman(3.1)} and computing at~$A\eqdef \car_\mbfa$,
\[
\probP_{\length{\mbfa}}[\car_\mbfa] = \frac{\theta (\length{\mbfa}-1)!}{\Poch{\theta}{\length{\mbfa}}} \binom{\length{\mbfa}}{\mbfa}\mom{\mbfa}{\nu} \comma \qquad \mbfa\in \BlockSupp{\probP_\bullet} = \Sigma_\nu \comma \mbfa\neq \zero \comma
\]
which determines all (non-zero) moments of~$\nu$.
Since~$\Delta^{q-1}\subset \R^q$ is compact, the multi-dimensional Hausdorff moment problem on~$\Delta^{q-1}$ is well-posed, which uniquely fixes~$\nu$.%
\end{proof}

\subsection{Polychromatic Deletion Chain and Hoppe Urn}\label{ss:Hoppe}
The transition kernels~$\BackwardK^\sym{q}$ in~\eqref{eq:Kernels} define a reverse Markov chain~$\hat\msA_\proc$ incoming from infinity into the state space~$\Part{q}{\cup}$, describing deletion of elements uniformly at random according to Remark~\ref{r:KernelInterpretation}.

\begin{defs}
We call~$\hat\msA_\proc$ the \emph{polychromatic deletion chain}.
\end{defs}

Conditionally on choosing the block from which an element is removed, the evolution of each block of~$\hat\msA_\proc$ is the reverse Markov chain of the~$\DRG{\mbfp}$ in~\S\ref{sss:DRG}, for any~$\mbfp\in\mathring\Delta^{q-1}$.
In this section, we construct a time-reversal~$\msA_\proc$ of~$\hat\msA_\proc$.
We provide a description of~$\msA_\proc$ via kernels, more suitable for our purposes in later sections.
We will comment in Remark~\ref{r:UrnInterpretation} on why~$\msA_\proc$ is a polychromatic analogue of the urn scheme of F.~M.~Hoppe~\cite{Hop84}.

Throughout this section, we fix~$\theta>0$ and~$\nu\in\msP(\Delta^{q-1})$ with~$\nu\mathring{\Delta}^{q-1}>0$.
By~\eqref{eq:p:SupportProp:2} this is equivalent to~$\Sigma_\nu=\N_0^q$ with~$\Sigma_\nu$ as in~\eqref{eq:SuppMomentsNu}, i.e.\ to
\begin{equation}\label{eq:StandingAssumption}
\mom{\mbfa}{\nu}>0 \comma \qquad \mbfa\in\N_0^q \comma
\end{equation}
and, by~\eqref{eq:p:SupportProp:1}, %
it implies
\begin{equation}\label{eq:PositivePESF}
\pEsf{n}{\theta,\nu}[A]>0 \comma \qquad A\in \Part{q}{n}\comma n \in\N_0 \fstop
\end{equation}
Everywhere in this section %
we regard~$\Part{q}{0}=\set{0}$ as the~$0^{\text{th}}$ level of~$\Part{q}{\cup}$ and we recall that
\begin{equation}\label{eq:LevelZero}
\pEsf{0}{\theta,\nu}[0]\eqdef 1 \comma \qquad \BackwardK^\sym{q}_{0,n}(A,0)\eqdef 1 \comma \qquad A \in \Part{q}{n}\comma n\in\N_0\fstop%
\end{equation}

We define a system~$\ForwardK^\sym{q}=\ForwardK^\sym{q}_{\theta,\nu, m,n}$, $n,m\in \N_0$, $m\leq n$, of probability kernels on~$\Part{q}{\cup}$ by
\begin{subequations}\label{eq:HoppeKernel}
\begin{align}
\label{eq:HoppeKernel:A}
\ForwardK^\sym{q}_{m,m+1}\tparen{B, B+\car_{\mbfe_j}}&\eqdef \frac{\theta\, \mom{\mbfe_j}{\nu}}{\theta+m} \comma
\\
\label{eq:HoppeKernel:B}
\ForwardK^\sym{q}_{m,m+1}\tparen{B, B-\car_\mbfb+\car_{\mbfb+\mbfe_j}}&\eqdef \frac{\length{\mbfb}\, B(\mbfb)}{\theta+m}\, \frac{\mom{\mbfb+\mbfe_j}{\nu}}{\mom{\mbfb}{\nu}} \comma \qquad \mbfb\in\supp B\comma
\end{align}
\end{subequations}
and~$\ForwardK^\sym{q}_{m,m+1}(B,A)\eqdef 0$ for every other~$A\in\Part{q}{m+1}$.
Let~$\ForwardK^\sym{q}=\ForwardK^\sym{q}_{m,n}$ be the unique system of kernels extending~\eqref{eq:HoppeKernel} and satisfying the cocycle relation
\begin{equation}\label{eq:CocycleForward}
\ForwardK^\sym{q}_{\ell, n}(C,A)= \sum_{B\in\Part{q}{m}} \ForwardK^\sym{q}_{\ell,m}(C,B)\, \ForwardK^\sym{q}_{m,n}(B,A)\comma \qquad 
\begin{gathered}
A\in\Part{q}{n}\comma C\in\Part{q}{\ell}\comma 
\\
\quad 0\leq \ell < m < n\fstop
\end{gathered}
\end{equation}

Note that each~$A\in\Part{q}{m+1}$ arises in at most one of the two lines of~\eqref{eq:HoppeKernel}, and for at most one pair~$(\mbfb,j)$ in the second one: indeed~$\size{B+\car_{\mbfe_j}}=\size{B}+1$ while~$\size{B-\car_\mbfb+\car_{\mbfb+\mbfe_j}}=\size{B}$, and the multiset identity~$\car_{\mbfb+\mbfe_j}+\car_{\mbfb'}=\car_{\mbfb'+\mbfe_{j'}}+\car_{\mbfb}$ forces~$\mbfb=\mbfb'$ and~$j=j'$.
Thus~\eqref{eq:StandingAssumption} guarantees that the ratio in~\eqref{eq:HoppeKernel} is finite, and~$\ForwardK^\sym{q}$ is well-defined.

\begin{lem}\label{l:HoppeChainWellDefined}
$\ForwardK^\sym{q}(B,\emparg)$ is a probability on~$\Part{q}{m+1}$ for every~$m\in\N_0$ and~$B\in\Part{q}{m}$.

\begin{proof}
Since~$\mom{\zero}{\nu}=1$, Lemma~\ref{l:HoppeTransition}\ref{i:l:HoppeTransition:2} applied to~$\mbfb=\zero$ gives~$\sum_{j=1}^q \mom{\mbfe_j}{\nu}=1$, whence~\eqref{eq:HoppeKernel:A} sums to~$\theta/(\theta+m)$ over~$j\in [q]$.
As for~\eqref{eq:HoppeKernel:B}, again by Lemma~\ref{l:HoppeTransition}\ref{i:l:HoppeTransition:2},
\[
\sum_{\mbfb\in\supp B} \frac{\length{\mbfb} B(\mbfb)}{\theta+m}\, \frac{1}{\mom{\mbfb}{\nu}}\sum_{j=1}^q \mom{\mbfb+\mbfe_j}{\nu} = \sum_{\mbfb\in\supp B} \frac{\length{\mbfb} B(\mbfb)}{\theta+m} = \frac{\length{\col{B}}}{\theta+m}=\frac{m}{\theta+m}\comma
\]
and the two contributions sum up to~$1$.
\end{proof}
\end{lem}

\begin{defs}\label{d:HoppeChain}
The \emph{polychromatic Hoppe chain} (in short: PHC) with parameters~$\theta,\nu$ is the Markov chain~$\msA_\proc\eqdef\seq{\msA_t}_{t\in\N_0}$ with state space~$\Part{q}{\cup}$ and transition kernel~$\ForwardK^\sym{q}$.
We write~$\probP^B$ for the law of~$\msA_\proc$ started at~$B$, and~$\ExpectP^B$ for the corresponding expectation.
\end{defs}

We now show that the kernels~\eqref{eq:Kernels} of~$\hat\msA_\proc$ are the co-transition kernels of the PHC~$\msA_\proc$.

\begin{lem}[Time reversal]\label{l:Reversal}
For every~$n\in\N_1$, every~$A\in\Part{q}{n}$, and every~$B\in\Part{q}{n-1}$,
\begin{equation}\label{eq:l:Reversal:0}
\ForwardK^\sym{q}_{n-1,n}(B,A)\, \pEsf{n-1}{\theta,\nu}[B] = \pEsf{n}{\theta,\nu}[A]\, \BackwardK^\sym{q}_{n-1,n}(A,B) \fstop
\end{equation}
In particular,~\eqref{eq:CocycleForward},~\eqref{eq:Cocycle} is a forward-backward pair of Chapman--Kolmogorov equations.

\begin{proof}
Both sides of~\eqref{eq:l:Reversal:0} vanish unless~$B=A_{\setminus \mbfa,j}$ for some (unique, cf.~\eqref{eq:AmenoEj})~$\mbfa\in\supp A$ and~$j\in [q]$ with~$a_j\geq 1$.
Indeed, this is the case for the right-hand side by~\eqref{eq:Kernels:B}; as for the left-hand side, $A=B+\car_{\mbfe_j}$ reads~$B=A_{\setminus \mbfe_j,j}$, while~$A=B-\car_\mbfb+\car_{\mbfb+\mbfe_j}$ reads~$B=A_{\setminus \mbfa,j}$ with~$\mbfa\eqdef\mbfb+\mbfe_j$.
Now, let~$B=A_{\setminus \mbfa,j}$ and set~$m\eqdef n-1$.
Making use of the definition~\eqref{d:eq:PolyESF} of~$\pEsf{n}{\theta,\nu}$, simple computations yield
\begin{equation}\label{eq:l:Reversal:1}
	\mom{\mbfe_j}{\nu} \frac{\theta}{\theta+m} \pEsf{m}{\theta,\nu}\bracket{A_{\setminus \mbfe_j,j}} = \pEsf{m+1}{\theta,\nu}[A] \frac{A(\mbfe_j)}{m+1}
\end{equation}
when~$\mbfa=\mbfe_j$ and~$A \ge \car_{\mbfe_j}$, and
\begin{equation}\label{eq:l:Reversal:2}
	(A(\mbfa - \mbfe_j)+1) \frac{\mom{\mbfa}{\nu}}{\mom{\mbfa - \mbfe_j}{\nu}} \frac{\length{\mbfa}-1}{\theta + m}  \pEsf{m}{\theta,\nu}\bracket{A_{\setminus \mbfa,j}} = \pEsf{m+1}{\theta,\nu}[A] \frac{a_j A(\mbfa)}{m+1}
\end{equation}
when~$\length{\mbfa} \ge 2$, $a_j \ge 1$, and~$A\geq \car_\mbfa$.
If~$\mbfa=\mbfe_j$, then~$\ForwardK^\sym{q}_{m,m+1}(B,A)=\theta\mom{\mbfe_j}{\nu}/(\theta+m)$ and $\BackwardK^\sym{q}_{n-1,n}(A,B)=A(\mbfe_j)/n$, so that~\eqref{eq:l:Reversal:0} is precisely~\eqref{eq:l:Reversal:1}.
If~$\length{\mbfa}\geq 2$, set~$\mbfb\eqdef \mbfa-\mbfe_j\in\supp B$ and note that~$B(\mbfb)=A(\mbfa-\mbfe_j)+1$ and~$\length{\mbfb}=\length{\mbfa}-1$, so that
\[
\ForwardK^\sym{q}(B,A)= \tparen{A(\mbfa-\mbfe_j)+1}\frac{\mom{\mbfa}{\nu}}{\mom{\mbfa-\mbfe_j}{\nu}} \frac{\length{\mbfa}-1}{\theta+m} \comma \qquad \BackwardK^\sym{q}_{n-1,n}(A,B)= \frac{a_j A(\mbfa)}{n} \comma
\]
and~\eqref{eq:l:Reversal:0} is precisely~\eqref{eq:l:Reversal:2}.
\end{proof}
\end{lem}

\begin{cor}\label{c:HoppeMarginals}
The one-dimensional time-marginals of~$\msA_\proc$ started at~$\msA_0=0$ are
\begin{equation}\label{eq:HoppeMarginals}
\probP^0[\msA_n=A]=\pEsf{n}{\theta,\nu}[A] \comma \qquad A\in\Part{q}{n}\comma n\in\N_0 \fstop
\end{equation}
Consequently~$\probP^0\tbracket{\msA_m=B \mid \msA_n=A}=\BackwardK^\sym{q}_{m,n}(A,B)$ for all~$m\leq n$, i.e.~$\hat\msA_\proc$ is the time reversal of~$\msA_\proc$ under~$\probP^0$.

\begin{proof}
We prove~\eqref{eq:HoppeMarginals} by induction on~$n$. The base step ($n=0$) holds by~\eqref{eq:LevelZero} since~$\msA_0=0$.
\paragraph{Inductive step}
Assuming~\eqref{eq:HoppeMarginals} for~$n-1$, and since~$\BackwardK^\sym{q}_{n-1,n}(A,\emparg)$ is a probability measure on~$\Part{q}{n-1}$, Lemma~\ref{l:Reversal} gives
\[
\probP^0[\msA_n=A] = \sum_{B\in\Part{q}{n-1}} \ForwardK^\sym{q}_{n-1,n}(B,A)\, \pEsf{n-1}{\theta,\nu}[B] = \pEsf{n}{\theta,\nu}[A] \sum_{B\in\Part{q}{n-1}} \BackwardK^\sym{q}_{n-1,n}(A,B) = \pEsf{n}{\theta,\nu}[A] \fstop
\]
The last assertion follows from~\eqref{eq:l:Reversal:0},~\eqref{eq:HoppeMarginals} and~\eqref{eq:PositivePESF} by Bayes' Formula for~$m=n-1$.%
\end{proof}
\end{cor}

\begin{rem}[The polychromatic Hoppe urn]\label{r:UrnInterpretation}
In~\cite{Hop84}, F.~M.~Hoppe showed that the ESF~$\Esf{n}{\theta}$ is the marginal distribution of a discrete-time Markov process~$\seq{\Pi_n}_n$ of integer partitions~$\Pi_n\vdash n$ obtained from the sampling process~$\seq{X_n}_n$ of what is now known as \emph{Hoppe's urn model}.
We adapt his construction to a similar urn model, resulting in a $\Part{q}{\cup}$-valued Markov process with marginal distribution~$\pEsf{n}{\theta,\nu}$ at time~$n$.

The PHC of Definition~\ref{d:HoppeChain} may be described as the polychromatic analogue of Hoppe's urn model, projected on~$\Part{q}{\cup}$.
We consider a process~$Y_\proc\eqdef\seq{Y_n}_n$ generated by sampling from an urn containing one cube and various numbers of labeled colored balls.
At time~$0$, the urn contains only the cube.
At every (integer) time~$m$, the labels are consecutive and ranging in~$\N_1$, while the colors range in~$[q]$.
The cube has mass~$\theta$ and every ball has mass~$1$.
At time~$m$, an object in the urn is drawn at random with a probability proportional to its mass.
If it is the cube, it is returned together with a ball with the smallest label not yet present in the urn and of a randomly chosen color.
For each~$j \in [q]$, the color~$j$ is chosen with probability~$\mom{\mbfe_j}{\nu}$.
If, instead, the drawn object is a ball, with label~$i$, it is returned together with one additional ball of the same label and of a color randomly chosen as follows: if~$a_j$ denotes, for each color~$j \in [q]$, the number of existing balls of label~$i$ and color~$j$ in the urn, the probability of the new ball having color~$j$ is~$\mom{\mbfa+\mbfe_j}{\nu}/\mom{\mbfa}{\nu}$, well-defined by~\eqref{eq:StandingAssumption}.

We define random variables~$r_m\in \N_0$ and~$Y_m\in \N_0\times [q]$ as the number of distinct labels (i.e.\ the maximal label) present in the urn, and the label and color of the additional ball returned after the $m^\text{th}$ drawing.
Observe that, for every~$n\in \N_1$, the process~$Y_\proc$ defines a random $q$-chromatic partition~$\msA_n$ by letting
\begin{gather} \label{eq:YtoA}
\mbfa_n(i)\eqdef \seq{a_{n,j}(i)}_{j\in [q]}\comma a_{n,j}(i)\eqdef\card{\set{m\in [n]: Y_m=(i,j)}}\comma
\qquad
\msA_n\eqdef\sum_{i=1}^{r_n}\car_{\mbfa_n(i)} \fstop
\end{gather}
As a consequence, in the notation of~\cite{Hop84}, the first component~$Y_{m,1}$ of~$Y_m$ satisfies~$Y_{m,1}=X_m$, while~$\shape(\msA_n)$ coincides with~$\Pi_n$.

At time~$m$, the urn contains the cube, of mass~$\theta$, and~$m$ balls, of mass~$1$ each, so that the total mass is~$\theta+m$, and that the multiset~$\msA_m$ in~\eqref{eq:YtoA} records, for each~$\mbfb\in\N^q_*$, the number~$\msA_m(\mbfb)$ of labels carried by exactly~$b_j$ balls of color~$j$, for every~$j\in[q]$.
Conditionally on the configuration of the urn at time~$m$, and writing~$B\eqdef\msA_m$:
\begin{itemize}
\item the cube is drawn with probability~$\theta/(\theta+m)$, and the ball returned with it carries a label not yet present and the color~$j$ with probability~$\mom{\mbfe_j}{\nu}$; the polychromatic partition thus passes from~$B$ to~$B+\car_{\mbfe_j}$, which accounts for~\eqref{eq:HoppeKernel:A};
\item since there are~$B(\mbfb)$ labels with color count~$\mbfb$, each carried by~$\length{\mbfb}$ balls, a ball whose label has color count~$\mbfb$ is drawn with probability~$\length{\mbfb}B(\mbfb)/(\theta+m)$, and the ball returned with it carries the same label and the color~$j$ with probability~$\mom{\mbfb+\mbfe_j}{\nu}/\mom{\mbfb}{\nu}$; the polychromatic partition thus passes from~$B$ to~$B-\car_\mbfb+\car_{\mbfb+\mbfe_j}$, which accounts for~\eqref{eq:HoppeKernel:B}.
\end{itemize}
\end{rem}

\subsection{Other models}\label{ss:Others}
We give here a \emph{combinatorial} interpretation of~\eqref{eq:t:Zn:0}, i.e.\ of the multinomial coefficient~$\Multi{\mbfn}$, in terms of %
\emph{necklaces}.
In turn, this provides a connection to pESF via the extension of the Chinese Restaurant Process which we now describe.%

For~$q\in\N_1$ denote by~$[q]^*$ the free monoid generated by~$[q]$.
This is the monoid whose elements are all the finite sequences of zero or more elements from~$[q]$, with string concatenation as the monoid operation and with the empty string as the identity element.
Elements of~$[q]^*$ are called ($q$-)\emph{words}.
Two words~$u,v$ are \emph{conjugate} if there exist words~$s,t$ so that~$u=st$ and~$v=ts$.
Two conjugate words are cyclic shifts of one another.
Thus, conjugacy is an equivalence relation on words.
Its equivalence classes are called ($q$-)\emph{necklaces}.

Let~$\varsigma=\class{w}$ be a necklace and~$w=c_1c_2\cdots c_\ell$ be any of its representatives.
The \emph{length}~$\ell_\varsigma$ of~$\varsigma$ is the total number~$\ell$ of characters in~$w$.
The \emph{period}~$p_\varsigma$ of~$\varsigma$ is the minimal integer~$p\geq 1$ with~$c_i=c_{(i+p-1 \pmod \ell) +1}$ for every~$i\in [\ell]$.
Clearly, %
$p_\varsigma$ divides $\ell_\varsigma$. %

\subsubsection{An extension of the Chinese Restaurant Process}\label{ss:CRP}
In~\cite{Ald85}, D.J.~Aldous introduced\footnote{In fact, Aldous credits the introduction of the CRP to J.~Pitman, who in turn acknowledges the contribution of L.~Dubins, see e.g.\ the attribution to Dubins and Pitman in~\cite[\S4.1]{Tav21}.}
the \emph{Chinese restaurant process} (CRP) as a sequential description of the sampling of random partitions distributed according to the Poisson--Dirichlet distribution.
The process (and many of its variations) has proven a very successful tool in the study of random partitions/permutations.
In the following example we state a variation of the CRP describing polychromatic partitions.

\paragraph{Dishes at the Chinese restaurant}
As for the standard CRP, \guillemotleft \emph{[customers]~$1,2,\dotsc, n$ arrive sequentially at an initially empty restaurant with a large number of large [circular] tables.
[Customer]~$j$ either sits at the same table as [customer]~$i$, with probability~$1/(j-1+\theta)$ for each~$i<j$, or else sits at an empty table, with probability~$\theta/(j-1+\theta)$}\guillemotright~\cite[(11.19), p.~91]{Ald85}.
Additionally to the standard CRP rules, %
however, each customer randomly chooses to order one out of the~$q$ proposed dishes, independently of the other customers and according to a fixed categorical distribution with parameter~$\mbfp$.
As a consequence, the polychromatic partition `customers at each table ordering each dish' is drawn according to~$\pEsf{n}{\theta,\delta_\mbfp}$.

\paragraph{Waiters and chefs at the Chinese restaurant}
On the one hand, waiters in our busy restaurant take care to remember, for every table, which clients order each dish.
The arrangement of the customers around the table is important in serving them efficiently.
All the information the waiters need about the customers' arrangement is thus contained in a $q$-chromatic necklace.
On the other hand, chefs in the restaurant only care about how many customers at each table order each dish, so that customers at the same table may be served at the same time.
All the information the chefs need about the customers' arrangement is thus contained in a $q$-chromatic partition.

We may therefore count $q$-chromatic partitions by collecting together $q$-chromatic necklaces with the same occurrences of each color (dish).
Indeed, let~$N$ be a multiset of necklaces describing an arrangement of the customers at each table and of the dish each of them ordered.
We write $N\in \mcN_A$ if the multiset of necklaces of clients ordering the same dish is encoded by~$A\vDash \mbfn$.
The number of permutations inducing a multiset~$N$ of $q$-chromatic necklaces with fixed occurrences of each dish is
$
\mbfn! \prod_{\varsigma\in\supp{N}} \frac{(p_\varsigma/\ell_\varsigma)^{N(\varsigma)}}{N(\varsigma)!}  %
$,
hence, the multinomial coefficient of~$A$ satisfies, cf.~Lemma~\ref{l:MultiA},
\[ 
\Multi{\mbfn}(A)=\sum_{N\in \mcN_A} \mbfn! \prod_{\varsigma\in\supp{N}} \frac{(p_\varsigma/\ell_\varsigma)^{N(\varsigma)}}{N(\varsigma)!} \fstop
\]

\section{Polychromatic Kingman Representation}\label{s:PolychromaticKingman}
In this section, we will state a polychromatic version of Kingman's Representation Theorem~\cite{Kin78} for partition structures:
in~\S\ref{ss:KingmanNotation} we collect the notation needed to state the theorem;
in~\S\ref{ss:Kingman} we state the Representation Theorem.
We then proceed to collect the many results necessary to restate and prove the Representation Theorem in terms of the Martin boundary of a reverse Markov chain;
in~\S\ref{ss:Functoriality} we discuss how the representation changes under color degeneracy;
in~\S\S\ref{ss:RepresentationPESF} and~\ref{sss:Strahov} we prove the representation of the pESF and specializations to the case of multipartitions.

\subsection{Notation}\label{ss:KingmanNotation}
Similarly to the monochromatic case, we will identify a Bauer simplex indexing the set~$\CPS{q}$ by classifying its extrema.
This simplex will be the space of probability measures over a $q$-version of Kingman's (extended ordered) simplex~$\nabla_0$ which we now define.

\begin{defs}
Fix~$q\in\N_1$.
We call \emph{polychromatic Kingman simplex~$\KSimplex{q}$} the space of infinite sequences~$\mbfx\eqdef \seq{\mbfx_i}_{i\geq 0}$ of vectors $\mbfx_i \eqdef\seq{x_{i,1},\dotsc, x_{i,q}}$ satisfying
\begin{equation}\label{eq:KingmanSimplex}
x_{i,j}\geq 0 \comma \quad i\geq 0\comma j\in [q] \comma \qquad \mbfx_i\geq_\shortlex \mbfx_{i+1} \comma \quad i\geq 1\comma \qquad \sum_{i=0}^\infty \length{\mbfx_i} =1 \fstop
\end{equation}
(Note that we allow~$\length{\mbfx_0}<\length{\mbfx_1}$.)
\end{defs}

For the moment, we endow~$\KSimplex{q}$ with the product topology~$\T_\mrmp$ inherited from~$\mbbI_q\eqdef ([0,1]^q)^{\N_0}$.
Note that~$\KSimplex{q}$ is a (convex) Borel subset of a Polish space, and thus a standard Borel space when endowed with its Borel $\sigma$-algebra.
As for Kingman's simplex, note that~$\nabla_0=\KSimplex{1}$.

\subsubsection{Extremal partition structures} \label{sec:extremal}
{Fix~$q\in\N_1$ and, %
for each~$r\in\N_0$, define the set 
\begin{equation}\label{eq:IndicesT}
T_r\eqdef\set{\iota\colon [r]\to\N_1 \quad \text{injective}}\comma \qquad T_0 = \set{\emp} \comma
\end{equation}
so that~$\seq{\iota(h)}_{h\in [r]}$ is any tuple of~$r$ distinct positive integer indices.

We will express a generic $\CPS{q}$ in terms of \emph{extremal} ones, labeled by points in~$\KSimplex{q}$.
Firstly, for each~$n\in\N_0$ and each~$A\in\Part{q}{n}$ with enumeration~$(\mbfa_h)_{h \in [\size{A}]}$} (cf.~Dfn.~\ref{d:Enumeration}), let us define the function~$\varphi_A\colon \KSimplex{q}\to [0,\infty)$, by
\begin{equation}\label{eq:DefPhi}
\varphi_0\equiv 1\comma\qquad\varphi_A\colon \mbfx \longmapsto n! \paren{\prod_{\mbfa\in \supp A}\frac{1}{\mbfa!^{A(\mbfa)} A(\mbfa)!}} \sum_{\iota \in T_{\size{A}}} \prod_{h=1}^{\size{A}} \mbfx_{\iota(h)}^{\mbfa_h} \comma \quad \mbfx\in\KSimplex{q} \fstop
\end{equation}
Note that~$\varphi_A(\mbfx)$ does not depend on the variables~$\mbfx_0$.
Indeed, we will show that~$\varphi_A(\mbfx)$ describes the extremal $\CPS{q}$ indexed by~$\mbfx\in\KSimplex{q}$ with~$\mbfx_0=\zero$.
The function~$\varphi_A$ is well-defined and Borel measurable by Corollary~\ref{c:PhiBorel} below. 

Now, for each~$\mbfm\in\N_0^q$ define the multiset
\begin{equation}
\uni{\mbfm}\eqdef m_1\car_{\mbfe_1}+\cdots + m_q\car_{\mbfe_q}\fstop
\end{equation}
A generic extremal $\CPS{q}$ is described by the function~$\psi_A\colon \KSimplex{q}\to [0,\infty)$ defined, for each~$n\in\N_0$ and each~$A\in\Part{q}{n}$, as
\begin{equation}\label{eq:DefPsi}
\psi_0\equiv 1\comma \qquad \psi_A(\mbfx)\eqdef n! \sum_{\substack{\mbfm\in \N_0^q\\ \uni{\mbfm} \leq A}} \frac{\mbfx_0^\mbfm}{\mbfm! (n-\length{\mbfm})!} \varphi_{A-\uni{\mbfm}}(\mbfx)\comma  \quad \mbfx\in\KSimplex{q}
\fstop
\end{equation}

We derive an alternative formula for~$\varphi_A(\mbfx)$ as follows. We denote by~$\boldnu\eqdef\seq{\boldnu_i}_{i\geq 1}$ vectors of indices
\begin{equation}\label{eq:l:TotalPhi:1.2}
	\boldnu_i\eqdef \seq{\nu_{i,1},\dotsc, \nu_{i,q}} \qquad \text{with} \qquad 0\leq \nu_{i,j}\leq n \qquad \text{and} \qquad \sum_{i=1}^\infty \length{\boldnu_i}=n \fstop
\end{equation}
To each sequence~$\boldnu$ as above we may assign a polychromatic partition~$A_\boldnu\in\Part{q}{n}$ by setting
\begin{equation}\label{eq:l:TotalPhi:1.5}
	A_\boldnu(\mbfa)\eqdef \card{\set{i : \boldnu_i=\mbfa}} \comma \qquad \mbfa\in\N^q_*\fstop
\end{equation}
For each~$\boldnu$ such that~$A_{\boldnu} = A$, there are exactly~$\prod_{\mbfa} A(\mbfa)!$ many injective maps~$\iota \in T_{\size{A}}$ such that~$\boldnu_{\iota(h)} = \mbfa_h$ for every~$h \in [\size{A}]$, and, moreover,~$\prod_{\mbfa \in \supp{A}} \mbfa!^{A(\mbfa)} = \prod_{i=1}^\infty \boldnu_i!$. We thus conclude that
\begin{equation} \label{eq:alterPhi}
	\varphi_A(\mbfx) = n! \sum_{\boldnu : A_\boldnu=A}  \prod_{i=1}^\infty \frac{\mbfx_i^{\boldnu_i}}{\boldnu_i!}  \fstop
\end{equation}
Consequently, we also have
\begin{equation} \label{eq:alterPsi}
	\psi_A(\mbfx) = n! \sum_{(\mbfm,\boldnu)\in {\mathcal I}_A} \frac{\mbfx_0^\mbfm}{\mbfm!} \prod_{i=1}^\infty \frac{\mbfx_i^{\boldnu_i}}{\boldnu_i!}\comma
\end{equation}
where the summation runs over the set~${\mathcal I}_A$ of all pairs~$(\mbfm,\boldnu)$ with~$\mbfm\in \N_0^q$ and~$\boldnu\eqdef\seq{\boldnu_i}_{i\geq 1}$,  $\boldnu_i\eqdef\seq{\nu_{i,j}}_{j\in [q]}$ such that~$\uni{\mbfm}+A_{\boldnu}=A$.
(Note that the factor~$(n-\length{\mbfm})!$ cancels the same factor in the definition of~$\varphi_{A-\uni{\mbfm}}$.)

\begin{rem}[$q=1$]
	Formula~\eqref{eq:alterPhi} shows that~$\varphi_A(\mbfx)$ is the correct polychromatic generalization of dust-free extremal partition structures in~\cite[Eqn.~(2.10)]{Kin78b}.
\end{rem}

\begin{rem} \label{rem:probKing}
	Fix~$n \in \N_0$ and~$\mbfx \in \KSimplex{q}$. Let~$(\xi_r)_{r \in \N_1}$ be independent random variables with
	\begin{equation} \label{eq:xir}
		\probP\tbraket{\xi_r=(i,j)}=x_{i,j}\comma \quad i\in\N_1 \comma \qquad
		\probP\tbraket{\xi_r=(-r,j)}=x_{0,j}\comma \qquad j\in [q]\fstop
	\end{equation}
	For fixed~$A \in \Part{q}{n}$ and~$(\mbfm,\boldnu) \in {\mathcal I}_A$, the quantity
	\[
		n! \, \frac{\mbfx_0^\mbfm}{\mbfm!} \prod_{i=1}^\infty \frac{1}{\boldnu_i!} \mbfx_i^{\nu_i} = \binom{n}{m_1, \dotsc, m_q, \nu_{1,1}, \dotsc, \nu_{1,q}, \nu_{2,1},\dotsc, \nu_{2,q},\dotsc } \, \prod_{j=1}^q x_{0,j}^{m_j} \, \prod_{i=1}^\infty \prod_{j=1}^q  x_{i,j}^{\nu_{i,j}}
	\]
	appearing in~\eqref{eq:alterPsi} equals the probability that, for every~$(i,j)$ with~$i \ge 1$, the number of indices~$r \in [n]$ with~$\xi_r = (i,j)$ equals~$\nu_{i,j}$, and that, for every~$j$, the number of indices~$r \in [n]$ with~$\xi_r = (-r,j)$ equals~$m_j$. This is the polychromatic analogue of Kingman's probabilistic description of extremal partition structures in~\cite[p.~375]{Kin78}. 
\end{rem}

\begin{lem} \label{l:TotalPhi}
	For every~$q\in\N_1$, every~$n\in\N_0$, %
	and every~$\mbfx\in\KSimplex{q}$,
	\begin{subequations}\label{eq:TotalPhiPsi:0}
		\begin{align} \label{eq:TotalPhiPsi:0a}
			(1-\length{\mbfx_0})^n &=\paren{\sum_{i=1}^\infty \length{\mbfx_i}}^n = \sum_{A\in\Part{q}{n}} \varphi_A(\mbfx) \comma
			\\
			\label{eq:TotalPhiPsi:0b}
			1&=\paren{\sum_{i=0}^\infty \length{\mbfx_i}}^n = \sum_{A\in\Part{q}{n}} \psi_A(\mbfx) \fstop
		\end{align}
	\end{subequations}
	In particular, for every~$q\in\N_1$,
	every~$n\in\N_0$, and every~$\mbfx\in\KSimplex{q}$, the map~$A \mapsto \psi_A(\mbfx)$ defines a probability measure on~$\Part{q}{n}$.
	
	\begin{proof}
		By the Multinomial Theorem and~\eqref{eq:alterPhi},
		\begin{equation}\label{eq:l:TotalPhiPsi:1}
			 \paren{\sum_{i=1}^\infty \sum_{j=1}^q x_{i,j}}^n = n!\sum_{\boldnu : \, \sum_i \length{\boldnu_i} = n} \prod_{i=1}^\infty \frac{\mbfx_i^{\boldnu_i}}{\boldnu_i!} = \sum_{A \in \Part{q}{n}} n!\sum_{\boldnu  : \, A_\boldnu = A} \prod_{i=1}^\infty \frac{\mbfx_i^{\boldnu_i}}{\boldnu_i!} = \sum_{A \in \Part{q}{n}} \varphi_A(\mbfx) \comma
		\end{equation}
		which proves~\eqref{eq:TotalPhiPsi:0a}. The proof of~\eqref{eq:TotalPhiPsi:0b} is done similarly, using~\eqref{eq:alterPsi} in place of~\eqref{eq:alterPhi}. The map~$A \mapsto \psi_A(\mbfx)$ is a probability measure on~$\Part{q}{n}$ by~\eqref{eq:TotalPhiPsi:0b} and by the non-negativity of~$\psi_A$, evident in both~\eqref{eq:DefPsi} and~\eqref{eq:alterPsi}.
	\end{proof}
\end{lem}

\begin{cor}\label{c:PhiBorel}
For every~$n\in\N_0$ and every~$A\in\Part{q}{n}$, the functions~$\varphi_A,\psi_A\colon \KSimplex{q}\to [0,\infty)$ are well-defined and $\T_\mrmp$-lower semicontinuous. %
\begin{proof}
Let~$\chi$ be any monomial in \emph{finitely many} coordinates~$x_{i,j}$ of~$\mbfx\in\KSimplex{q}$ ---as are the products in~\eqref{eq:DefPhi}.
Since all~$\chi$'s are non-negative, the series in~\eqref{eq:DefPhi} converges absolutely by Lemma~\ref{l:TotalPhi}, and thus~$\varphi_A$ is well-defined.
Since all~$\chi$'s are continuous on~$\KSimplex{q}$, $\varphi_A$ is l.s.c., being a monotone increasing limit of a family of continuous functions.
That~$\psi_A$ too is well-defined and l.s.c.\ holds since it is a (finite) linear combination of~$\varphi_B$'s.
\end{proof}
\end{cor}

\subsection{Polychromatic Kingman Representation Theorem}\label{ss:Kingman}
We have the following:%

\begin{thm}[Polychromatic Kingman Representation]\label{t:KingmanRepresentation}
For every Borel probability measure~$\mu\in\msP(\KSimplex{q})$ and every~$n \in \N_0$, the assignment
\begin{equation}\label{eq:p:KingmanRepresentationConverse:0}
\probP_n\colon A\longmapsto \int_{\KSimplex{q}} \psi_A(\mbfx) \diff\mu(\mbfx) \comma \qquad A\in\Part{q}{n}\comma
\end{equation}
defines a probability measure~$\probP_n$ on~$\Part{q}{n}$, and the family~$\seq{\probP_n}_{n\geq 0}$ is in $\CPS{q}$.

Conversely, every $\CPS{q}$ is represented by~\eqref{eq:p:KingmanRepresentationConverse:0} for one and only one~$\mu\in\msP(\KSimplex{q})$.
\end{thm}

Before discussing the proof, let us comment further on the assumptions in Characterization Theorem~\ref{t:CharacterizationNew} for the pESF.
\begin{rem}\label{r:NecessityAssChar}
Let~$\mbfp\in\mathring\Delta^{q-1}$.
\begin{enumerate*}[$(1)$]
\item According to Corollary~\ref{c:AsymptoticsPESF}, the~$\CPS{q}$ $\pEsf{\bullet}{0,\delta_\mbfp}$ is concentrated only on polychromatic partitions consisting of a single block with colors drawn from a multinomial distribution of parameters~$n$ and~$\mbfp\in \Delta^{q-1}$, and thus it is represented by~$\mu=\delta_{(\zero,\mbfp,\zero,\zero,\dotsc)}$.
It is not difficult to show that it satisfies assumptions~\ref{i:t:CharacterizationNew:1} and~\ref{i:t:CharacterizationNew:3} in Theorem~\ref{t:CharacterizationNew}, as well as~\ref{i:t:CharacterizationNew:4} with~$b(n,\mbfa)= n \binom{n}{\mbfa}\mbfp^\mbfa$ when~$\length{\mbfa}=n$ and~$b(n,\mbfa)=0$ otherwise, but it does not satisfy~\ref{i:t:CharacterizationNew:2}.

\item Analogously, the~$\CPS{q}$ $\pEsf{\bullet}{\infty,\delta_\mbfp}$ is concentrated only on polychromatic partitions consisting of all singletons, with i.i.d.\ colors drawn from~$\Cat_\mbfp$, and thus it is represented by~$\delta_{(\mbfp,\zero,\zero,\dotsc)}$.
It is not difficult to show that it satisfies~\ref{i:t:CharacterizationNew:1} and~\ref{i:t:CharacterizationNew:2}, and also~\ref{i:t:CharacterizationNew:4} with~$b(n,\mbfe_j)= n p_j$ (and vanishing otherwise), but it does not satisfy~\ref{i:t:CharacterizationNew:3}.
\end{enumerate*}
\end{rem}

\subsubsection{Proof structure}\label{sss:KingmanRepresentationDiscussion}
We will start in~\S\ref{sss:MartinCompactness} with the definition of a compact Polish topology~$\T_*$ on~$\KSimplex{q}$, different from~$\T_\mrmp$, obtained by `compactifying vanishing sequences at pure-dust points'.
This topology will make the functions~$\mbfx\mapsto \psi_A(\mbfx)$ (with fixed~$A\in\Part{q}{\cup}$) continuous (rather than only l.s.c.\ as w.r.t.~$\T_\mrmp$) and their linear span a unital point-separating sub-algebra of~$\Cb(\KSimplex{q},\T_*)$, which will be used for the uniqueness part of Theorem~\ref{t:KingmanRepresentation}.

In~\S\ref{sss:AuxCombinatorics} we show that~$\psi_\bullet(\mbfx)$ is a~$\CPS{q}$. %
In~\S\ref{sss:PHCRevisited} we move on to revisit the PHC by characterizing its harmonic densities in a one-to-one correspondence with polychromatic partition structures.  %
We then address (\S\ref{sss:MartinConvergence}) the convergence of the Martin kernels of the PHC by identifying~$\KSimplex{q}$ as the Martin boundary of~$\Part{q}{\cup}$.
This step is essential in that the convergence of the Martin kernels allows us to prove the convergence of block-color frequencies, which cannot be addressed by the original martingale argument in Kingman's proof~\cite{Kin78b}.

In~\S\ref{sss:MainProofs} we prove the Representation Theorem by randomizing the convergence of the Martin kernels over (what will turn out to be) the mixing measure~$\mu$.

{In~\S\ref{sss:Minimality}, we complete the proof of the triple correspondence stated in the Introduction.}

We present auxiliary results in full detail.
In comparison with~\cite{Kin78}, this is necessary:
\begin{itemize*}[]
\item due to the additional presence of the coloring;
\item since integer partitions are replaced by multisets on $\N^q_*$-vectors, resulting in scalar quantities in~\cite{Kin78,Kin78b} being here vectorial;
\item since the total order on~$\N_1$ (partition sizes) underlying the reasoning in~\cite{Kin78} must be replaced by the shortlex order~\eqref{eq:ShortlexOrder} on~$\N_*^q$;
\item because it provides a reference for proofs left to the reader in the monochromatic case in~\cite{Kin78,Kin78b}.
\end{itemize*}

\subsubsection{A compact topology on the polychromatic Kingman simplex}\label{sss:MartinCompactness}
Recall the notation established in the beginning of \S\ref{ss:KingmanNotation}.
We start by recording why the product topology on~$\KSimplex{q}$ is not suitable for our purposes.

\begin{rem}
\begin{enumerate*}[$(1)$]
\item The space~$(\KSimplex{q},\T_\mrmp)$ is not compact (in fact not closed in~$\mbbI_q$), since any sequence~$\seq{\mbfx_\sym{n}}_n\subset \KSimplex{q}$ with~$\length{\mbfx_{\sym{n},0}}=0$ and~$\length{\mbfx_{\sym{n},i}}=1/n$ for all~$i\leq n$ (which forces~$\mbfx_{\sym{n},i}=\zero$ for~$i>n$) converges to the zero sequence in~$\mbbI_q$.
\item Furthermore, it is possible to show that~$\mbfx\mapsto \psi_A(\mbfx)$ does not extend continuously to the $\T_\mrmp$-closure of~$\KSimplex{q}$ %
whenever~$A$ contains non-singleton blocks.
In the case~$q=1$ this was already noted by Kingman in~\cite[after Eqn.~(A9)]{Kin77}.
\item Finally, it is possible to show that the shortlex reordering~$\sort_\shortlex$ too is not continuous.
\end{enumerate*}
\end{rem}

It is nonetheless possible to endow~$\KSimplex{q}$ with a metrizable compact topology, as follows.

The $q$-dimensional \emph{corner simplex} is
\[
\CSimpl^q\eqdef \set{\mbft\in [0,1]^q : \length{\mbft}\leq 1}\comma \qquad \CSimpl^q_*\eqdef \CSimpl^q\setminus\set{0} \fstop
\]

Let~$\msS\eqdef \Delta^{q-1}\times \msP(\CSimpl^q)$ be endowed with the product topology inherited from the Euclidean topology on~$\Delta^{q-1}$ and the weak topology on~$\msP(\CSimpl^q)$, which turns~$\msS$ into a compact Polish space.
Further let~$\msS_*$ be the topological subspace of all $(\DustV,\rho)\in \msS$ satisfying
\begin{enumerate}[$({\mrmS}_1)$]
\item\label{i:Sstar:1} $\Bulk$ is purely atomic and~$\Bulk\!\set{\mbft}\in\length{\mbft}\N_0$ for every~$\mbft\in\CSimpl^q_*$;
\item\label{i:Sstar:2} $\DustV \geq_\had \BulkDustMap{\Bulk}\eqdef \int_{\CSimpl^q_*} \frac{\mbft}{\length{\mbft}}\diff\Bulk(\mbft)$.
\end{enumerate}

\begin{lem}\label{l:IsoSimplex}
The map~$\IsoSimplex \colon \KSimplex{q} \to \msS$ defined by~$\IsoSimplex\eqdef \DustV\otimes\IsoSimplex_2$ with
\begin{equation}\label{eq:IsoSimplex}
\DustV \colon \mbfx \longmapsto \DustV(\mbfx)\eqdef \sum_{i=0}^\infty\mbfx_i\comma \qquad \IsoSimplex_2\colon\mbfx\longmapsto \Bulk_\mbfx\eqdef \length{\mbfx_0}\delta_{\zero} + \sum_{i=1}^\infty \length{\mbfx_i}\delta_{\mbfx_i}
\end{equation}
is a Borel isomorphism of measurable spaces~$\IsoSimplex\colon \KSimplex{q}\to \msS_*$.

\begin{proof}
In light of the ordering in the definition of~$\KSimplex{q}$, the map~$\IsoSimplex$ is injective and thus invertible on its image.
The set~$\IsoSimplex_2(\KSimplex{q})$ exactly consists of those purely atomic~$\Bulk\in\msP(\CSimpl^q)$ satisfying~\ref{i:Sstar:1} above; the block frequencies~$\mbfx_i$, with~$i\geq 1$, are the atoms of~$\Bulk\mrestr{\CSimpl^q_*}$, each listed with multiplicity~$\Bulk\!\set{\mbft}/\length{\mbft}$; and
$%
\mbfx_0= \DustV - \BulkDustMap{\Bulk} \geq_\had\zero %
$, %
which shows that~$\IsoSimplex(\KSimplex{q})=\msS_*$.
Recall from~\S\ref{ss:KingmanNotation} that~$\KSimplex{q}$ is a standard Borel space.
Also,~$\msS_*$ is a Borel subset of the Polish space~$\msS$ and thus itself a standard Borel space.
Since~$\IsoSimplex$ is a bijective Borel map between standard Borel spaces it is a Borel isomorphism.
\end{proof}
\end{lem}

From now on, endow~$\KSimplex{q}$ with the pull-back topology~$\T_*$ via~$\IsoSimplex$ of the topology on~$\msS$.
Since~$\IsoSimplex$ is bi-measurable, the Borel $\sigma$-algebra of~$(\KSimplex{q},\T_\mrmp)$ coincides with the Borel $\sigma$-algebra of~$(\KSimplex{q},\T_*)$.
Thus, there is no ambiguity about Borel measurability on~$\KSimplex{q}$.

\begin{lem}[Compactness of~$\KSimplex{q}$]\label{l:Compactness}
The space~$(\KSimplex{q},\T_*)$ is compact metrizable.
\begin{proof}
Metrizability is trivial since~$\msS$ is compact Polish.
In order to show compactness, it suffices to show that~$\msS_*$ is closed in~$\msS$.
Let~$\seq{(\DustV_n,\Bulk_n)}_n\subset \msS_*$ be converging to~$(\DustV,\Bulk)$.
Let~$\eps\in(0,1)$ be such that~$\Bulk\set{\length{\mbft}=\eps}=0$, which excludes at most countably many values of~$\eps$ since~$\Bulk$ is a finite measure.
In particular, there exists a sequence of such~$\eps$'s tending to~$0$.
Set~$K_\eps\eqdef\set{\mbft \in\CSimpl^q : \length{\mbft}\geq\eps}$, so that~$\Bulk_n\mrestr{K_\eps}\rightharpoonup\Bulk\mrestr{K_\eps}$ by the Portmanteau Theorem.
By~\ref{i:Sstar:1} each~$\Bulk_n\mrestr{K_\eps}$ is purely atomic with at most $\floor{1/\eps}$ atoms, each of mass at least~$\eps$ and of the form~$k\length{\mbft}$ for some~$k\in\N_1$ with~$k\leq1/\eps$.
Along a subsequence the atoms' locations converge in the compact set~$K_\eps$ and the multiplicities~$k$ are eventually constant; atoms converging to a common location contribute masses that add up.
Hence~$\Bulk\mrestr{K_\eps}$ is again purely atomic and each atom~$\mbft$ has mass in~$\length{\mbft}\N_1$. %
Letting~$\eps\to 0$ along an admissible sequence gives~\ref{i:Sstar:1} for~$\Bulk$.
As for~\ref{i:Sstar:2}, the function~$\mbft\mapsto \mbft/\length{\mbft}$ is continuous on~$K_\eps$, and therefore
\begin{equation*}
\int_{K_\eps}\frac{\mbft}{\length{\mbft}}\diff\Bulk(\mbft) = \nlim \int_{K_\eps} \frac{\mbft}{\length{\mbft}}\diff\Bulk_n(\mbft) \leq_\had \liminf_n \BulkDustMap{\Bulk_n} \leq_\had \liminf_n \DustV_n=\DustV \fstop
\end{equation*}
Letting~$\eps\downarrow0$ in the left-most integral,~$\BulkDustMap{\Bulk}\leq_\had \DustV$ follows by Monotone Convergence.
\end{proof}
\end{lem}

For~$\mbfr\in\N^q_*$ and~$\mbfx\in\KSimplex{q}$ set now
\begin{equation}\label{eq:KSimplexCoordinates}
\Pmom{\mbfr}(\mbfx) \eqdef \sum_{i=1}^\infty \mbfx_i^\mbfr = \int_{\CSimpl^q_*} \frac{\mbft^\mbfr}{\length{\mbft}} \diff \Bulk_\mbfx(\mbft) \comma \qquad \Dust_j(\mbfx)\eqdef \sum_{i=0}^\infty x_{i,j}\comma \quad \DustV(\mbfx)= \seq{\Dust_j(\mbfx)}_{j\in [q]} \fstop
\end{equation}
Note that, as soon as~$\length{\mbfr}\geq 2$, the function~$\mbft\mapsto \mbft^\mbfr/\length{\mbft}$ extends continuously from~$\CSimpl^q_*$ to~$\CSimpl^q$ by setting it equal to~$\zero$ at~$\zero$.
Thus,~$\Pmom{\mbfr}$ is continuous on~$\KSimplex{q}$ for every~$\mbfr\in\N^q_*$ with~$\length{\mbfr}\geq 2$, while $\Pmom{\mbfe_j}$~is merely lower semicontinuous for every~$j\in [q]$.
Furthermore,~$\mbfx\mapsto \DustV(\mbfx)$ too is continuous on~$\KSimplex{q}$ by definition of the topology.

\begin{lem} \label{l:poly}
	For every~$m \in \N_0$,~$B \in \Part{q}{m}$, and~$\mbfa \in \N^q_*$, we have
	\begin{equation} \label{eq:l:poly1:00}
		\frac{\mbfa! \, (B(\mbfa)+1)}{(m+\length{\mbfa})!}\varphi_{B+\car_{\mbfa}} = \frac{1}{m!}\Pmom{\mbfa}\,\varphi_B  -  \sum_{\mbfb \in \supp{B}}\! \frac{(\mbfb+\mbfa)! \, (B(\mbfb+\mbfa)+1)}{(m+\length{\mbfa})! \, \mbfb!} \varphi_{B-\car_{\mbfb}+\car_{\mbfb+\mbfa}} \fstop
	\end{equation}
	As a consequence, for every~$n\in\N_0$ and~$A \in \Part{q}{n}$, there exists a polynomial~$\Psi_A$ such that
	\begin{equation} \label{eq:l:poly1:0}
		\varphi_A(\mbfx) = \Psi_A\left(\Pmom{\mbfe_1}(\mbfx), \dotsc, \Pmom{\mbfe_q}(\mbfx), \seq{\Pmom{\mbfr}(\mbfx)}_{2\leq \length{\mbfr}\leq n}\right) \comma \qquad \mbfx \in \KSimplex{q} \fstop
	\end{equation}

\begin{proof}
	Fix~$B \in \Part{q}{\cup}$ and~$\mbfa \in \N^q_*$, define~$A \coloneqq B + \car_{\mbfa}$, let~$(\mbfa_h)_{h \in [\size{A}]}$ be the enumeration of~$A$, let~$\bar h \in [\size{A}]$ be an index such that~$\mbfa_{\bar h}=\mbfa$, let~$H \coloneqq [\size{A}] \setminus \set{\bar h}$. For~$\mbfx \in \KSimplex{q}$,
	\begin{align*}
		\sum_{\iota \in T_{\size{A}}} \prod_{h=1}^{\size{A}} \mbfx_{\iota(h)}^{\mbfa_h} &= \sum_{i=1}^\infty \mbfx_i^\mbfa \sum_{\substack{\iota \in T_{\size{A}} \\ \iota(\bar h) = i}} \prod_{h \in H} \mbfx_{\iota(h)}^{\mbfa_h} = \sum_{i=1}^\infty \mbfx_i^\mbfa \sum_{\substack{\iota \colon H \to \N_1 \setminus \set{i} \\ \text{injective}}} \, \prod_{h \in H} \mbfx_{\iota(h)}^{\mbfa_h} 
		\\
		&=\Pmom{\mbfa}(\mbfx) \sum_{\substack{\iota \colon H \to \N_1 \\ \text{injective}}} \, \prod_{h \in H} \mbfx_{\iota(h)}^{\mbfa_h} - \sum_{k \in H} \, \sum_{\substack{\iota \colon H \to \N_1 \\ \text{injective}}} \mbfx_{\iota(k)}^{\mbfa_{k}+\mbfa} \prod_{ h \in H \setminus \set{k}} \mbfx_{\iota(h)}^{\mbfa_{h}} \\
		&=\Pmom{\mbfa}(\mbfx) \sum_{\substack{\iota \colon H \to \N_1 \\ \text{injective}}} \, \prod_{h \in H} \mbfx_{\iota(h)}^{\mbfa_h} - \sum_{\mbfb \in \supp{B}} \sum_{\substack{k \in H \\ \mbfa_h=\mbfb}} \, \sum_{\substack{\iota \colon H \to \N_1 \\ \text{injective}}} \mbfx_{\iota(k)}^{\mbfb+\mbfa} \prod_{ h \in H \setminus \set{k}} \mbfx_{\iota(h)}^{\mbfa_{h}} \fstop
	\end{align*}
	Up to the constants (independent of~$\mbfx$) in the definition~\eqref{eq:DefPhi}, the left-hand side in the last display is~$\varphi_{A}(\mbfx)$, the first term on the third line is~$\Pmom{\mbfa}(\mbfx) \varphi_B(\mbfx)$, and each of the summands indexed by~$\mbfb$ in the last term is~$\varphi_{B-\car_{\mbfb}+\car_{\mbfb+\mbfa}}(\mbfx)$, each counted with multiplicity~$B(\mbfb)$. Plugging in the precise constants from~\eqref{eq:DefPhi} proves~\eqref{eq:l:poly1:00}.
	
	The second statement is proven by induction on~$\size{A}$.
	The base step~$\size{A} = 0$ is trivial.
	Given~$A$ with~$\size{A} = r+1 \in \N_1$, we pick~$\mbfa \in \supp{A}$ and set~$B \coloneqq A - \car_{\mbfa}$.
	Since~$\size{B} = \size{B-\car_{\mbfb}+\car_{\mbfb+\mbfa}} = r$ for every~$\mbfb \in \supp{B}$, the inductive step follows from~\eqref{eq:l:poly1:00}.
\end{proof}

\end{lem} 

\begin{lem} \label{l:limphi}
	For every~$\mbfx \in \KSimplex{q}$ there exists a sequence~$(\mbfy^{\sym{k}})_{k \in \N_1} \subset \KSimplex{q}$ with~$\mbfy^{\sym{k}}_0 = \zero$ such that
	\begin{equation} \label{eq:l:limphi:00}
		\lim_{k \to \infty} \varphi_A(\mbfy^\sym{k}) = \psi_A(\mbfx) \qquad \text{for every } A \in \Part{q}{\cup} \comma
	\end{equation}
	and
	\begin{equation}\label{eq:l:limphi:01}
		\lim_{k \to \infty} \Pmom{\mbfr}(\mbfy^\sym{k}) = \begin{cases}
			\Dust_j(\mbfx) &\text{if } \mbfr = \mbfe_j \text{~~for some~~} j \in [q]  \\
			\Pmom{\mbfr}(\mbfx) &\text{if } \length{\mbfr} \ge 2 
		\end{cases} \fstop
	\end{equation}

\begin{proof}
	If~$\mbfx_0 = \zero$, we may simply take the constant sequence~$\mbfy^\sym{k} = \mbfx$.
	Let us thus assume that~$\length{\mbfx_0} > 0$.
	For~$k \in \N_1$, we define~$\mbfy^\sym{k} \in \KSimplex{q}$ adding~$k$ blocks~$k^{-1} \mbfx_0$ and replacing~$\mbfx_0$ with~$\zero$. More precisely,
	\[
	\mbfy^\sym{k}_i \coloneqq \begin{cases}
		\zero &\text{if } i=0  \\
		\mbfx_i &\text{if }  i \ge 1 \text{ and } \mbfx_i \geq_\shortlex  k^{-1} \mbfx_0  \\
		\mbfx_{i-k} &\text{if }  i \ge k+1 \text{ and } \mbfx_{i-k} \not\geq_\shortlex k^{-1} \mbfx_0  \\
		k^{-1} \mbfx_0 &\text{otherwise}
	\end{cases} \fstop
	\]
	We let~$I \subset \N_1$ be the set of (consecutive) indices~$i$ for which the fourth case in the latter display holds.
	Since~$\length{\mbfx_0} > 0$ and~$\length{\mbfx_i} \to 0$ as~$i \to \infty$, eventually,~$\mbfx_i \not\geq_\shortlex k^{-1}\mbfx_0$; hence,~$\card{I} = k$.
	
	The limit~\eqref{eq:l:limphi:01} is straightforward.
	Fix~$A \in \Part{q}{\cup}$ with enumeration~$(\mbfa_h)_{h \in [\size{A}]}$ and focus on~\eqref{eq:l:limphi:00}.
	Observe that we may partition~$T_{\size{A}}$ according to the set~$\iota^{-1}(I)$,
	and thus write,
	\begin{align*}
		&\sum_{\iota \in T_{\size{A}}} \prod_{h=1}^{\size{A}} \tparen{\mbfy_{\iota(h)}^\sym{k}}^{\mbfa_{h}} = \sum_{\substack{H \subset [\size{A}] \\ \card{H} \le k}} (k^{-1}\mbfx_0)^{\sum_{h \in H} \mbfa_h} \sum_{\substack{\iota \in T_{\size{A}} \\ \iota^{-1}(I)=H}} \prod_{h \in [\size{A}]\setminus H} \tparen{\mbfy_{\iota(h)}^\sym{k}}^{\mbfa_{h}} \\
		&\quad = \sum_{\substack{H \subset [\size{A}] \\ \card{H} \le k}} (k^{-1}\mbfx_0)^{\sum_{h \in H} \mbfa_h} \frac{k!}{(k-\card{H})!}\sum_{\substack{\iota \colon [\size{A}]\setminus H \to \N_1 \setminus I \\ \text{injective}}} \, \prod_{h \in [\size{A}]\setminus H} \tparen{\mbfy_{\iota(h)}^\sym{k}}^{\mbfa_{h}}
	\end{align*}
	Define~$A_{H^c} \coloneqq \sum_{h \in [\size{A}]\setminus H} \car_{\mbfa_h}$. From the previous formula and the definition of~$\varphi_A$:
	\[ 
	\varphi_A(\mbfy^\sym{k}) = \sum_{\substack{H \subset [\size{A}] \\ \card{H} \le k}}  \frac{k! \, \length{\col{A}}! \, (k^{-1}\mbfx_0)^{\sum_{h \in H} \mbfa_h}}{(k-\card{H})! \, \length{\col{A_{H^c}}}!} \paren{\prod_{\mbfa \in \supp{A}} \frac{\mbfa!^{A_{H^c}(\mbfa)}  A_{H^c}(\mbfa)! }{\mbfa!^{A(\mbfa)} A(\mbfa)! }}   \varphi_{A_{H^c}}(\mbfx) \fstop
	\]
	In the limit~$k \to \infty$, we have~$\frac{k!}{(k-\card{H})!} \sim k^{\card{H}}$. 
	Therefore, the only non-negligible terms in the sum above are those for which~$\length{\mbfa_h} = 1$ for every~$h \in H$, i.e.,~$A_{H^c} = A-\uni{\mbfm}$ for some~$\mbfm \in \N_0^q$. We find
	\[
	\lim_{k \to \infty} \varphi_A(\mbfy^\sym{k}) = \sum_{\substack{\mbfm \in \N_0^q \\ \uni{\mbfm} \le A}} c_{\mbfm} \frac{\length{\col{A}}! \, \mbfx_0^{\mbfm}}{(\length{\col{A}}-\length{\mbfm})!} \paren{\prod_{j=1}^q \frac{(A(\mbfe_j)-m_j)! }{A(\mbfe_j)!} }   \varphi_{A-\uni{\mbfm}}(\mbfx) \comma
	\]
	where~$c_{\mbfm} = \prod_{j=1}^q \binom{A(\mbfe_j)}{m_j}$ is the number of sets~$H \subset [\size{A}]$ such that~$A_{H^c} = A - \uni{\mbfm}$. The conclusion follows from the definition of~$\psi_A$.
\end{proof}

\end{lem}

\begin{cor} \label{cor:poly2}
	For every~$n\in\N_0$ and~$A \in \Part{q}{n}$
	\begin{equation} \label{eq:cor:poly2:0}
		\psi_A(\mbfx) = \Psi_A\paren{\Dust_1(\mbfx), \dotsc, \Dust_q(\mbfx), \seq{\Pmom{\mbfr}(\mbfx)}_{2\leq \length{\mbfr}\leq n} } \comma \qquad \mbfx \in \KSimplex{q} \comma
	\end{equation}
	where~$\Psi_A$ is the polynomial in~\eqref{eq:l:poly1:0}.

\begin{proof}
	This is a direct consequence of Lemmas~\ref{l:poly} and~\ref{l:limphi}.
\end{proof}
\end{cor}

\begin{cor}\label{c:MeasureDetermining}
The function~$\psi_A$ is continuous on~$(\KSimplex{q},\T_*)$ for every~$A\in\Part{q}{\cup}$.
Furthermore, the linear span~$\msL$ of~$\set{\psi_A}_{A\in\Part{q}{\cup}}$ is a unital point-separating sub-algebra of~$\Cb(\KSimplex{q},\T_*)$,  dense in~$\Cb(\KSimplex{q},\T_*)$.
In particular,~$\set{\psi_A}_{A\in \Part{q}{\cup}}$ is a measure-determining class on~$\KSimplex{q}$ and~$\T_*$ is the coarsest topology making all functions in~$\msL$ continuous.

\begin{proof}
	Recall that~$\Pmom{\mbfr}$ is continuous for every~$\mbfr \in \N^q_*$ with~$\length{\mbfr} \ge 2$, and~$\Dust_j$ is continuous for every~$j \in [q]$. Therefore,~$\psi_A$ is continuous for every~$A \in \Part{q}{\cup}$ by Corollary~\ref{cor:poly2}. It is also clear that~$\msL$ contains the constants, since~$\psi_0 \equiv 1$.
	
	\paragraph{Claim:~$\msL$ is an algebra}
	The linear span~$\msL_\varphi$ of the~$\varphi_A$'s is closed under multiplication by~$\Pmom{\mbfa}$ for all~$\mbfa$'s in $\N^q_*$ by~\eqref{eq:l:poly1:00}. In light of~\eqref{eq:l:poly1:0}, it follows that~$\msL_\varphi$ is an algebra. With~\eqref{eq:l:limphi:00} in Lemma~\ref{l:limphi}, we extend the result to~$\msL$.

	\paragraph{Claim:~$\msL$ is point-separating} Fix~$\mbfx,\mbfy \in \KSimplex{q}$ and assume that~$\psi_A(\mbfx) = \psi_A(\mbfy)$ for every~$A \in \Part{q}{\cup}$. For every~$\mbfr \in \N_*^q$, since~$\length{\mbfr+\mbfe_j} \ge 2$ for every~$j \in [q]$, we have
	\[
		\int_{\CSimpl^q} \mbft^\mbfr \diff \Bulk_\mbfx(\mbft) = \sum_{j=1}^q \int_{\CSimpl^q_*} \frac{\mbft^{\mbfr+\mbfe_j}}{\length{\mbft}} \diff \Bulk_\mbfx(\mbft) = \sum_{j=1}^q \Pmom{\mbfr+\mbfe_j}(\mbfx) = \sum_{j=1}^q \binom{\length{\mbfr}+1}{\mbfr+\mbfe_j}^{-1} \psi_{\car_{\mbfr+\mbfe_j}}(\mbfx) \fstop
	\]
	It follows from this computation and our assumption that~$\int_{\CSimpl^q} \mbft^\mbfr \diff \Bulk_\mbfx(\mbft) = \int_{\CSimpl^q} \mbft^\mbfr \diff \Bulk_\mbfy(\mbft)$ for every~$\mbfr \in \N_*^q$. Since~$\CSimpl^q$ is compact, and polynomials form a point-separating sub-algebra of~$\Cb(\CSimpl^q)$, we infer that~$\Bulk_\mbfx = \Bulk_\mbfy$. Furthermore, we have
	\[
		\Dust_j(\mbfx) = \psi_{\car_{\mbfe_j}}(\mbfx) = \psi_{\car_{\mbfe_j}}(\mbfy) = \Dust_j(\mbfy)\comma \qquad j \in [q] \comma
	\]
	and we conclude that~$\mbfx = \mbfy$ by Lemma~\ref{l:IsoSimplex}.
	
	\paragraph{Conclusion}
	$\msL$ is dense in~$\Cb(\KSimplex{q},\T_*)$ by Stone--Weierstra\ss\ Theorem.
	Since~$\KSimplex{q}$ is metrizable, it is completely regular Hausdorff and thus~$\Cb(\KSimplex{q},\T_*)$ (therefore~$\msL$) identifies~$\T_*$.
\end{proof}
\end{cor}

\subsubsection{Consistency of extremal structures}\label{sss:AuxCombinatorics}
In order to prove the polychromatic Kingman Representation Theorem, let us show that the family~$\psi_\proc(\mbfx)$ is indeed a $\CPS{q}$.%

\begin{lem}\label{l:PhiConsistency}
For every~$q\in\N_1$, every~$m,n \in \N_0$ with $m \leq n$, and every~$\mbfx\in\KSimplex{q}$,
\begin{equation}\label{eq:l:PhiConsistency2:0}
	(1-\length{\mbfx_0})^{n-m} \varphi_B(\mbfx) = \sum_{A\in \Part{q}{n}} \varphi_A(\mbfx)\, \BackwardK^\sym{q}_{m,n}(A,B) \comma \qquad B\in\Part{q}{m} \qquad (0^0\eqdef 1)\fstop%
\end{equation}

\begin{proof}
It suffices to verify~\eqref{eq:l:PhiConsistency2:0} with~$m=n-1$. Analogously to~\eqref{eq:ConsistencySpecial}, this reduces to
\begin{equation}\label{eq:ConsistencyPhiSingleStep}
\begin{aligned}
(1-\length{\mbfx_0}) \, n\, \varphi_B(\mbfx) &= \sum_{j=1}^q \tparen{B(\mbfe_j)+1} \varphi_{B+\car_{\mbfe_j}}(\mbfx)
\\
&\qquad + \sum_{j=1}^q \sum_{\mbfb\in\supp B} (b_j+1) \tparen{B(\mbfb+\mbfe_j)+1} \varphi_{B-\car_\mbfb+\car_{\mbfb+\mbfe_j}}(\mbfx) \fstop
\end{aligned}
\end{equation}
The latter can be easily derived by summing~\eqref{eq:l:poly1:00} with~$\mbfa = \mbfe_j$ over~$j \in [q]$ and observing that~$\sum_{j=1}^q P_{\mbfe_j}(\mbfx) = 1-\length{\mbfx_0}$.
\end{proof}
\end{lem}

\begin{cor}\label{l:PsiConsistency}
	For every~$q\in\N_1$, every~$m,n \in \N_0$ with~$m\leq n$, and every~$\mbfx\in\KSimplex{q}$,
	\begin{equation}\label{eq:l:PsiConsistency:0}
		\psi_B(\mbfx) = \sum_{A\in \Part{q}{n}} \psi_A(\mbfx)\, \BackwardK^\sym{q}_{m,n}(A,B) \comma \qquad B\in\Part{q}{m} \fstop
	\end{equation}

	\begin{proof}
	This is a direct consequence of Lemmas~\ref{l:PhiConsistency} and~\ref{l:limphi}.
	\end{proof}
	
\end{cor}

\subsubsection{Potential theory for the polychromatic Hoppe chain}\label{sss:PHCRevisited}
Throughout this section, we fix~$\theta>0$ and~$\nu\in\msP(\Delta^{q-1})$ with~$\nu\mathring{\Delta}^{q-1}>0$.
Recall the construction of the PHC~$\msA_\proc$ in~\S\ref{ss:Hoppe}. %
For~$A\in\Part{q}{\cup}$ further let~$\tau^A\eqdef\min\set{t\in\N_0 : \msA_t=A}$, with~$\min\emp\eqdef\infty$.
The \emph{potential kernel} of~$\msA_\proc$ is
\[
\PotK^\sym{q}(B,A)\eqdef \sum_{t=0}^\infty \probP^B[\tau^A=t] \fstop
\]
Since~$\length{\col{\msA_{t+1}}}=\length{\col{\msA_t}}+1$, the chain~$\msA_\proc$ is transient and every state is visited at most once.
In particular, its Green kernel coincides with its potential kernel, so that we may restrict ourselves to only use the latter.
For notational convenience, in the following we will occasionally relabel~$\ForwardK^\sym{q}_t\eqdef \ForwardK^\sym{q}_{m,n}$ whenever~$n-m=t$.

\begin{lem}[Potential kernel]\label{l:Green}
For~$m,n\in\N_0$, and every~$B\in\Part{q}{m}$ and~$A\in\Part{q}{n}$,
\begin{equation}\label{eq:l:Green:0}
\PotK^\sym{q}(B,A)= \sum_{t=0}^\infty \ForwardK^\sym{q}_t(B,A) = \begin{cases} \displaystyle \BackwardK^\sym{q}_{m,n}(A,B)\, \frac{\pEsf{n}{\theta,\nu}[A]}{\pEsf{m}{\theta,\nu}[B]} & \text{if } m\leq n \\[.6em] 0 & \text{otherwise}\end{cases} \fstop
\end{equation}

\begin{proof}
Since each step increases~$\length{\col{\msA_t}}$ by one, $A$ is accessible from~$B$ only if~$m\leq n$, in which case it can be reached only at time~$n-m$ and at most once, so that~$\PotK^\sym{q}(B,A)=\ForwardK^\sym{q}_{n-m}(B,A)$.
We argue by induction on~$n-m$.
If~$n=m$, then~$\ForwardK^\sym{q}_0(B,A)=\car_{A=B}=\BackwardK^\sym{q}_{n,n}(A,B)$ by~\eqref{eq:Kernels}, which is~\eqref{eq:l:Green:0}.
Assume~\eqref{eq:l:Green:0} for~$n-1-m$.
By the cocycle relation~\eqref{eq:CocycleForward}, by the inductive hypothesis, by Lemma~\ref{l:Reversal}, and by the cocycle relation~\eqref{eq:Cocycle},
\begin{align*}
\ForwardK^\sym{q}_{n-m}(B,A) &= \sum_{C\in\Part{q}{n-1}} \ForwardK^\sym{q}_{n-1-m}(B,C)\, \ForwardK^\sym{q}_{n-1,n}(C,A)
\\
&= \sum_{C\in\Part{q}{n-1}} \BackwardK^\sym{q}_{m,n-1}(C,B) \frac{\pEsf{n-1}{\theta,\nu}[C]}{\pEsf{m}{\theta,\nu}[B]} \, \BackwardK^\sym{q}_{n-1,n}(A,C) \frac{\pEsf{n}{\theta,\nu}[A]}{\pEsf{n-1}{\theta,\nu}[C]}
\\
&= \frac{\pEsf{n}{\theta,\nu}[A]}{\pEsf{m}{\theta,\nu}[B]} \sum_{C\in\Part{q}{n-1}} \BackwardK^\sym{q}_{m,n-1}(C,B)\, \BackwardK^\sym{q}_{n-1,n}(A,C) = \frac{\pEsf{n}{\theta,\nu}[A]}{\pEsf{m}{\theta,\nu}[B]}\, \BackwardK^\sym{q}_{m,n}(A,B)\comma
\end{align*}
which concludes the proof. (All the denominators are non-zero by~\eqref{eq:PositivePESF}.)
\end{proof}
\end{lem}

\begin{defs}[Martin kernel]\label{d:MartinKernelPartitions}
The \emph{Martin kernel} (\emph{based at~$0\in \Part{q}{0}$}) of the PHC is
\begin{equation}\label{eq:d:MartinKernelPartitions:0}
\MartinK^\sym{q}(B,A)\eqdef \frac{\PotK^\sym{q}(B,A)}{\PotK^\sym{q}(0,A)} = \begin{cases} \displaystyle \frac{\BackwardK^\sym{q}_{m,n}(A,B)}{\pEsf{m}{\theta,\nu}[B]} & \text{if } m\leq n \\[.6em] 0 & \text{otherwise}\end{cases}\comma \qquad B\in \Part{q}{m}\comma A\in\Part{q}{n}\fstop
\end{equation}
\end{defs}

\begin{rem}[Weak irreducibility]\label{r:WeakIrreducibilityPartitions}
As is the case for~$\DRG{\mbfp}$ in~\S\ref{sss:DRG}, the chain~$\msA_\proc$ is not irreducible: the accessible component of~$B\in\Part{q}{m}$ is~$\tset{A: \BackwardK^\sym{q}_{m,\length{\col{A}}}(A,B)>0}$.
It is however a \emph{weakly irreducible} combinatorial Markov chain in the sense of~\cite[\S3]{GerGruHag16}, since by~\eqref{eq:PositivePESF} and Lemma~\ref{l:Green} we have~$\PotK^\sym{q}(0,A)=\pEsf{n}{\theta,\nu}[A]>0$ for every~$A$, i.e.\ every state is accessible from~$0$, which grants that the Martin kernel based at~$0$ is everywhere \emph{finite}.
\end{rem}

\begin{rem}[Independence of the reference chain]\label{r:IndependenceReference}
By~\eqref{eq:d:MartinKernelPartitions:0}, the Martin kernel of the PHC differs from the deletion kernel~$\BackwardK^\sym{q}$ only by the factor~$\pEsf{m}{\theta,\nu}[B]^{-1}$, which is positive and depends on~$B$ alone.
Consequently, for a sequence~$\seq{A_k}_k\subset\Part{q}{\cup}$, the pointwise convergence of~$\MartinK^\sym{q}(\emparg,A_k)$ is equivalent to that of~$\BackwardK^\sym{q}_{\emparg,\length{\col{A_k}}}(A_k,\emparg)$, and the Martin compactification of~$\Part{q}{\cup}$ constructed below, as well as the notion of convergence to the boundary, is independent of the choice of the parameters~$(\theta,\nu)$ satisfying~\eqref{eq:StandingAssumption}.
Considering~$\msA_\proc$ allows us to turn the co-transition kernels~\eqref{eq:Kernels} into a Markov chain, as opposed to the \emph{reverse} Markov chain~$\hat\msA_\proc$.
\end{rem}

We denote by~$\mcH^\tgeq_{\theta,\nu}$ the space of all non-negative \emph{$\ForwardK^\sym{q}$-harmonic functions} on~$\Part{q}{\cup}$, i.e.\ of all~$h\colon\Part{q}{\cup}\to[0,\infty)$ with
\begin{equation}\label{eq:HarmonicPartitions}
\sum_{A\in\Part{q}{m+1}} \ForwardK^\sym{q}_{m,m+1}(B,A)\, h(A) = h(B)\comma \qquad B\in\Part{q}{m}\comma \quad m\in\N_0\comma
\end{equation}
and we call \emph{minimal harmonic} any~$h_0\in \mcH^\tgeq_{\theta,\nu}$ with~$h_0(0)=1$ and such that~$h_0\geq h\in \mcH^\tgeq_{\theta,\nu}$ implies~$h=\tau h_0$ for some~$\tau\in[0,1]$.

\begin{lem}[Harmonic functions and~$\CPS{q}$'s]\label{l:HarmonicVsCPS}
The assignment
\begin{equation}\label{eq:l:HarmonicVsCPS:0}
h\longmapsto \probP^h_\proc \comma \qquad \probP^h_n[A]\eqdef h(A)\, \pEsf{n}{\theta,\nu}[A]\comma \qquad A\in\Part{q}{n}\comma\quad n \in\N_0 \comma
\end{equation}
is an affine bijection from~$\tset{h\in \mcH^\tgeq_{\theta,\nu}: h(0)=1}$ onto~$\CPS{q}$.
Under this bijection,~$\probP^h_\proc$ is the $\BackwardK^\sym{q}$-harmonic density corresponding to~$h$.
In particular, for every~$\mbfx\in\KSimplex{q}$, the function
\begin{equation}\label{eq:l:HarmonicVsCPS:1}
h_\mbfx\colon A\longmapsto \frac{\psi_A(\mbfx)}{\pEsf{n}{\theta,\nu}[A]}\comma \qquad A \in \Part{q}{n}\comma \quad n\in\N_0\comma
\end{equation}
belongs to~$\mcH^\tgeq_{\theta,\nu}$ and satisfies~$h_\mbfx(0)=1$.

\begin{proof}
Let~$h\colon \Part{q}{\cup}\to [0,\infty)$ and~$\probP^h_\proc$ be as in~\eqref{eq:l:HarmonicVsCPS:0}.
By Lemma~\ref{l:Reversal}, for every~$B\in\Part{q}{m}$,
\begin{align*}
\sum_{A\in\Part{q}{m+1}} \ForwardK^\sym{q}_{m,m+1}(B,A) h(A) &= \frac{1}{\pEsf{m}{\theta,\nu}[B]} \sum_{A\in\Part{q}{m+1}} \BackwardK^\sym{q}_{m,m+1}(A,B)\, \probP^h_{m+1}[A] 
\\
&= \frac{\tparen{\sigma^\sym{q}_{m,m+1}\probP^h_{m+1}}[B]}{\pEsf{m}{\theta,\nu}[B]} \comma
\end{align*}
cf.~\eqref{eq:Consistency}.
Thus~\eqref{eq:HarmonicPartitions} holds if and only if~$\probP^h_m=\sigma^\sym{q}_{m,m+1}\probP^h_{m+1}$ for every~$m\in\N_0$, which in turn is equivalent to~\eqref{eq:d:Consistency:0} by the cocycle relation~\eqref{eq:CocycleSigma}.
Moreover, summing the consistency identity~$\probP^h_m[B]=\tparen{\sigma^\sym{q}_{m,n}\probP^h_{n}}[B]$ over~$B\in\Part{q}{m}$ and using that~$\BackwardK^\sym{q}_{m,n}(A,\emparg)$ is a probability shows that~$m\mapsto \sum_{B\in\Part{q}{m}}\probP^h_m[B]$ is constant, hence identically equal to~$\probP_0^h[0]=h(0)$.
Therefore~$h\in\mcH^\tgeq_{\theta,\nu}$ with~$h(0)=1$ if and only if~$\probP^h_\proc\in\CPS{q}$.
Since~$\pEsf{n}{\theta,\nu}[A]>0$ for every~$A\in\Part{q}{\cup}$ by~\eqref{eq:PositivePESF}, the assignment~\eqref{eq:l:HarmonicVsCPS:0} is injective and onto, and it is clearly affine.
Finally,~$\psi_\proc(\mbfx)$ is a probability density on~$\Part{q}{n}$ for every~$n$ by Lemma~\ref{l:TotalPhi}, and it is $\BackwardK^\sym{q}$-harmonic by Corollary~\ref{l:PsiConsistency}; that is,~$\psi_\proc(\mbfx)\in\CPS{q}$, and the last assertion follows by the first part of the proof since~$\psi_0\equiv 1$.
\end{proof}
\end{lem}

\subsubsection{Convergence of the Martin kernel}\label{sss:MartinConvergence}
We start by recording the (deterministic) map assigning to a polychromatic partition its empirical block-color frequencies.

\begin{defs}[Empirical frequencies]\label{d:EmpiricalMap}
For~$A\in\Part{q}{n}$ with~$n\geq 1$ and enumeration~$\seq{\mbfa_h}_{h\in[\size{A}]}$ as in Definition~\ref{d:Enumeration}, let
\begin{equation}\label{eq:d:EmpiricalMap:0}
\mbfX(A)\eqdef \seq{\mbfX(A)_i}_{i\geq 0} \comma \qquad \mbfX(A)_0\eqdef \zero\comma \qquad \mbfX(A)_i\eqdef\begin{cases} \mbfa_i/n & i\in [ \size{A}] \\ \zero & i>\size{A}\end{cases} \fstop
\end{equation}
\end{defs}

Since~$\mbfa_i\geq_\shortlex \mbfa_{i+1}$ and~$\sum_{i}\length{\mbfa_i}=\length{\col{A}}=n$, we have~$\mbfX(A)\in\KSimplex{q}$.

\begin{lem}[Sampling with and without replacement]\label{l:SamplingEstimate}
For~$m\in\N_0$,~$n \in \N_1$ with~$m\leq n$,
\begin{equation}\label{eq:l:SamplingEstimate:0}
\abs{\BackwardK^\sym{q}_{m,n}(A,B)-\psi_B\tparen{\mbfX(A)}} \leq 1- \frac{\AntiPoch{n}{m}}{n^m} %
\comma\qquad A\in\Part{q}{n}\comma B\in\Part{q}{m}\fstop
\end{equation}
Furthermore,
\begin{equation}\label{eq:l:SamplingEstimate:1}
\sum_{B\in\Part{q}{m}} \abs{\BackwardK^\sym{q}_{m,n}(A,B)-\psi_B\tparen{\mbfX(A)}} \leq 2\paren{1- \frac{\AntiPoch{n}{m}}{n^m}} \comma \qquad A \in\Part{q}{n} \fstop
\end{equation}

\begin{proof}
Set~$\mbfn\eqdef\col{A}$ and fix an $\mbfn$-coloring~$\coloring_\mbfn$ as in Definition~\ref{d:Category}.
By Lemma~\ref{l:MultiA}, we have~$\card{\colore^{-1}(A)}=\Multi{\mbfn}(A)\geq 1$, hence there exists~$\pi\in\mfS_n$ with~$\colore(\pi)=A$; let~$\Lambda\eqdef \Lambda_\pi$ be the corresponding set partition~\eqref{eq:StdSetPartitionPi} of~$[n]$, with blocks~$\seq{L_h}_{h\in [\size{A}]}$ labeled so that~$L_h$ has color count~$\mbfa_h$.

Let~$\seq{V_\ell}_{\ell\in [m]}$ be a uniform ordered sample \emph{without} replacement from~$[n]$, and let $\seq{W_\ell}_{\ell\in [m]}$ be i.i.d.\ uniform on~$[n]$, i.e.\ a sample \emph{with} replacement.
By Remark~\ref{r:KernelInterpretation}, the polychromatic partition induced by~$(\Lambda,\coloring_\mbfn)$ on~$\set{V_\ell}_{\ell\in[m]}$ has law~$\BackwardK^\sym{q}_{m,n}(A,\emparg)$.
On the other hand, a single draw~$W_\ell$ belongs to the block~$L_i$ and has color~$j$ with probability~$a_{i,j}/n=\mbfX(A)_{i,j}$.
Thus, denoting by~$\boldnu_i$ the vector of the colors of the draws falling into~$L_i$, the probability of observing a given sequence~$\boldnu=\seq{\boldnu_i}_{i\geq 1}$ as in~\eqref{eq:l:TotalPhi:1.2} is~$m!\prod_{i\geq 1}\mbfX(A)_i^{\boldnu_i}/\boldnu_i!$.
By~\eqref{eq:alterPhi}, summing over all~$\boldnu$'s with~$A_\boldnu=B$, cf.~\eqref{eq:l:TotalPhi:1.5}, shows that the polychromatic partition induced by~$(\Lambda,\coloring_\mbfn)$ on the sample~$\seq{W_\ell}_{\ell\in [m]}$ has law~$\varphi_\proc\tparen{\mbfX(A)}=\psi_\proc\tparen{\mbfX(A)}$, where the last equality holds since~$\mbfX(A)_0=\zero$.

Conditionally on the event~$E\eqdef\set{W_\ell \text{ pairwise distinct}}$, of probability~$\probP[E]=\AntiPoch{n}{m} /n^m$, the sample~$\seq{W_\ell}_\ell$ is a uniform ordered sample without replacement.
The %
samples~$\seq{W_\ell}_\ell$ and~$\seq{V_\ell}_\ell$ may therefore be coupled so as to coincide on~$E$. %
Hence,~$\abs{\BackwardK^\sym{q}_{m,n}(A,B)-\psi_B\tparen{\mbfX(A)}}\leq \probP[E^\complement]$ for every~$B\in\Part{q}{m}$, which is~\eqref{eq:l:SamplingEstimate:0}, and the total-variation distance between the two laws is at most~$2\,\probP[E^\complement]$, which proves~\eqref{eq:l:SamplingEstimate:1}.
\end{proof}
\end{lem}

\begin{lem}[Density of the empirical frequencies]\label{l:DensityEmpirical}
For every~$\mbfx\in\KSimplex{q}$ there exists a sequence~$\seq{A_k}_k\subset \Part{q}{\cup}$ with~$\length{\col{A_k}}\to\infty$ and~$\mbfX(A_k)\to\mbfx$ in~$(\KSimplex{q},\T_*)$.
\end{lem}

\begin{proof}
Since~$(\KSimplex{q},\T_*)$ is compact metrizable by Lemma~\ref{l:Compactness}, it suffices, by a diagonal argument, to prove the assertion for~$\mbfx$ in a $\T_*$-dense subset of~$\KSimplex{q}$. 
We proceed in steps: in Step~1 we reduce ourselves to the subset of~$\mbfx$'s with~$\mbfx_i=\zero$ eventually as~$i\to\infty$;
in Step~2 we reduce further to the subset of~$\mbfx$'s with rational coordinates; and in Step~3 we construct polychromatic partitions with frequencies in the latter set.
Recall that, by definition of~$\T_*$ and by Lemma~\ref{l:IsoSimplex}, a sequence~$\seq{\mbfy_\sym{k}}_k$ converges to~$\mbfy$ in~$\T_*$ if and only if~$\DustV(\mbfy_\sym{k})\to \DustV(\mbfy)$ in~$\Delta^{q-1}$ and~$\Bulk_{\mbfy_\sym{k}}\rightharpoonup\Bulk_{\mbfy}$ in~$\msP(\CSimpl^q)$, cf.~\eqref{eq:IsoSimplex}.

\paragraph{Step~1: Truncation}
For~$k\in\N_1$, let~$\mbfy_\sym{k}\in\KSimplex{q}$ be defined by
\[
\mbfy_{\sym{k},0}\eqdef \mbfx_0+\sum_{i>k}\mbfx_i \comma \qquad \mbfy_{\sym{k},i}\eqdef \mbfx_i \quad \text{for}\quad 1\leq i \leq k\comma \qquad \mbfy_{\sym{k},i}\eqdef\zero\quad \text{for} \quad i>k\fstop
\]
Then~$\DustV(\mbfy_\sym{k})=\DustV(\mbfx)$ for every~$k$, while, for every~$f\in\Cb(\CSimpl^q)$,
\[
\abs{\int_{\CSimpl^q} f\diff \Bulk_{\mbfy_\sym{k}} - \int_{\CSimpl^q} f \diff \Bulk_\mbfx} = \abs{\sum_{i>k} \length{\mbfx_i}\tparen{f(\zero)-f(\mbfx_i)}} \leq 2 \norm{f}_\infty \sum_{i>k} \length{x_i} \xrightarrow{k\to\infty} 0 \fstop
\]
Thus~$\mbfy_\sym{k}\to\mbfx$ in~$\T_*$, and we may and shall assume that~$\mbfx_i=\zero$ for all~$i>k$, for some~$k\in\N_1$.

\paragraph{Step~2: Rationalization}
Let~$\mbfx$ be as at the end of Step~1.
For~$m\in\N_1$ set~$\mbfc_i\eqdef \seq{c_{i,j}}_{j\in [q]}$~with
\[
c_{i,j}\eqdef \begin{cases} \floor{m x_{0,1}} + m-\sum_{i=0}^k\sum_{j=1}^q \floor{m x_{i,j}} & i=0\comma j=1
\\
\floor{m x_{i,j}} & \text{otherwise} 
\end{cases}
\comma
\]
so that
\[
\sum_{i=0}^k\length{\mbfc_i}=m \comma \qquad \abs{\mbfc_i/m-\mbfx_i}\leq (k+1)q/m \comma \qquad 0\leq i \leq k \fstop
\]
Let~$\mbfz_\sym{m}$ be obtained by leaving~$\mbfc_0/m$ in place, by reordering~$\seq{\mbfc_i/m}_{i\in [k]}$ decreasingly in the shortlex order, and by completing with~$\zero$'s, and note that~$\mbfz_\sym{m}\in\KSimplex{q}$ for every~$m\in\N_1$.
Since
\[
\IsoSimplex\colon \seq{\mbfx_i}_{0\leq i\leq k}\longmapsto \tparen{\DustV(\mbfx),\Bulk_\mbfx}=\paren{\sum_{i=0}^k \mbfx_i\, ,\ \length{\mbfx_0}\delta_{\zero}+\sum_{i=1}^k\length{\mbfx_i}\delta_{\mbfx_i}}
\]
is continuous, %
and since~$\DustV$ and~$\Bulk$ are invariant under reordering of~$\seq{\mbfx_i}_{i\geq 1}$, we may compute~$\tparen{\DustV(\mbfz_\sym{m}),\Bulk_{\mbfz_\sym{m}}}$ from the unordered family~$\seq{\mbfc_i/m}_{0\leq i \leq k}$ and conclude that~$\mbfz_\sym{m}\to\mbfx$ in~$\T_*$ as~$m\to\infty$.
In the following, we may thus additionally assume that~$x_{i,j}=c_{i,j}/m$ for some~$m\in\N_1$ and~$c_{i,j}\in\N_0$ with~$\sum_{i=0}^k \length{\mbfc_i}=m$.

\paragraph{Step~3: Realization}
Let~$\mbfx$ be as at the end of Step~2, and set, for~$\ell\in\N_1$ and~$n\eqdef \ell m$,
\[
A_n\eqdef \sum_{i\in[k]:\, \mbfc_i\neq \zero} \car_{\ell\mbfc_i} + \sum_{j=1}^q \ell\, c_{0,j}\, \car_{\mbfe_j} \in \Part{q}{n} \fstop
\]
The non-zero entries of~$\mbfX(A_n)$ are~$\ell\mbfc_i/n=\mbfx_i$, for~$i\in [k]$ with~$\mbfc_i\neq\zero$, together with%
~$\ell c_{0,j}$ entries equal to~$\mbfe_j/n$, for each~$j\in[q]$.
Hence
\[
\DustV\tparen{\mbfX(A_n)} = \sum_{i=1}^k \mbfx_i + \sum_{j=1}^q \frac{\ell c_{0,j}}{n}\mbfe_j = \sum_{i=1}^k \mbfx_i +\mbfx_0 = \DustV(\mbfx) \comma
\]
while
\[
\Bulk_{\mbfX(A_n)} = \sum_{i=1}^k \length{\mbfx_i}\delta_{\mbfx_i} + \sum_{j=1}^q x_{0,j}\, \delta_{\mbfe_j/n} \xrightharpoonup{\ \ell\to\infty\ } \sum_{i=1}^k \length{\mbfx_i}\delta_{\mbfx_i} + \length{\mbfx_0}\delta_{\zero} = \Bulk_\mbfx
\]
narrowly, since~$\mbfe_j/n\to\zero$.
Thus~$\mbfX(A_{\ell m})\to\mbfx$ in~$\T_*$ as~$\ell\to\infty$, which concludes the proof.
\end{proof}

As in~\S\ref{sss:PHCRevisited}, let us fix~$\theta>0$ and~$\nu\in\msP(\Delta^{q-1})$ with~$\nu\mathring{\Delta}^{q-1}>0$.
We say that a sequence~$\seq{A_k}_k\subset\Part{q}{\cup}$ \emph{converges to the Martin boundary} if~$A_k\to\infty$ (i.e.~$\seq{A_k}_k$ eventually leaves every finite subset of~$\Part{q}{\cup}$) and~$\MartinK^\sym{q}(\emparg,A_k)$ converges pointwise on~$\Part{q}{\cup}$.
(Note that~$A_k\to\infty$ if and only if~$\length{\col{A_k}}\to\infty$, since~$\Part{q}{n}$ is a finite set for every~$n$.) %
Two such sequences are \emph{equivalent} if the corresponding kernel limits coincide at every point of~$\Part{q}{\cup}$.
The \emph{Martin boundary}~$\msM$ of the PHC is the set of the corresponding equivalence classes, identified with the set of the corresponding limit functions and endowed with the topology of pointwise convergence.

For each~$\mbfx\in\KSimplex{q}$ we let~$\MartinK^\sym{q}(\emparg,\mbfx)$ be the \emph{(extended) Martin kernel}, cf.~\eqref{eq:l:HarmonicVsCPS:1},
\begin{equation}\label{eq:ExtendedMartinKernel}
\MartinK^\sym{q}(\emparg,\mbfx)\colon B\longmapsto \frac{\psi_B(\mbfx)}{\pEsf{m}{\theta,\nu}[B]}=h_\mbfx(B) \comma \qquad B\in\Part{q}{m}\comma m\in\N_0\fstop
\end{equation}

We may now identify~$\KSimplex{q}$ as the Martin boundary~$\msM$.

\begin{prop}[Martin boundary of the PHC]\label{p:MartinBoundaryPartitions}
Assume~$\theta > 0$ and~$\nu\mathring{\Delta}^{q-1}>0$.
Let~$\seq{A_k}_k\subset\Part{q}{\cup}$ with~$A_k\to\infty$.
Then, the following are equivalent:
\begin{enumerate}[$(i)$]
\item\label{i:p:MartinBoundaryPartitions:1} $\seq{A_k}_k$ converges to the Martin boundary;
\item\label{i:p:MartinBoundaryPartitions:2} $\seq{\mbfX(A_k)}_k$ converges in~$(\KSimplex{q},\T_*)$.
\end{enumerate}
When this is the case, denoting by~$\mbfx$ the limit in~\ref{i:p:MartinBoundaryPartitions:2},
\begin{equation}\label{eq:p:MartinBoundaryPartitions:0}
\lim_{k\to\infty} \MartinK^\sym{q}(B,A_k) = \MartinK^\sym{q}(B,\mbfx) \comma \qquad B\in\Part{q}{\cup} \comma
\end{equation}
and two sequences converging to the Martin boundary are equivalent if and only if the corresponding limits~\ref{i:p:MartinBoundaryPartitions:2} coincide.
Furthermore, the map
\[
\Theta\colon \KSimplex{q}\longrar \msM \comma \qquad \Theta\colon \mbfx\longmapsto \MartinK^\sym{q}(\emparg,\mbfx)\comma
\]
is a homeomorphism; that is, the Martin boundary of the PHC is homeomorphic to~$(\KSimplex{q},\T_*)$.

\begin{proof}
Set~$n_k\eqdef \length{\col{A_k}}\to\infty$, and fix~$m\in\N_0$ and~$B\in\Part{q}{m}$.
By~\eqref{eq:d:MartinKernelPartitions:0} and Lemma~\ref{l:SamplingEstimate},
\begin{equation}\label{eq:p:MartinBoundaryPartitions:1}
\abs{\MartinK^\sym{q}(B,A_k) - \frac{\psi_B\tparen{\mbfX(A_k)}}{\pEsf{m}{\theta,\nu}[B]}} \leq \frac{1}{\pEsf{m}{\theta,\nu}[B]}\, \paren{1-\frac{\AntiPoch{n_k}{m}}{n_k^m}} \xrightarrow{k\to\infty} 0 \fstop
\end{equation}
Thus,~$\MartinK^\sym{q}(B,A_k)$ converges if and only if~$\psi_B\tparen{\mbfX(A_k)}$ does, and the two limits differ by the factor~$\pEsf{m}{\theta,\nu}[B]>0$ (cf.~\eqref{eq:PositivePESF}).

\paragraph{\ref{i:p:MartinBoundaryPartitions:2}$\implies$\ref{i:p:MartinBoundaryPartitions:1} and~\eqref{eq:p:MartinBoundaryPartitions:0}} 
If~$\mbfX(A_k)\to\mbfx$ in~$\T_*$, then~$\psi_B\tparen{\mbfX(A_k)}\to\psi_B(\mbfx)$, since~$\psi_B$ is $\T_*$-continuous by Corollary~\ref{c:MeasureDetermining}.
The conclusion follows by~\eqref{eq:p:MartinBoundaryPartitions:1} and~\eqref{eq:ExtendedMartinKernel}.

\paragraph{\ref{i:p:MartinBoundaryPartitions:1}$\implies$\ref{i:p:MartinBoundaryPartitions:2}} 
By~\eqref{eq:p:MartinBoundaryPartitions:1}, the sequence~$\psi_B\tparen{\mbfX(A_k)}$ converges for every~$B\in\Part{q}{\cup}$.
By Lemma~\ref{l:Compactness} the space~$(\KSimplex{q},\T_*)$ is compact metrizable. 
Let~$\mbfx,\mbfy$ be the limits of two convergent subsequences of~$\seq{\mbfX(A_k)}_k$.
By the $\T_*$-continuity of~$\psi_B$ we have~$\psi_B(\mbfx)=\lim_k \psi_B\tparen{\mbfX(A_k)}=\psi_B(\mbfy)$ for every~$B$, whence~$\mbfx=\mbfy$, since~$\set{\psi_B}_{B\in\Part{q}{\cup}}$ is point-separating by Corollary~\ref{c:MeasureDetermining}.
Since~~$(\KSimplex{q},\T_*)$ is metrizable, the whole sequence~$\seq{\mbfX(A_k)}_k$ converges to~$\mbfx$.

The assertion on equivalent sequences, as well as the injectivity of~$\Theta$, follow again from the fact that~$\set{\psi_B}_{B\in\Part{q}{\cup}}$ is point-separating, and~$\Theta$ is onto~$\msM$ by Lemma~\ref{l:DensityEmpirical} together with~\eqref{eq:p:MartinBoundaryPartitions:0}.
Since each~$\MartinK^\sym{q}(B,\emparg)$ is $\T_*$-continuous on~$\KSimplex{q}$ by Corollary~\ref{c:MeasureDetermining} and~\eqref{eq:ExtendedMartinKernel}, the map~$\Theta$ is continuous when~$\msM$ carries the topology of pointwise convergence.
As the latter is Hausdorff and~$(\KSimplex{q},\T_*)$ is compact,~$\Theta$ is a homeomorphism onto its image~$\msM$, which is therefore compact as well.
\end{proof}
\end{prop}

\begin{rem}[Martin compactification]\label{r:MartinCompactification}
No boundary point is a Martin kernel based at a state: if~$A\in \Part{q}{n}$ and~$\mbfx\in\KSimplex{q}$, choose~$m>n$ and, by~\eqref{eq:TotalPhiPsi:0b}, some~$B\in\Part{q}{m}$ with~$\psi_B(\mbfx)>0$.
Then~$\MartinK^\sym{q}(B,\mbfx)>0=\MartinK^\sym{q}(B,A)$ by~\eqref{eq:d:MartinKernelPartitions:0}.
Thus, by Proposition~\ref{p:MartinBoundaryPartitions}, the Martin compactification of~$\Part{q}{\cup}$ is the disjoint union~$\Part{q}{\cup}\sqcup\KSimplex{q}$, in which a sequence of states~$\seq{A_k}_k$ converges to~$\mbfx\in\KSimplex{q}$ precisely when~$\mbfX(A_k)\to\mbfx$ in~$\T_*$, i.e.\ precisely when the empirical block-color frequencies of~$A_k$ converge to~$\mbfx$.
\end{rem}

\subsubsection{Proof of Theorem~\ref{t:KingmanRepresentation}}\label{sss:MainProofs}
Let~$\probP_n$ be defined by~\eqref{eq:p:KingmanRepresentationConverse:0} for some~$\mu\in\msP(\KSimplex{q})$.
The measurability of~$\psi_A$ was shown in Corollary~\ref{c:PhiBorel}.
Since~$\mu$ is a probability measure, it follows from Lemma~\ref{l:TotalPhi} that~$\probP_n$ too is a probability measure for every~$n\in\N_0$, being a mixture of the probability measures~$A\mapsto \psi_A(\mbfx)$ for every~$\mbfx\in \KSimplex{q}$.
Exchanging the integration w.r.t.~$\mu$ in~\eqref{eq:p:KingmanRepresentationConverse:0} with the sum in~\eqref{eq:l:PsiConsistency:0} shows that~$\probP_n$ satisfies~\eqref{eq:d:Consistency:0}, which is enough to conclude the first assertion.

\medskip

In order to prove the converse, fix a $\CPS{q}$~$\probP_\bullet$.

\paragraph{Probability space}
Let~$\Omega$ be the Cartesian product~$\Omega = \prod_{n=0}^\infty \Part{q}{n}$, compact and separable in the product topology, and~$\msF$ be the product $\sigma$-algebra, which coincides with the Borel $\sigma$-algebra.
For each~$m\in\N_0$, let~$\omega_m$ be the $m^\text{th}$-coordinate function on~$\Omega$, and~$\msF_m$ be the $\sigma$-algebra on~$\Omega$ generated by~$\seq{\omega_n}_{n\geq m}$.
Finally, let~$\msF_\infty\eqdef \bigcap_m \msF_m$ be the tail $\sigma$-algebra on~$\Omega$.
Combining the consistency~\eqref{eq:d:Consistency:0} of~$\probP_\bullet$ with the Daniell--Kolmogorov Theorem proves the existence of a probability~$\probP$ on~$(\Omega,\msF)$ satisfying, for the kernels~\eqref{eq:Kernels},
\[
\probP[\omega_m = B] = \probP_m[B] \qquad \text{and} \qquad \probP[\omega_m=B \mid \omega_n,\omega_{n+1},\dotsc]= \BackwardK^\sym{q}_{m,n}(\omega_n,B)
\]
for every~$0\leq m<n<\infty$ and every~$B\in\Part{q}{m}$.

We thus have, for every~$m\in\N_0$, for every~$B\in\Part{q}{m}$,
\begin{equation}\label{eq:t:RepresentationProofsNuova:1}
\BackwardK^\sym{q}_{m,n}(\omega_n,B)=\probP[\omega_m=B \mid \omega_n,\omega_{n+1},\dotsc] = \ExpectP[\car_{\set{\omega_m=B}} \mid \msF_n] \comma
\end{equation}
which shows that~$\seq{\BackwardK^\sym{q}_{m,n}(\omega_n,B)}_{n > m}$ is a reverse $\msF_\bullet$-martingale under~$\probP$, converging to
\begin{equation}\label{eq:t:RepresentationProofsNuova:2}
\mathcal L_B\eqdef \nlim \BackwardK^\sym{q}_{m,n}(\omega_n,B) = \ExpectP[\car_{\set{\omega_m=B}} \mid \msF_\infty] 
\end{equation}
by Lévy's Downward Theorem~\cite[Thm.~14.4, p.~136]{Wil91}, since~$\car_{\set{\omega_m=B}}\in L^1(\probP)$.
Note that%
\begin{equation} \label{eq:t:RepresentationProofsNuova:3}
\ExpectP[\mathcal L_B]=\probP_m[B] \qquad \text{and} \qquad \sum_{B\in\Part{q}{m}} \mathcal L_B=1 \as{\probP} \comma\qquad m \in \N_0 \fstop
\end{equation}

\paragraph{Sampling with and without replacement}
For each~$n\in\N_1$, let~$\mbfX_\sym{n}\eqdef \mbfX(\omega_n)$ be as in~\eqref{eq:d:EmpiricalMap:0}.
By Lemma~\ref{l:SamplingEstimate},
\begin{equation}\label{eq:t:RepresentationProofsNuova:5}
\abs{\BackwardK^\sym{q}_{m,n}(\omega_n,B)- \psi_B(\mbfX_\sym{n})} \leq 1- \AntiPoch{n}{m}/n^m \comma \qquad B\in \Part{q}{\cup} %
\fstop
\end{equation}
Combining~\eqref{eq:t:RepresentationProofsNuova:1} and~\eqref{eq:t:RepresentationProofsNuova:5} shows that
\begin{equation}\label{eq:t:RepresentationProofsNuova:6}
\nlim \psi_B(\mbfX_\sym{n}) = \mathcal L_B \as{\probP}\comma \qquad B\in \Part{q}{\cup} \fstop
\end{equation}

Let us now focus on the case where~$B$ consists of a single block.
Recall the notation~\eqref{eq:KSimplexCoordinates}.
For fixed (deterministic)~$\mbfx\eqdef \ttseq{\mbfx_i}_{i\geq 0}\in \KSimplex{q}$, \eqref{eq:DefPhi}-\eqref{eq:DefPsi} give
\begin{equation}\label{eq:t:RepresentationProofsNuova:8}
\psi_{\car_\mbfr}(\mbfx) = \binom{\length{\mbfr}}{\mbfr} \Pmom{\mbfr}(\mbfx)\comma \quad \mbfr\in\N^q_*:\length{\mbfr}\geq 2 \comma \qquad \psi_{\car_{\mbfe_j}}(\mbfx) = \Dust_j(\mbfx)\comma \quad j\in [q] \fstop
\end{equation}
Choosing~$\mbfx=\mbfX_\sym{n}$ as above and applying~\eqref{eq:t:RepresentationProofsNuova:6} shows that the following limits exist
\begin{equation}\label{eq:t:RepresentationProofsNuova:9}
\Pmom{\mbfr}\eqdef \nlim \Pmom{\mbfr}(\mbfX_\sym{n}) \comma \qquad \Dust_j\eqdef \nlim \Dust_j(\mbfX_\sym{n})\comma
\end{equation}
and that~$\sum_{j=1}^q \Dust_j=1$.

\paragraph{Reconstruction via moments}
Recall the notation in~\eqref{eq:IsoSimplex} and define random probabilities
$
\Bulk_n \eqdef \IsoSimplex_2(\mbfX_\sym{n}) %
$. 
For every~$\mbfr\in \N^q_*$ we have
\begin{equation}\label{eq:t:RepresentationProofsNuova:11}
\int_{\CSimpl^q} \mbft^\mbfr \diff\Bulk_n(\mbft) = \sum_{j=1}^q \Pmom{\mbfr+\mbfe_j}(\mbfX_\sym{n}) \comma
\end{equation}
and we conclude by~\eqref{eq:t:RepresentationProofsNuova:9} that all moments of~$\Bulk_n$ converge~$\probP$-a.s.\ as~$n\to\infty$.
By compactness of~$\CSimpl^q$, hence of~$\msP(\CSimpl^q)$ in the weak topology, we conclude that~$\Bulk_n$ converges $\probP$-a.s., weakly as~$n\to\infty$, to a random probability~$\Bulk$ on~$\CSimpl^q$.
Together with the convergence of~$\Dust_j(\mbfX_\sym{n})$, this shows the existence of the limit~$(\DustV,\Bulk)\eqdef \nlim (\DustV(\mbfX_\sym{n}),\Bulk_n)$ $\probP$-a.s.\ in~$\msS$.
Since~$\seq{(\DustV(\mbfX_\sym{n}),\Bulk_n)}_n\subset \msS_*$ and~$\msS_*$ is closed by the proof of Lemma~\ref{l:Compactness},~$(\DustV,\Bulk)$~is $\probP$-a.s.\ an element of~$\msS_*$.
Since~$\msS_*=\IsoSimplex(\KSimplex{q})$ by Lemma~\ref{l:IsoSimplex}, there exists a random~$\mbfX\eqdef\seq{\mbfX_i}_{i\geq 0}\in\KSimplex{q}$ with~$\IsoSimplex(\mbfX)=(\DustV,\Bulk)$.
Finally, since~$\IsoSimplex\colon (\KSimplex{q},\T_*)\to \msS_*$ is a homeomorphism by definition of~$\T_*$, we have
\begin{equation}\label{eq:t:RepresentationProofsNuova:9.5}
\nlim \mbfX_\sym{n}=\mbfX \as{\probP} \quad \text{in} \quad (\KSimplex{q},\T_*)\fstop
\end{equation}

\paragraph{Measurability}
Since~$\mbfX_\sym{n}$ is~$\sigma(\omega_n)$-measurable and~$\mbfX=\nlim \mbfX_\sym{n}$ by~\eqref{eq:t:RepresentationProofsNuova:9.5}, the random variable~$\mbfX$ is~$\msF_m$-measurable for every~$m\in\N_1$, hence~$\msF_\infty$-measurable.

\paragraph{Representing measure}
By $\T_*$-continuity of~$\psi_A$ shown in Corollary~\ref{c:MeasureDetermining} and by~\eqref{eq:t:RepresentationProofsNuova:9.5},
\[
\nlim \psi_B(\mbfX_\sym{n}) = \psi_B(\mbfX) \as{\probP}\comma \qquad B\in\Part{q}{\cup}\comma
\]
and further combining this with~\eqref{eq:t:RepresentationProofsNuova:6} gives
\begin{equation}\label{eq:t:RepresentationProofsNuova:20}
\mathcal L_B=\psi_B(\mbfX) \as{\probP}\comma \qquad B\in\Part{q}{\cup}\fstop
\end{equation}

Let~$\mu$ be the law of~$\mbfX$ on~$\KSimplex{q}$.
By the first equality in~\eqref{eq:t:RepresentationProofsNuova:3}, taking expectations in~\eqref{eq:t:RepresentationProofsNuova:20} proves that
\[
\probP_m[B] = \int_{\KSimplex{q}} \psi_B(\mbfx) \diff\mu(\mbfx) \comma \qquad B\in \Part{q}{m}\comma\qquad m\in\N_0\comma 
\]
which is the representation formula~\eqref{eq:p:KingmanRepresentationConverse:0} for the probability~$\mu$ on~$\KSimplex{q}$ constructed above.

\paragraph{Uniqueness}
Since~$\set{\psi_B}_{B\in\Part{q}{\cup}}$ is a measure-determining class on~$\KSimplex{q}$ by Corollary~\ref{c:MeasureDetermining}, the mixing measure~$\mu$ is uniquely determined by the~$\CPS{q}$ $\probP_\bullet$, which completes the proof.

\subsubsection{Minimality and the Martin correspondence}\label{sss:Minimality}
We assume that~$\theta>0$ and~$\nu\in\msP(\Delta^{q-1})$ with~$\nu\mathring{\Delta}^{q-1}>0$.
We start with a lemma.
\begin{lem}%
\label{l:MinimalityPartitions}
$h_\mbfx=\MartinK^\sym{q}(\emparg,\mbfx)$ in~\eqref{eq:ExtendedMartinKernel} is minimal harmonic for all~$\mbfx\in\KSimplex{q}$.

\begin{proof}
We already know that~$h_\mbfx\in\mcH^\tgeq_{\theta,\nu}$ and~$h_\mbfx(0)=1$ by Lemma~\ref{l:HarmonicVsCPS}.
Let~$h\in\mcH^\tgeq_{\theta,\nu}$ with~$h\leq h_\mbfx$, and let~$\probP^h_\proc$ be as in~\eqref{eq:l:HarmonicVsCPS:0}, so that~$0\leq \probP^h_n[A]\leq \psi_A(\mbfx)$ for every~$A\in\Part{q}{n}$ and every~$n$.
As shown in the proof of Lemma~\ref{l:HarmonicVsCPS}, the quantity
\[
\tau\eqdef \sum_{A\in\Part{q}{n}} \probP^h_n[A] = h(0) \in [0,1]
\]
does not depend on~$n$.
Set~$R_n[A]\eqdef \psi_A(\mbfx)-\probP^h_n[A]\geq 0$, and note that~$R_\proc$, too, satisfies~\eqref{eq:d:Consistency:0}, with total mass~$1-\tau$ at each level, by linearity and Corollary~\ref{l:PsiConsistency}.

If~$\tau=0$, then~$\probP^h_\proc\equiv 0$ and thus~$h\equiv 0=0\cdot h_\mbfx$.
If~$\tau=1$, then~$R_\proc\equiv 0$ and thus~$h= h_\mbfx$.
Assume therefore~$\tau\in(0,1)$, so that~$\tau^{-1}\probP^h_\proc$ and~$(1-\tau)^{-1}R_\proc$ are both in~$\CPS{q}$ and
\[
\psi_\proc(\mbfx) = \tau\cdot \tau^{-1}\probP^h_\proc + (1-\tau)\cdot (1-\tau)^{-1}R_\proc \fstop
\]
Let~$\mu_1,\mu_2\in\msP(\KSimplex{q})$ be representing~$\tau^{-1}\probP^h_\proc$ and~$(1-\tau)^{-1}R_\proc$ in the sense of Theorem~\ref{t:KingmanRepresentation}.
By linearity of~\eqref{eq:p:KingmanRepresentationConverse:0}, the measure~$\tau\mu_1+(1-\tau)\mu_2$ represents~$\psi_\proc(\mbfx)$.
However, so does~$\delta_\mbfx$, again by~\eqref{eq:p:KingmanRepresentationConverse:0}.
By the uniqueness in Theorem~\ref{t:KingmanRepresentation} we conclude that~$\tau\mu_1+(1-\tau)\mu_2=\delta_\mbfx$.
Since Dirac masses are extremal in~$\msP(\KSimplex{q})$, we conclude that~$\mu_1=\mu_2=\delta_\mbfx$.
Thus~$\probP^h_\proc=\tau\, \psi_\proc(\mbfx)$, that is~$h=\tau h_\mbfx$, which concludes the proof.
\end{proof}
\end{lem}

\begin{cor}[Martin correspondence]\label{c:MartinCorrespondence}
Assume~$\theta > 0$ and~$\nu\mathring{\Delta}^{q-1}>0$.
The Martin boundary of the PHC coincides with its minimal Martin boundary, and the correspondences
\begin{equation*}
\begin{tikzcd}[row sep=-5pt]
\msP(\KSimplex{q}) \ar[r, leftrightarrow] &  \set{h\in \mcH^\tgeq_{\theta,\nu}: h(0)=1} \ar[r, leftrightarrow] & \CPS{q}
\\
\mu \ar[r, leftrightarrow]& \quad h_\mu\eqdef \displaystyle\int_{\KSimplex{q}} \MartinK^\sym{q}(\emparg,\mbfx)\diff\mu(\mbfx) \ar[r, leftrightarrow] & \quad \probP^{h}_\proc\quad
 \end{tikzcd}
 \end{equation*}
are one-to-one and onto, the second one being the bijection~\eqref{eq:l:HarmonicVsCPS:0}.
Their composition is the polychromatic Kingman representation~\eqref{eq:p:KingmanRepresentationConverse:0}.

\begin{proof}
	We know from Proposition~\ref{p:MartinBoundaryPartitions} that~$	\msM = \set{\MartinK^\sym{q}(\emparg,\mbfx) \, : \, \mbfx \in \KSimplex{q}}$. All such functions are minimal harmonic by Lemma~\ref{l:MinimalityPartitions}, whence the first assertion.

	By~\eqref{eq:ExtendedMartinKernel}, the map~$\mu \mapsto h_\mu$ may be written as
	\[
		h_\mu[B] = \int_{\KSimplex{q}} \MartinK^\sym{q}(B,\mbfx)\diff\mu(\mbfx) =  \frac{1}{\pEsf{m}{\theta,\nu}[B]} \int_{\KSimplex{q}} \psi_B(\mbfx)\diff\mu(\mbfx)
	\]
	for~$B \in \Part{q}{m}$ and~$m \in \N_0$. That is,~$\mu \mapsto h_\mu$ is the composition of the bijection~\eqref{eq:p:KingmanRepresentationConverse:0} with the inverse of~\eqref{eq:l:HarmonicVsCPS:0}. It follows, in particular, that all nonnegative harmonic functions~$h$ with~$h(0)=1$ are of the form~$h_\mu$.
\end{proof}
\end{cor}

The next corollary is a straightforward consequence of the fact that the functions~$\psi_A$ span a dense subset of~$\Cb(\KSimplex{q},\T_*)$ (Cor.~\ref{c:MeasureDetermining}). %

\begin{cor}[Continuity]\label{c:Continuity}
For each~$r\in\N_1$ let~$\probP^\sym{r}_\proc\eqdef \seq{\probP^\sym{r}_n}_n$ be a $\CPS{q}$ with representing measure~$\mu_r$ in the sense of Theorem~\ref{t:KingmanRepresentation}.
Then, the following are equivalent:
\begin{enumerate}[$(i)$]
\item for every~$n\in\N_1$, for every~$A\in\Part{q}{n}$, there exists
\[
\probQ_n[A]\eqdef \lim_{r\to\infty} \probP^\sym{r}_n[A] \fstop
\]
\item $\seq{\mu_r}_r$ converges weakly (relative to~$\T_*$) to a measure~$\mu\in \msP(\KSimplex{q})$.
\end{enumerate}

Furthermore, if either condition holds, then~$\seq{\probQ_n}_n$ is a $\CPS{q}$.
\end{cor}

Endow both~$\mcH^\tgeq_{\theta,\nu}$ and~$\CPS{q}$ with the topology of pointwise convergence on~$\Part{q}{n}$.
Recall a \emph{Bauer simplex} is a \emph{Choquet simplex} in which extremal points form a closed set.

\begin{cor}
The correspondences in Corollary~\ref{c:MartinCorrespondence} are affine homeomorphisms of Bauer simplices.
\begin{proof}
The affine property is clear. That~$\mu\mapsto h_\mu$ is continuous follows from the continuity of~$\psi_A$ (Cor.~\ref{c:MeasureDetermining}) together with the definition~\eqref{eq:ExtendedMartinKernel} of~$\MartinK^\sym{q}(\emparg,\mbfx)$.
Since~$(\KSimplex{q},\T_*)$ is compact by Lemma~\ref{l:Compactness}, so is~$\msP(\KSimplex{q})$ with the induced weak topology, by standard arguments.
Since~$\mu\mapsto h_\mu$ is an injective (Cor.~\ref{c:MartinCorrespondence}) map from a compact space into a Hausdorff space, it is a homeomorphism.
That~$h\mapsto \probP^h_\bullet$ is a homeomorphism is immediate by definition of the topologies of pointwise convergence.
Thus, all three spaces are Bauer simplices since so is~$\msP(\KSimplex{q})$.
\end{proof}
\end{cor}

\begin{rem}[Exit boundary]
Let~$\mbfX$ be the $\KSimplex{q}$-valued random variable constructed in the proof of Theorem~\ref{t:KingmanRepresentation}. By Proposition~\ref{p:MartinBoundaryPartitions} and Remark~\ref{r:MartinCompactification}, $\mbfX$~is exactly the~$\probP$-a.s.\ limit of~$\seq{\omega_n}_n$ in the Martin compactification of~$\Part{q}{\cup}$.
Since every boundary point is minimal by Lemma~\ref{l:MinimalityPartitions}, standard Doob--Martin theory, cf.~\cite[\S24]{Woe00}, further gives~$\msF_\infty=\sigma(\mbfX)$ up to~$\probP$-negligible sets. 
In other words, the Martin boundary of the PHC is also its exit boundary, and the representing measure~$\mu$ of~$\probP_\proc$ is the law of the almost sure limits of the empirical block-color frequencies~$a_{n,i,j}/n$ in the enumeration of~$\omega_n$ as in~\eqref{eq:d:EmpiricalMap:0}.
\end{rem}

\begin{rem}[Logical dependencies]\label{r:LogicalDependencies}
The proof of Lemma~\ref{l:MinimalityPartitions} uses the uniqueness part of Theorem~\ref{t:KingmanRepresentation}. 
Accordingly, Corollary~\ref{c:MartinCorrespondence} is to be read as an \emph{interpretation} of Theorem~\ref{t:KingmanRepresentation} in the language of Doob--Martin theory, \emph{not} as an independent proof of it.
Conversely, a direct proof of the minimality of~$\MartinK^\sym{q}(\emparg,\mbfx)$ would yield, together with Proposition~\ref{p:MartinBoundaryPartitions} and~\cite[Thms.~24.7 and~24.9]{Woe00}, an alternative proof of Theorem~\ref{t:KingmanRepresentation}.
We do not pursue this here.
\end{rem}

\subsection{Functoriality of the polychromatic Kingman representation}\label{ss:Functoriality}
We show that the polychromatic Kingman Representation Theorem~\ref{t:KingmanRepresentation} is functorial in an appropriate sense.

\subsubsection{Color degeneracy}\label{sss:Degeneracy}
Fix~$p,q\in\N_1$ with~$p\leq q$.
Any surjective map~$g\colon [q]\to [p]$ induces a push-forward map~$g_\pfwd\colon {[0,\infty)}^q\to {[0,\infty)}^p$ when regarding~${[0,\infty)}^\ell$ as the space of non-negative measures on~$[\ell]$ with~$\ell=p,q$.
Note that~$g_\pfwd$ maps $\Delta^{q-1}$ to $\Delta^{p-1}$.
Similarly, the restriction of~$g_\pfwd$ to~$\N^q_*$ takes values in~$\N^p_*$.

We call \emph{degeneracy map} the push-forward of multisets~$g^\deg\eqdef {g_\pfwd}_\mpfwd$ via~$g_\pfwd$, defined as in~\eqref{eq:MultiPush} with~$g_\pfwd$ in place of~$f$ there.
In other words, the degeneracy map~$g^\deg$ is the map on $q$-chromatic partitions resulting in the $p$-chromatic partition obtained by identifying colors~$j,j'\in [q]$ whenever~$g(j)=g(j')$.
In turn, $g^\deg$~induces a push-forward map~${g^\deg}_\pfwd$ on measures over polychromatic partitions, in the standard way.

Let us list some elementary facts:
\begin{enumerate*}[$(i)$]
\item $g_\pfwd$ is linear and mass preserving, i.e.~$\length{g_\pfwd\mbfa}=\length{\mbfa}$;
\item $g^\deg$ recolors blocks in~$A$ without merging them;
\item thus $\size{g^\deg A}=\size{A}$ and $g_\pfwd\supp A= \supp(g^\deg A)$;
\item with the notation of~\eqref{eq:AmenoEj},~$g^\deg(A_{\setminus \mbfa,j})= g^\deg(A)_{\setminus g_\pfwd \mbfa, g(j)}$;
\item $h^\deg\circ g^\deg= (h\circ g)^\deg$ for every surjective~$g\colon [q]\to [p]$ and surjective~$h\colon [p]\to [r]$;
\item and $\shape=\car^\deg$ is the degeneracy map of the constant map~$\car\colon [q] \to [1]$.
\item as a consequence, $\shape(g^\deg A)=\shape(A)$.
\end{enumerate*}

\begin{prop}\label{p:ConsistencyDegeneracy}
For~$m,n\in \N_0$ and~$p,q\in\N_1$, with~$p\leq q$, and %
surjective~$g\colon [q]\to [p]$,
\[
{g^\deg}_\pfwd\circ  \sigma^\sym{q}_{m,n}= \sigma^\sym{p}_{m,n}\circ {g^\deg}_\pfwd \fstop
\]

\begin{proof}
By the cocycle property for~$\sigma^\sym{q}$ and~$\sigma^\sym{p}$, it suffices to prove the assertion with~$m=n-1$.
For any~$\probP_n\in\msP(\Part{q}{n})$ we have
\begin{align*}
\tparen{{g^\deg}_\pfwd (\sigma^\sym{q}_{n-1, n} \probP_n)} [D] &= \sum_{\substack{B \in [g^\deg]^{-1}(D)}} \sum_{A\in\Part{q}{n}} \BackwardK^\sym{q}_{n-1, n}(A,B) \probP_n[A] 
\\
&= \sum_{A\in\Part{q}{n}} \probP_n[A] \sum_{\substack{B \in [g^\deg]^{-1}(D)}} \BackwardK^\sym{q}_{n-1, n}(A,B)\comma
\\
\tparen{\sigma^\sym{p}_{n-1, n} ({g^\deg}_\pfwd \probP_n)}[D] &= \sum_{C\in \Part{p}{n}} \BackwardK^\sym{p}_{n-1, n}(C,D)\sum_{\substack{A \in [g^\deg]^{-1}(C)}} \probP_n[A] 
\\
&= \sum_{A\in\Part{q}{n}} \probP_n[A]\, \BackwardK^\sym{p}_{n-1, n}\tparen{g^\deg(A),D} \fstop
\end{align*}
Thus, it suffices to show that
\begin{equation}\label{eq:p:ConsistencyDegeneracy:1}
\sum_{\substack{B \in [g^\deg]^{-1}(D)}} \BackwardK^\sym{q}_{n-1, n}(A,B) = \BackwardK^\sym{p}_{n-1, n}\tparen{g^\deg(A),D} \comma \qquad A\in\Part{q}{n}\comma D\in\Part{p}{n-1} \fstop
\end{equation}

\paragraph{Conclusion}
The left-hand side of~\eqref{eq:p:ConsistencyDegeneracy:1} is
\begin{align*}
\sum_{\substack{B \in [g^\deg]^{-1}(D)}} \BackwardK^\sym{q}_{n-1, n}(A,B) &= \sum_{\mbfa\in\supp A}\sum_{j=1}^q \BackwardK^\sym{q}_{n-1, n}\tparen{A,A_{\setminus\mbfa,j}} \car_{g^\deg(A_{\setminus\mbfa,j})=D} \fstop
\end{align*}
Let~$C \coloneqq g^\deg(A)$. By respectively reindexing the sums over~$\mbfa\in \supp A$ and~$j\in [q]$ by $g_\pfwd\mbfa=\mbfc\in\supp C$ and~$g(j)=k\in [p]$,
\begin{align*}
\sum_{\substack{B \in [g^\deg]^{-1}(D)}} \BackwardK^\sym{q}_{n-1, n}(A,B) &= \sum_{\mbfc\in\supp C} \, \sum_{\substack{\mbfa\in \N^q_* \\ g_\pfwd \mbfa=\mbfc}} \, \sum_{k=1}^p \, \sum_{j\in g^{-1}(k)} \frac{a_j A(\mbfa)}{n} \car_{C_{\setminus\mbfc,k}=D} \\ %
&=\sum_{\mbfc\in\supp C} \sum_{k=1}^p \frac{\car_{C_{\setminus\mbfc,k}=D}}{n} \sum_{\substack{\mbfa\in \N^q_* \\ g_\pfwd \mbfa=\mbfc}} A(\mbfa) \sum_{j\in g^{-1}(k)} a_j
\\
&=\sum_{\mbfc\in\supp C} \sum_{k=1}^p \frac{\car_{C_{\setminus\mbfc,k}=D}}{n} \, C(\mbfc)\, c_k
= \BackwardK^\sym{p}_{n-1, n}(C,D)\fstop
\end{align*}
Recalling that~$C\eqdef g^\deg(A)$ shows~\eqref{eq:p:ConsistencyDegeneracy:1}, which proves the assertion.
\end{proof}
\end{prop}

\begin{cor}\label{c:ShapePartitionPartition}
Fix~$p,q\in\N_1$, with~$p\leq q$, and a surjective~$g\colon [q]\to [p]$.
If~$\probP_\bullet$ is a~$\CPS{q}$, then~${g^\deg}_\sharp\probP_\bullet$ is a~$\CPS{p}$.
In particular,~$\seq{\shape_\pfwd \probP_n}_n$ is a partition structure.
\end{cor}

\subsubsection{Functoriality}
We now move to maps on the polychromatic Kingman simplex.

\begin{defs}\label{d:LiftedDegeneracy}
For~$g\colon [q]\to [p]$, the \emph{lifted degeneracy map}~$G_g\colon \KSimplex{q}\to\KSimplex{p}$ is
\[
G_g(\mbfx)_0\eqdef g_\pfwd\mbfx_0 \comma \qquad \seq{G_g(\mbfx)_i}_{i\geq 1} \eqdef \sort_\shortlex\tparen{\seq{g_\pfwd\mbfx_i}_{i\geq 1}}\comma
\]
where~$\sort_\shortlex$ denotes the reordering of the argument decreasing in the shortlex order.
\end{defs}
We note that~$G_g$ is well-defined, since~$\length{g_\pfwd\mbfx_i}=\length{\mbfx_i}$ for every~$i\in\N_1$ and thus $\sum_{i\geq 0} \length{G_g(\mbfx)_i}=1$, and further that it is functorial, in the sense that it satisfies
\[
G_{\id_{[q]}}=\id_{\KSimplex{q}} \qquad \text{and} \qquad G_g\circ G_h = G_{g\circ h}\fstop
\]
The reordering only permutes vectors of equal size, and there are finitely many such vectors for each positive size, so the map is Borel measurable.

\begin{lem}[Intertwining]\label{l:Intertwining}
For every~$n\geq 0$, every~$B\in\Part{p}{n}$, and every~$\mbfx\in\KSimplex{q}$,
\begin{equation}\label{eq:l:Intertwining:0}
\sum_{A\in [g^\deg]^{-1}(B)} \psi_A(\mbfx) = \psi_B\tparen{G_g(\mbfx)}\fstop
\end{equation}

\begin{proof}
By~\eqref{eq:alterPsi}, we have
\begin{equation}\label{eq:l:Intertwining:1}
\psi_A(\mbfx) = n! \sum_{(\mbfm,\boldnu)\in {\mathcal I}_A} \frac{\mbfx_0^\mbfm}{\mbfm!} \prod_{i = 1}^\infty \frac{\mbfx_i^{\boldnu_i}}{\boldnu_i!}\comma
\end{equation}
where the summation runs over the set~${\mathcal I}_A$ of all pairs~$(\mbfm,\boldnu)$ with~$\mbfm\in \N_0^q$ and~$\boldnu\eqdef\seq{\boldnu_i}_{i\geq 1}$,  $\boldnu_i\eqdef\seq{\nu_{i,j}}_{j\in [q]}$ such that~$\uni{\mbfm}+A_{\boldnu}=A$. 
By the Multinomial Theorem,
\begin{equation}\label{eq:l:Degeneracy:1}
\frac{(g_\pfwd\mbfu)^\mbfb}{\mbfb!} = \frac{1}{\mbfb!} \prod_{k=1}^p \bracket{\sum_{j\in g^{-1}(k)}u_j}^{b_k}= \sum_{\mbfa\in \N_0^q: g_\pfwd\mbfa=\mbfb} \frac{\mbfu^\mbfa}{\mbfa!}\comma \qquad \mbfu\in [0,\infty)^q\fstop
\end{equation}
Applying now~\eqref{eq:l:Degeneracy:1} to each~$\mbfx_i$ and expanding the product,
\[
\frac{(g_\pfwd\mbfx_0)^{\mbfm'}}{\mbfm'!} \prod_{i\geq 1} \frac{(g_\pfwd\mbfx_i)^{\boldnu_i'}}{\boldnu_i'!} = \sum_{\substack{(\mbfm,\boldnu): g_\pfwd\mbfm=\mbfm'\\ g_\pfwd\had\boldnu= \boldnu'}} \frac{\mbfx_0^\mbfm}{\mbfm!} \prod_{i\geq 1} \frac{\mbfx_i^{\boldnu_i}}{\boldnu_i!} \comma
\]
where we set~$g_\pfwd\had\boldnu\eqdef \seq{g_\pfwd\boldnu_i}_{i\geq 1}$.
Using that~$g^\deg$ preserves supports, that it is additive on multisets, and that~$g^\deg\uni{\mbfm}=\uni{g_\pfwd\mbfm}$ and~$g^\deg A_\boldnu= A_{g_\pfwd\had\boldnu}$, we conclude that~$(\mbfm,\boldnu)\in {\mathcal I}_A$ for some~$A\in (g^\deg)^{-1}(B)$ if and only if~$(g_\pfwd\mbfm,g_\pfwd\had\boldnu)\in {\mathcal I}_B$.
Writing~\eqref{eq:l:Intertwining:1} for~$B$ computed at~$G_g(\mbfx)$ and reindexing gives~\eqref{eq:l:Intertwining:0}, since the sets~$({\mathcal I}_A)_{A\in \Part{q}{n}}$ are pairwise disjoint and since the formula~\eqref{eq:l:Intertwining:1} for~$\psi_A$ is symmetric in the variables~$(\mbfx_i)_{i \ge 1}$, so that the reordering in Definition~\ref{d:LiftedDegeneracy} is immaterial.
\end{proof}
\end{lem}

\begin{prop}\label{p:Functoriality}
Let~$g\colon [q]\to [p]$ be surjective.
If~$\probP_\bullet\in\CPS{q}$ is represented by~$\mu\in\msP(\KSimplex{q})$ in the sense of Theorem~\ref{t:KingmanRepresentation}, then~${g^\deg}_\pfwd\probP_\bullet\in\CPS{p}$ as in Corollary~\ref{c:ShapePartitionPartition} is represented by~$(G_g)_\pfwd\mu$.
In particular, ${g^\deg}_\pfwd$~maps extremal~$\CPS{q}$'s to extremal~$\CPS{p}$'s, viz.
\[
{g^\deg}_\pfwd\psi_\proc(\mbfx)=\psi_\proc\tparen{G_g(\mbfx)}\fstop
\]
\begin{proof}
The assertion for extremal structures is an immediate consequence of Lemma~\ref{l:Intertwining}.
For general structures it follows by linearity of~$(G_g)_\pfwd$ and disintegration.
\end{proof}
\end{prop}

\subsection{Special cases}
In this section we collect some specializations of Theorem~\ref{t:KingmanRepresentation} to particular classes of $\CPS{q}$'s, including those \emph{of partially exchangeable type} which we now define, as well as multipartition structures.

\subsubsection{Partially exchangeable partition structures}\label{sss:PEPS}
For~$\mbfn\in\N_0^q$ and~$j\in [q]$ with~$n_j\geq 1$, let~$\BackwardK^{\sym{q},j}_\mbfn$ be the
kernel from~$\mcA_\mbfn$ to~$\mcA_{\mbfn-\mbfe_j}$ describing the deletion of one element of
color~$j$ chosen uniformly at random among the~$n_j$ available, viz.
\begin{equation}\label{eq:WithinGroupKernel}
\BackwardK^{\sym{q},j}_\mbfn(A,B)\eqdef \begin{cases} \frac{a_j A(\mbfa)}{n_j} & \text{if } B=A_{\backslash \mbfa,j}\\ 0 & \text{otherwise}\end{cases}\comma \qquad A\in\mcA_\mbfn\comma B\in \mcA_{\mbfn-\mbfe_j} \comma
\end{equation}
and let~$\sigma^{\sym{q},j}_\mbfn$ be the induced map on measures.
Comparison with~\eqref{eq:Kernels:B} gives
\begin{equation}\label{eq:KernelFactorization}
\BackwardK^\sym{q}_{n-1,n}(A,B)= \frac{n_j}{n}\, \BackwardK^{\sym{q},j}_\mbfn(A,B) \comma \qquad A\in\mcA_\mbfn\comma B\in \mcA_{\mbfn-\mbfe_j}\fstop
\end{equation}

\begin{defs}\label{d:PEPS}
A family~$\probQ_\proc\eqdef\seq{\probQ_\mbfn}_{\mbfn\in\N_0^q}$, with~$\probQ_\mbfn\in\msP(\mcA_\mbfn)$, is called a
\emph{partially exchangeable partition structure} if
$\probQ_{\mbfn-\mbfe_j}=\sigma^{\sym{q},j}_\mbfn\probQ_\mbfn$ for every~$\mbfn$ and~$j\in[q]$
with~$n_j\geq 1$.
We say that~$\probP_\proc\in\CPS{q}$ is \emph{of partially exchangeable type} if there exist
$\mbfp\in\Delta^{q-1}$ and a partially exchangeable partition structure~$\probQ_\proc$ with
\begin{equation}\label{eq:PEType}
\probP_n=\sum_{\mbfn\in\N_0^q\,:\,\length{\mbfn}=n} \binom{n}{\mbfn}\, \mbfp^\mbfn\, \probQ_\mbfn \comma \qquad n\in\N_0\fstop
\end{equation}
\end{defs}

We present the polychromatic analogue of Kingman's partition in~\cite[p.~375]{Kin78}.
Fix~$\mbfx\in\KSimplex{q}$ and let~$\seq{\xi_r}_{r\in\N_1} = \seq{(\beta_r,\gamma_r)}_{r\in\N_1}$ be independent with~\eqref{eq:xir}.
We denote by~$K_\mbfx\in\Part{q}{n}$ the random $q$-chromatic partition of~$[n]$ induced by grouping elements of~$[n]$ according to~$\seq{\beta_r}_{r\in [n]}$ and coloring them with~$\seq{\gamma_r}_{r \in [n]}$. The following is an immediate consequence of Remark~\ref{rem:probKing}.

\begin{lem}[Paintbox sampling]\label{l:Paintbox}
$K_\mbfx$ has law~$\psi_n(\mbfx)$.
\end{lem}

\begin{lem}[Color refinement of~\eqref{eq:TotalPhiPsi:0b}]\label{l:PsiColor}
For every~$\mbfx\in\KSimplex{q}$ and every~$\mbfn\in\N_0^q$,
\begin{equation}\label{eq:PsiColor}
\sum_{A\vDash \mbfn} \psi_A(\mbfx) = \binom{\length{\mbfn}}{\mbfn} \DustV(\mbfx)^{\mbfn}\fstop
\end{equation}
\begin{proof}
The random variables~$\seq{\gamma_r}_r$ are i.i.d.\ $\Cat_{\DustV(\mbfx)}$,
since a single draw has color~$j$ with probability~$\sum_{i\geq 0}x_{i,j}=\Dust_j(\mbfx)$. By Lemma~\ref{l:Paintbox}, the
left-hand side of~\eqref{eq:PsiColor} is the probability that their counts equal~$\mbfn$.
\end{proof}
\end{lem}

\begin{prop}\label{p:DiracCriterion}
Let~$\probP_\proc\in\CPS{q}$ with representing measure~$\mu\in\msP(\KSimplex{q})$ as in
Theorem~\ref{t:KingmanRepresentation}.
Then, the following are equivalent:
\begin{enumerate}[$(i)$]
\item\label{i:p:DiracCriterion:1} $\probP_\proc$ is of partially exchangeable type;
\item\label{i:p:DiracCriterion:2} $\DustV_\pfwd\mu=\delta_\mbfp$ for some~$\mbfp\in\Delta^{q-1}$.
\end{enumerate}
In this case,~$\mbfp$ is uniquely determined by
\begin{equation}
	p_j=\probP_1\tbraket{\car_{\mbfe_j}}\comma \qquad j \in [q] \fstop
\end{equation}
If, additionally,~$\mbfp\in\mathring{\Delta}^{q-1}$, then~$\probQ_\proc$ in~\eqref{eq:PEType} satisfies
\begin{equation}\label{eq:PEDisintegration}
\probQ_\mbfn[A]=\binom{\length{\mbfn}}{\mbfn}^{-1}\mbfp^{-\mbfn}\, \probP_{\length{\mbfn}}[A] \comma \qquad A\vDash \mbfn \comma \mbfn \in \N_0^q \fstop
\end{equation}

\begin{proof}
\paragraph{\ref{i:p:DiracCriterion:1}$\implies$\ref{i:p:DiracCriterion:2}}
Summing~\eqref{eq:PEType} over~$A\vDash\mbfn$ shows that~$\probP_{\length{\mbfn}}\tbraket{\mcA_\mbfn}=\binom{\length{\mbfn}}{\mbfn}\mbfp^\mbfn$,
while~\eqref{eq:p:KingmanRepresentationConverse:0} and Lemma~\ref{l:PsiColor} give~$\probP_{\length{\mbfn}}\tbraket{\mcA_\mbfn}=\binom{\length{\mbfn}}{\mbfn}\int \DustV^\mbfn\diff\mu$.
Hence, all moments of~$\DustV_\pfwd\mu$ agree with those of~$\delta_\mbfp$.
Since~$\Delta^{q-1}$ is compact, the multi-dimensional Hausdorff moment problem gives~$\DustV_\pfwd\mu=\delta_\mbfp$.
Taking~$\mbfn=\mbfe_j$ identifies~$\mbfp$, and~\eqref{eq:PEType} then forces~\eqref{eq:PEDisintegration}.

\paragraph{Reduction} Assume~\ref{i:p:DiracCriterion:2} and, with no loss of generality,~that colors are ordered so that~$\mbfp = (\mbfp',\zero)$ with~$\mbfp' \in \mathring{\Delta}^{q'-1}$,~$q' \le q$. If~$j \in [q] \setminus [q']$, then~$(D_j)_\pfwd \mu = \delta_0$, and thus~$x_{i,j}=0$ for all~$i$'s, for~$\mu$-a.e.~$\mbfx$; by~\eqref{eq:DefPhi}-\eqref{eq:DefPsi} and~\eqref{eq:p:KingmanRepresentationConverse:0} we get~$\probP_n[A]=0$ whenever color~$j$ occurs in~$A$, i.e.,~$\probP_\proc$ is the trivial extension of a~$\CPS{q'}$~$\probP_\proc'$ with~$q'$ colors.
This~$\probP_\proc'$ still satisfies~\ref{i:p:DiracCriterion:2} (with~$q'$ in place of~$q$).
Furthermore, if~$\probP_\proc'$ is of partially exchangeable type for ($\mbfp'$ and) some~$\probQ_\proc'$, one can show that~$\probP_\proc$ is of partially exchangeable type with ($\mbfp$ and)
\begin{multline*}
	\probQ_{(\mbfn',\mbfn'')}\big[A' + \uni{(\zero,\mbfn'')}\big] \coloneqq \probQ'_{\mbfn'}[A'] \comma \qquad A' \in \Part{q'}{\mbfn'} \comma \quad \mbfn' \in \N_0^{q'} \comma \quad \mbfn'' \in \N_0^{q-q'} \comma
\end{multline*}
and~$\probQ_{\mbfn}[A] \coloneqq 0$ for all other~$A$'s. In the proof of \ref{i:p:DiracCriterion:2}$\implies$\ref{i:p:DiracCriterion:1}, we may thus assume~$\mbfp\in\mathring{\Delta}^{q-1}$. 

\paragraph{\ref{i:p:DiracCriterion:2}$\implies$\ref{i:p:DiracCriterion:1}}
For~$\mbfx \in \KSimplex{q}$ and set~$\psi^\mbfn_A(\mbfx)\eqdef\binom{\length{\mbfn}}{\mbfn}^{-1}\mbfp^{-\mbfn}\psi_A(\mbfx)$
for~$A\vDash\mbfn$, a probability on~$\mcA_\mbfn$ by Lemma~\ref{l:PsiColor}. 

\paragraph{Claim:~$\psi^{\col{\proc}}_\proc(\mbfx)$ is a partially exchangeable partition structure}
For every~$j \in [q]$ and~$B \in \Part{q}{\mbfn-\mbfe_j}$,~\eqref{eq:l:poly1:00} with~$\mbfa = \mbfe_j$ gives
{\begin{align*}
	&n \, P_{\mbfe_j}(\mbfx) \, \varphi_B(\mbfx) \\
	&\quad= (B(\mbfe_j)+1) \varphi_{B+\car_{\mbfe_j}}(\mbfx) + \! \sum_{\mbfb \in \supp B} \! (b_j+1) (B(\mbfb+\mbfe_j)+1) \varphi_{B-\car_{\mbfb}+\car_{\mbfb+\mbfe_j}}(\mbfx) \\
	&\quad = \sum_{A \in \Part{q}{\mbfn}} n_j \, \BackwardK^{\sym{q},j}_\mbfn(A,B) \varphi_A(\mbfx) \fstop
\end{align*}}
Using Lemma~\ref{l:limphi}, we find
\begin{align*}
	n \, D_j(\mbfx) \, \psi_B(\mbfx) = \sum_{A \in \Part{q}{\mbfn}} n_j \, \BackwardK^{\sym{q},j}_\mbfn(A,B) \psi_A(\mbfx) \comma \qquad \mbfx \in \KSimplex{q} \fstop
\end{align*}
For~$\mu$-almost every~$\mbfx$, we may replace~$D_j(\mbfx)$ with~$p_j$ in the latter formula.
Replacing~$\psi_B(\mbfx)=\binom{\length{\mbfn}-1}{\mbfn-\mbfe_j}\mbfp^{\mbfn-\mbfe_j} \psi_B^{\mbfn-\mbfe_j}$ and~$\psi_A(\mbfx)=\binom{\length{\mbfn}}{\mbfn}\mbfp^{\mbfn} \psi_A^{\mbfn}$ proves the claim.

\paragraph{Conclusion} Define
\[
\probQ_\mbfn[A]=\int_{\KSimplex{q}} \psi^\mbfn_A(\mbfx)\diff\mu(\mbfx)\comma \qquad A\vDash\mbfn\fstop
\]
Since Definition~\ref{d:PEPS} is an affine condition, it passes to such mixtures, and~$\probQ_\proc$ is a partially exchangeable partition structure.
Then~\eqref{eq:PEType} holds by the definitions of~$\psi_A^\mbfn$ and~$Q_\mbfn$.
\end{proof}
\end{prop}

\begin{rem}[Comparison with~\cite{FraLijPruReb25}]\label{r:Comparison}
It is shown in~\cite[Thm.~5]{FraLijPruReb25} that the partially exchangeable partition law induced by a vector of possibly dependent random probability measures as in~\cite[p.~17]{FraLijPruReb25} depends on the vector only through the law of its atomic weights. In particular, it is reproduced by a multivariate species sampling process, whose atoms are i.i.d.\ from an arbitrary non-atomic base measure, independently of the weights.
To match their notation, let~$q=J$, fix~$\mbfp\in\mathring{\Delta}^{q-1}$, and assign to each individual an independent color (group) with
law~$\Cat_\mbfp$.
Deleting an individual chosen uniformly from a sample with color counts~$\mbfn$ amounts to first choosing a color~$j \in [q]$ with probability~$n_j / \length{\mbfn}$, and then deleting an individual uniformly within that color group. This turns a
partially exchangeable array with~$q$ groups into a $q$-chromatic partition structure of the form~\eqref{eq:PEType}.
Theorem~\ref{t:KingmanRepresentation} then provides a unique representing measure~$\mu\in\msP(\KSimplex{q})$ and Proposition~\ref{p:DiracCriterion} gives~${\DustV}_\pfwd\mu=\delta_\mbfp$.
Thus~$\pi_{j,i}\eqdef x_{i,j}/p_j$, with~$i\in\N_1$ and~$j\in[q]$, defines for $\mu$-a.e.~$\mbfx$ a random
sub-probability array; the dust~$\mbfx_0$ accounts for the improper part
$1-\sum_{i\geq 1}\pi_{j,i}=x_{0,j}/p_j$.
This is precisely the \emph{multivariate species sampling process}
of~\cite[Dfn.~1]{FraLijPruReb25}.
The underlying space plays no role, since the induced partition records only which sampled values coincide, and any two non-atomic Borel probability
measures on standard Borel spaces are isomorphic.
\end{rem}

\subsubsection{Kingman representation of the pESF}\label{ss:RepresentationPESF}
Fix~$\theta>0$ and~$\nu\in\msP(\Delta^{q-1})$. Further set~$\OSimplex{q} \eqdef \set{\mbfx\in \KSimplex{q}: \mbfx_0=\zero}$.
Let~$\boldell\eqdef\seq{\ell_i}_{i\geq 0}\in \OSimplex{1}$ be distributed according to the Poisson--Dirichlet distribution~$\PD{\theta}$ on the standard ordered infinite simplex~$\OSimplex{1}$, let~$\seq{\mbfp_i}_{i\geq 1}\subset \Delta^{q-1}$ be i.i.d.\ $\nu$-distributed, and set
\[
\mbfx_0\eqdef \zero\comma \quad \mbfx_i\eqdef \ell_i \mbfp_i\comma \quad i\geq 1\fstop
\]
Since~$\length{\mbfx_i}=\ell_i>\ell_{i+1}$ $\PD{\theta}$-a.s. for each~$i\in\N_1$ and~$\sum_{i=1}^\infty \ell_i=1$, we have~$\mbfx\eqdef\seq{\mbfx_i}_{i\geq 0}\in \KSimplex{q}$.

\begin{defs}\label{d:polyPD}
We call \emph{polychromatic Poisson--Dirichlet measure with parameters~$\theta>0$ and $\nu\in\msP(\Delta^{q-1})$ the measure on~$\KSimplex{q}$}
\[
\PD{\theta,\nu}\eqdef \law\tparen{(\mbfx_0,\mbfx_1,\mbfx_2,\dotsc)} \fstop
\]
\end{defs}

\begin{prop}\label{p:RepresentationPESF}
$\pEsf{\bullet}{\theta,\nu}$ is a $\CPS{q}$ represented by~$\mu\eqdef \PD{\theta,\nu}$ in the sense of Theorem~\ref{t:KingmanRepresentation}.
\end{prop}

In other words, $\pEsf{\bullet}{\theta,\nu}$ is the `chromatic paintbox' in which the block frequencies $\seq{\length{\mbfx_i}}_i$ are $\PD{\theta}$-distributed and, conditionally on them, the normalized color vectors $\mbfx_i/\length{\mbfx_i}$ are i.i.d.~$\nu$.

\begin{proof}[Proof of Proposition~\ref{p:RepresentationPESF}]
Since~$\mbfx_0=\zero$ $\PD{\theta,\nu}$-a.s., $\psi_A=\varphi_A$ $\PD{\theta,\nu}$-a.s.
Fix~$A\in\Part{q}{n}$ enumerated by~$\seq{\mbfa_h}_h$ and set~$\boldlambda\eqdef \shape(A)$.
Since the indices~$\iota(h)$ in~\eqref{eq:DefPhi} are pairwise distinct (by injectivity of~$\iota$) and the marks~$\mbfp_i$ are i.i.d.,
\[
\CondExpectP{\prod_{h=1}^{\size{A}} \mbfx_{\iota(h)}^{\mbfa_h}}{\boldell} = \prod_{h=1}^{\size{A}} \ell_{\iota(h)}^{\length{\mbfa_h}} \mom{\mbfa_h}{\nu} \comma \qquad \iota\in T_{\size{A}}\comma
\]
so that, by definition~\eqref{eq:DefPhi}  of~$\varphi_A$, and setting~$\Sigma_{\boldlambda}(\boldell)\eqdef \sum_{\iota\in T_{\size{A}} } \prod_{h=1}^{\size{A}} \ell_{\iota(h)}^{\length{\mbfa_h}}$, we have
\[
\int_{\KSimplex{q}} \varphi_A \diff\PD{\theta,\nu} = n! \prod_{\mbfa\in\supp A}\frac{\mom{\mbfa}{\nu}^{A(\mbfa)}}{\mbfa!^{A(\mbfa)}A(\mbfa)!}  \, \ExpectP_{\PD{\theta}}[\Sigma_{\boldlambda}(\boldell)]\fstop
\]
By Kingman's monochromatic representation~\cite{Kin77},~$\ExpectP_{\PD{\theta}}[\Sigma_{\boldlambda}(\boldell)] = \Esf{n}{\theta}[\boldlambda] (n!)^{-1} \prod_{i=1}^n (i!)^{\lambda_i}\lambda_i!$.
Since~$\prod_{i=1}^n (i!)^{\lambda_i}= \prod_{\mbfa\in\supp A} \length{\mbfa}!^{A(\mbfa)}$,
\[
\int_{\KSimplex{q}} \varphi_A \diff\PD{\theta,\nu} = \Esf{n}{\theta}[\boldlambda] \prod_{i=1}^n \lambda_i! \prod_{\mbfa\in\supp A} \frac{1}{A(\mbfa)!} \paren{\binom{\length{\mbfa}}{\mbfa}\mom{\mbfa}{\nu}}^{A(\mbfa)} = \pEsf{n}{\theta, \nu}[A] \fstop
\]
It then follows from Theorem~\ref{t:KingmanRepresentation} that~$\pEsf{\bullet}{\theta,\nu}$ is a $\CPS{q}$ represented by~$\mu\eqdef \PD{\theta,\nu}$.
\end{proof}

Let~$\Bary{\nu}\in\Delta^{q-1}$ be the barycenter of~$\nu\in\msP(\Delta^{q-1})$ and~$\iota_1\colon \Delta^{q-1}\to\KSimplex{q}$ be the map
\[
\iota_1\colon \mbft \longmapsto \seq{\zero,\mbft,\zero,\zero,\dotsc}\fstop
\]

\begin{cor}[Asymptotics]\label{c:Asymptotics}
We have
\[
\lim_{\theta\downarrow 0} \PD{\theta,\nu} = (\iota_1)_\pfwd \nu \qquad \text{and} \qquad  \lim_{\theta\uparrow \infty} \PD{\theta,\nu} = \delta_{\seq{\Bary{\nu}, \zero,\zero,\dotsc}}\fstop
\]
\begin{proof}
Combine Proposition~\ref{p:RepresentationPESF} with Corollaries~\ref{c:AsymptoticsPESF} and~\ref{c:Continuity}.
\end{proof}
\end{cor}

\begin{cor}\label{c:AggregationPESF}
Let~$g\colon [q]\to [p]$ be surjective.
Then,~${g^\deg}_\pfwd \pEsf{\bullet}{\theta,\nu} = \pEsf{\bullet}{\theta, {g_\pfwd}_\pfwd\nu}$.
\begin{proof}
By Propositions~\ref{p:RepresentationPESF} and~\ref{p:Functoriality},~${g^\deg}_\pfwd \pEsf{\bullet}{\theta,\nu}$ is the unique~$\CPS{p}$ represented by~$(G_g)_\pfwd \PD{\theta,\nu}$, with~$\PD{\theta,\nu}$ as in Definition~\ref{d:polyPD}.
It is then straightforward that~$(G_g)_\pfwd\PD{\theta,\nu}=\PD{\theta,{g_\pfwd}_\pfwd\nu}$ and the conclusion follows by a second application of Proposition~\ref{p:RepresentationPESF}.
\end{proof}
\end{cor}

We now determine when~$\pEsf{\bullet}{\theta,\nu}$ is of partially exchangeable type, cf.~Definition~\ref{d:PEPS}.

\begin{cor}\label{c:PESFNotPE}
For every~$\theta>0$, every~$\nu\in\msP(\Delta^{q-1})$ and every~$j\in [q]$,
\begin{equation}\label{eq:VarDust}
\ExpectP_{\PD{\theta,\nu}}\tbraket{\Dust_j}=\mom{\mbfe_j}{\nu}\comma \qquad
\Var_{\PD{\theta,\nu}}\tbraket{\Dust_j}=\frac{\Var_\nu[u_j]}{\theta+1}\fstop
\end{equation}
Consequently,~$\pEsf{\proc}{\theta,\nu}$ is of partially exchangeable type if and only if~$\nu=\delta_\mbfp$ for some~$\mbfp\in\Delta^{q-1}$.
\begin{proof}
By Lemma~\ref{l:PsiColor} and Proposition~\ref{p:RepresentationPESF}, $\int \Dust_j^n\diff\PD{\theta,\nu}=\pEsf{n}{\theta,\nu}\tbraket{\mcA_{n\mbfe_j}}$.
For~$n=1$ this is~$\momentum{\mbfe_j}{\theta,\nu}/\theta=\mom{\mbfe_j}{\nu}$. For~$n=2$ we have~$\mcA_{2\mbfe_j}=\set{\car_{2\mbfe_j},2\car_{\mbfe_j}}$, whence by~\eqref{d:eq:PolyESF}
\[
\int_{\KSimplex{q}} \Dust_j^2\diff\PD{\theta,\nu}=\frac{\mom{2\mbfe_j}{\nu}+\theta\,\mom{\mbfe_j}{\nu}^2}{\theta+1}\comma
\]
and~\eqref{eq:VarDust} follows.
Thus~$\DustV$ is $\PD{\theta,\nu}$-a.s.\ constant if and only if~$\Var_\nu[u_j]=0$ for every~$j$, i.e.\ if and only if~$\nu$ is a Dirac mass.
We conclude by Proposition~\ref{p:DiracCriterion}.
\end{proof}
\end{cor}

\subsubsection{Multipartition structures and the refined Ewens Sampling Formula}\label{sss:Strahov}
Recall from Remark~\ref{r:ColoredVsPolychromatic} that multipartitions are polychromatic partitions of a special type: in our language, a $q$-chromatic partition~$A\in\Part{q}{n}$ represents a \emph{multi(ple) partition of~$n$ into~$q$ components} (see \cite[\S2]{Str26}) if and only if
\begin{equation}\label{eq:Multipartition}
A(\mbfa)>0 \quad \implies \quad \mbfa= a\, \mbfe_j \; \text{ for some } \; a\in [n] \comma j\in [q] \fstop
\end{equation}
We denote by~$\PartM{q}{n}\subset \Part{q}{n}$ the subset of all multipartitions of~$n$ into $q$ components.

\paragraph{Multipartition structures}
In~\cite[Dfn.~2.1]{Str24a}, E.~Strahov introduced \emph{multipartition structures} as a natural generalization to multipartitions of Kingman's partition structures.
Since~$\PartM{q}{n}\subset \Part{q}{n}$, and since deleting any element from a multipartition~$B_n\in \PartM{q}{n}$ results in a multipartition~$B_{n-1}\in\PartM{q}{n-1}$, the kernels~\eqref{eq:Kernels} \emph{preserve} multipartitions; that is, whenever~$A\in\PartM{q}{n}$ we have~$\BackwardK^\sym{q}_{\ell,n}(A,C)>0$ only if~$C\in\PartM{q}{\ell}$.
Thus, our notion of consistency (Dfn.~\ref{d:Consistency}) for random polychromatic partitions extends the one for random multipartitions in~\cite[Dfn.~2.1]{Str24a}.
Therefore, multipartition structures in the sense of~\cite{Str24a} are a special type of $\CPS{q}$, and our Representation Theorem~\ref{t:KingmanRepresentation} extends the analogue for multipartitions in~\cite[Thm.~2.2]{Str24a}.
Recalling that extremal~$\CPS{q}$'s are described by their asymptotic block-color frequencies, note that each non-dust block is asymptotically monochromatic if and only if its frequency has only one non-zero entry, while dust blocks are always singletons and therefore trivially monochromatic. 
Translating~\cite{Str24a} into our language, we thus have the following characterization of multipartition structures.

\begin{prop}\label{p:RepresentationStrahov}
Let~$\probP_\proc\in\CPS{q}$ with representing measure~$\mu\in\msP(\KSimplex{q})$ as in
Theorem~\ref{t:KingmanRepresentation}.
Then, the following are equivalent:
\begin{enumerate}[$(i)$]
\item\label{i:p:RepresentationStrahov:1} $\probP_\proc$ is a multipartition structure (i.e., it is concentrated on~$\PartM{q}{\cup} \coloneqq \bigcup_{n \in \N_0} \PartM{q}{n}$);
\item\label{i:p:RepresentationStrahov:2} $\mbfx_i$ has only one non-zero entry for every~$i\geq 1$ for $\mu$-a.e.~$\mbfx\eqdef \seq{\mbfx_i}_{i\geq 0}\in\KSimplex{q}$.
\end{enumerate}
\end{prop}

\paragraph{Refined Ewens Sampling Formula}
In the recent work~\cite{Str26}, E.~Strahov considered a \emph{refined Ewens Sampling Formula with mutation rate\emph{s}~$\boldtheta\eqdef\seq{\theta_1,\dotsc,\theta_q}\in [0,\infty)^q\setminus\set{\zero}$} (rESF) on multipartitions.
On~$\PartM{q}{n}$, Strahov's rESF~\cite[Eqn.~(2.6)]{Str26} reads
\begin{equation}\label{eq:rESF}
\rEsf{n}{\boldtheta}[A] = \frac{n!}{\Poch{\length{\boldtheta}}{n}} \prod_{j=1}^q \prod_{a=1}^n \paren{\frac{\theta_j}{a}}^{\lambda_a^\sym{j}} \frac{1}{\lambda_a^\sym{j}!} \comma \quad \lambda_a^\sym{j}\eqdef A(a\, \mbfe_j)\comma \qquad A\in\PartM{q}{n}\fstop
\end{equation}

Let~$V_q\eqdef\set{\mbfe_1,\dotsc, \mbfe_q}$ be the set of vertices of~$\Delta^{q-1}$.
The pESF is an extension of the rESF.

\begin{prop}\label{p:StrahovReduction}
$\pEsf{n}{\theta,\nu}[\PartM{q}{n}]=1$ for all~$n\in\N_1$ if and only if~$\nu[V_q]=1$.
If so, %
\[
\pEsf{n}{\theta,\nu_{\boldtheta}}=\rEsf{n}{\boldtheta} \quad \text{with} \quad \boldtheta\in [0,\infty)^q \setminus \set{\zero} \comma \quad \theta\eqdef\length{\boldtheta}\comma \quad \nu_{\boldtheta}\eqdef \theta^{-1}\sum_{j=1}^q \theta_j\delta_{\mbfe_j} \fstop
\]
\begin{proof}
Assume that~$\nu[V_q]=1$, i.e.,~$\nu = \nu_{\boldtheta}$ for some~$\boldtheta$.
With the convention~$0^0=1$, we have~$\mbfe_j^\mbfa=1$ if and only if~$\mbfa= a\mbfe_j$ for some~$a\geq 1$.
Thus~$\mom{\mbfa}{\nu}=\theta^{-1}\theta_j$ if~$\mbfa=a \mbfe_j$ with~$a\geq 1$ and~$\mom{\mbfa}{\nu}=0$ otherwise.
For all such~$\mbfa$'s we also have~$\binom{\length{\mbfa}}{\mbfa}=1$.
Substituting these in~\eqref{d:eq:PolyESF} yields~\eqref{eq:rESF}.

Conversely, assume~$\pEsf{n}{\theta,\nu}[\PartM{q}{n}]=1$ for every~$n\in\N_1$.
Since we have in particular~$\pEsf{2}{\theta,\nu}[\car_{\mbfe_j+\mbfe_k}]=0$ for~$j \neq k$, we get~$\mom{\mbfe_j+\mbfe_k}{\nu}= \int u_ju_k\diff\nu(\mbfu)=0$, which forces $\nu$-a.e.~$\mbfu$ to have a single non-zero coordinate.
\end{proof}
\end{prop}

Together with Theorem~\ref{t:CharacterizationNew}, Proposition~\ref{p:StrahovReduction} readily provides the following characterization of the rESF, which does not appear in~\cite{Str24a,Str26}.

\begin{prop}[rESF characterization]\label{p:CharacterizationRESF}
Let~$\probP_\bullet$ be a multipartition structure %
with $\probP_n[A]>0$ for every~$A\in\PartM{q}{n}$ for every~$n\geq 1$, and satisfying
\begin{equation}
a \lambda_a^\sym{j} \probP_n[A] = b(n,a,j) \probP_{n-a} [A-\car_{a\mbfe_j}] \comma \qquad A\in\PartM{q}{n}\comma n\in\N_1\comma j\in [q]\comma
\end{equation}
with~$\lambda_a^\sym{j} \eqdef A(a\mbfe_j)\geq 1$ for some constant~$b(n,a,j)$ independent of~$A$.
Then,~$\probP_\bullet=\rEsf{\bullet}{\boldtheta}$ for a unique~$\boldtheta\in (0,\infty)^q$.
Furthermore,~$b(n,a,j)=b_*(n,a\mbfe_j)$ as in~\eqref{eq:KingmanDefinitionB}.
\end{prop}

\subsubsection{Further examples}\label{sss:FurtherExamples}
We collect some further~$\CPS{q}$'s.

\begin{lem}[Color decoration]\label{l:Decoration}
Let~$\probP^\sym{1}_\proc$ be a partition structure represented by~$\mu^\sym{1}\in\msP(\KSimplex{1})$ in the sense of~\cite{Kin78}, and let~$\nu\in\msP(\Delta^{q-1})$.
Let~$\mu\in\msP(\KSimplex{q})$ be the law of
\[
\mbfx_0\eqdef \ell_0\, \Bary{\nu} \comma \qquad \mbfx_i\eqdef \ell_i\, \mbfp_i \comma \quad i\geq 1\comma
\]
where~$\boldell\eqdef\seq{\ell_i}_{i\geq 0}\sim\mu^\sym{1}$ and~$\seq{\mbfp_i}_{i\geq 1}$ are i.i.d.~$\nu$-dis\-tri\-buted and independent of~$\boldell$.
Then, the~$\CPS{q}$ represented by~$\mu$ in the sense of Theorem~\ref{t:KingmanRepresentation} is
\begin{equation}\label{eq:l:Decoration:0}
\probP_n[A]= \probP^\sym{1}_n[\boldlambda]\ \frac{\prod_{i=1}^n \lambda_i!}{\prod_{\mbfa\in\supp A} A(\mbfa)!}\ \prod_{\mbfa\in \supp A} \paren{\binom{\length{\mbfa}}{\mbfa} \mom{\mbfa}{\nu}}^{A(\mbfa)} \comma \qquad \boldlambda\eqdef \shape(A) \fstop
\end{equation}
In particular~$\shape_\pfwd\probP_\proc=\probP^\sym{1}_\proc$.

\begin{proof}
The proof of Proposition~\ref{p:RepresentationPESF} only uses that the marks~$\seq{\mbfp_i}_i$ are i.i.d.\ and independent of~$\boldell$, together with the fact that~$\mu^\sym{1}$ represents~$\probP^\sym{1}_\proc$; it therefore applies \emph{verbatim}.
The last assertion follows summing~\eqref{eq:l:Decoration:0} over all~$A$ with~$\shape(A)=\boldlambda$ and applying, for each block size~$i$, the Multinomial Theorem together with~$\sum_{\length{\mbfa}=i}\binom{i}{\mbfa}\mom{\mbfa}{\nu}=\int_{\Delta^{q-1}}\length{\mbfu}^i\diff\nu(\mbfu)=1$, exactly as in Proposition~\ref{p:ESFProbab}.
\end{proof}
\end{lem}

Choosing~$\mu^\sym{1}=\PD{\theta}$ in Lemma~\ref{l:Decoration} returns the pESF, cf.~Proposition~\ref{p:RepresentationPESF}.
It is clear that this construction has limited scope: the color composition of the dust is always the barycenter~$\Bary{\nu}$, since in the paintbox~\eqref{eq:xir} each dust element forms its own block and thus acquires color~$j$ with probability~$\mom{\mbfe_j}{\nu}$; and the marks are independent of the block frequencies.

\begin{ese}[Polychromatic Pitman--Yor]\label{ese:PitmanYor}
Let~$0\leq \alpha<1$ and~$\theta>-\alpha$, and let~$\PD{\alpha,\theta}$ be the two-parameter Poisson--Dirichlet distribution~\cite{PitYor97}.
Further fix~$\nu\in\msP(\Delta^{q-1})$ with~$\nu\mathring{\Delta}^{q-1}>0$.
Decorating it as in Lemma~\ref{l:Decoration} and using~$\binom{\length{\mbfa}}{\mbfa}\big/\length{\mbfa}!=1/\mbfa!$ gives
\begin{equation}\label{eq:ese:PitmanYor:0}
\probP_n[A]=\frac{n!\, \displaystyle\prod_{j=1}^{\size{A}-1}(\theta+j\alpha)}{\Poch{\theta+1}{n-1}} \prod_{\mbfa\in\supp A} \frac{1}{A(\mbfa)!}\paren{\frac{\Poch{1-\alpha}{\length{\mbfa}-1}\, \mom{\mbfa}{\nu}}{\mbfa!}}^{A(\mbfa)} \fstop
\end{equation}
For~$\alpha=0$ this is~\eqref{d:eq:PolyESF}; for~$\theta=0$ and~$\alpha\in(0,1)$ it is the polychromatic `normalized $\alpha$-stable structure', with $\prod_{j<\size{A}}(j\alpha)/\Poch{1}{n-1}=\alpha^{\size{A}-1}(\size{A}-1)!/(n-1)!$.
The corresponding urn model is obtained from~\eqref{eq:HoppeKernel} by deforming the block-\emph{size} dynamics alone, the color mechanism being unchanged:
\begin{align*}
\ForwardK^\sym{q}_{m,m+1}\tparen{B, B+\car_{\mbfe_j}}&= \frac{\theta+\alpha\,\size{B}}{\theta+m}\, \mom{\mbfe_j}{\nu} \comma
\\
\ForwardK^\sym{q}_{m,m+1}\tparen{B, B-\car_\mbfb+\car_{\mbfb+\mbfe_j}}&= \frac{\tparen{\length{\mbfb}-\alpha} B(\mbfb)}{\theta+m}\, \frac{\mom{\mbfb+\mbfe_j}{\nu}}{\mom{\mbfb}{\nu}} \comma \qquad \mbfb\in\supp B \fstop
\end{align*}
\end{ese}

\begin{ese}[Dirichlet colors]\label{ese:DirichletColors}
Let~$\boldalpha\in(0,\infty)^q$, %
so that
\[
\mom{\mbfa}{\Dir{\boldalpha}} = \frac{\Poch{\boldalpha}{\mbfa}}{\Poch{\length{\boldalpha}}{\length{\mbfa}}} \comma \qquad \mbfa\in\N_0^q \fstop
\]
Choosing~$\nu\eqdef \Dir{\boldalpha}$ in the pESF gives
\begin{equation}\label{eq:ese:DirichletColors:0}
\probP_n[A]= \frac{n!\,\theta^{\size{A}}}{\Poch{\theta}{n}} \prod_{\mbfa\in\supp A}\frac{1}{A(\mbfa)!} \paren{\frac{\tparen{\length{\mbfa}-1}!}{\mbfa!}\, \frac{\Poch{\boldalpha}{\mbfa}}{\Poch{\length{\boldalpha}}{\length{\mbfa}}}}^{A(\mbfa)} \fstop
\end{equation}
Since~$\mom{\mbfb+\mbfe_j}{\nu}\big/\mom{\mbfb}{\nu}=(\alpha_j+b_j)\big/(\length{\boldalpha}+\length{\mbfb})$, the kernels~\eqref{eq:HoppeKernel} read
\[
\ForwardK^\sym{q}_{m,m+1}\ttparen{B, B+\car_{\mbfe_j}}= \frac{\theta}{\theta+m}\, \frac{\alpha_j}{\length{\boldalpha}} \comma \qquad
\ForwardK^\sym{q}_{m,m+1}\ttparen{B, B-\car_\mbfb+\car_{\mbfb+\mbfe_j}}= \frac{\length{\mbfb} B(\mbfb)}{\theta+m}\, \frac{\alpha_j+b_j}{\length{\boldalpha}+\length{\mbfb}} \fstop
\]
In other words, blocks are governed by a standard Hoppe urn with parameter~$\theta$, while for each block an independent P\'olya urn is run with initial composition~$\boldalpha$ on the~$q$ colors.
The family~\eqref{eq:ese:DirichletColors:0} interpolates between the two distinguished subclasses of~\S\ref{sss:PEPS} and~\S\ref{sss:Strahov}.
Indeed, for fixed~$\mbfp\eqdef \boldalpha/\length{\boldalpha}\in\Delta^{q-1}$,
\[
\lim_{\beta\downarrow 0}\ \Dir{\beta\mbfp}=\sum_{j=1}^q p_j\delta_{\mbfe_j} \qquad \text{and} \qquad \lim_{\beta\uparrow\infty}\ \Dir{\beta \mbfp}=\delta_\mbfp \comma
\]
weakly, %
so that~$\probP_\proc$ converges to Strahov's rESF~$\rEsf{\proc}{\theta\mbfp}$ of Proposition~\ref{p:StrahovReduction} as~$\beta\downarrow0$, and to the pESF~$\pEsf{\proc}{\theta,\delta_\mbfp}$ as~$\beta\uparrow\infty$.
The parameter~$\length{\boldalpha}$ thus measures how far each block is from being monochromatic.
\end{ese}

\begin{ese}[A single color urn shared by all blocks]\label{ese:GlobalPolya}
Let~$\boldalpha\in(0,\infty)^q$ and let~$\mu\eqdef\int_{\Delta^{q-1}} \PD{\theta,\delta_\mbfp}\diff\Dir{\boldalpha}(\mbfp)$, with~$\PD{\theta,\nu}$ as in Definition~\ref{d:polyPD}.
In other words, all blocks share the same color composition~$\mbfp$, itself~$\Dir{\boldalpha}$-distributed.
The corresponding~$\CPS{q}$ is
\begin{equation}\label{eq:ese:GlobalPolya:0}
\probP_n[A] = \binom{n}{\col{A}} \frac{\Poch{\boldalpha}{\col{A}}}{\Poch{\length{\boldalpha}}{n}} \cdot \frac{\theta^{\size{A}}\, \Multi{\col{A}}(A)}{\Poch{\theta}{n}} \fstop
\end{equation}
The first factor is the Dirichlet--multinomial law of the color counts~$\col{A}$.
By Corollary~\ref{c:Multi} the second one is a probability on~$\mcA_{\col{A}}$, and it is the law of the coloring of a $\theta$-biased random permutation of~$[n]$.
In words, we color the~$n$ elements by a single P\'olya urn, independently of everything else, then read off the coloring of a $\theta$-biased random permutation.
Here~$\DustV_\pfwd\mu=\Dir{\boldalpha}$.
In particular~$\mu$ is not a Dirac mass, so~$\probP_\proc$ is \emph{neither} extremal %
nor of partially exchangeable type (by Proposition~\ref{p:DiracCriterion}) although every~$\pEsf{\proc}{\theta,\delta_\mbfp}$ in the mixture is of partially exchangeable type.
\end{ese}

\begin{ese}[Unmarked dust]\label{ese:Dust}
Fix~$s\in (0,1)$, $\theta>0$, $\nu\in\msP(\Delta^{q-1})$, and let~$\mu$ be the law of
\[
\mbfx_0\eqdef (1-s)\, \mbfe_1\comma \qquad \seq{\mbfx_i}_{i\geq 1}\eqdef s\, \seq{\mbfy_i}_{i\geq 1} \comma \qquad \mbfy\sim \PD{\theta,\nu}\fstop
\]
Since~$\varphi_B$ is positively homogeneous of degree~$\length{\col{B}}$ in the variables~$\seq{\mbfx_i}_{i\geq 1}$, the definition~\eqref{eq:DefPsi} of~$\psi_A$ and Proposition~\ref{p:RepresentationPESF} give
\begin{equation}\label{eq:ese:Dust:0}
\probP_n[A]= \sum_{m=0}^{A(\mbfe_1)} \binom{n}{m} (1-s)^m s^{\,n-m}\ \pEsf{n-m}{\theta,\nu}\bracket{A- m\, \car_{\mbfe_1}} \fstop
\end{equation}
Equivalently, independently of everything else, each of the~$n$ elements is with probability~$1-s$ a singleton of color~$1$, and with probability~$s$ it is sampled from~$\pEsf{\proc}{\theta,\nu}$.
Whenever~$\Bary{\nu}\neq \mbfe_1$, the measure~$\mu$ is \emph{not} of the same form as in Lemma~\ref{l:Decoration}, since the direction of the dust vector~$\mbfx_0$ differs from that of the blocks.
This is the simplest~$\CPS{q}$ exhibiting the vectorial nature of the dust, void when~$q=1$.
(In the interpretation of Appendix~\ref{ss:PopGen} below, new allelic types arise unmarked, whereas markers accumulate only along established lineages.)
\end{ese}

\section{Multivariate Dirichlet and Gamma moments}\label{s:ApplMoments}
We present applications of polychromatic permutations to multivariate moments of some well-known distributions.

Fix~$k\in\N_1$.
For vectors~$\mbfs_1,\dotsc,\mbfs_q\in\R^k$, let~$\mbfS\in\R^{k\times q}$ be the matrix with columns~$\mbfs_1,\dotsc,\mbfs_q$.

\subsection{Multivariate moments of Dirichlet measures}\label{sss:Moments}
Finding useful representations for the moments of the Dirichlet distribution~$\Dir{\boldalpha}$ is a difficult problem.
We present a brief historical account in~\S\ref{ss:MomentsOverview}.
As an application of our computations on polychromatic partitions, we provide an %
elementary, closed formula for all multivariate moments of~$\Dir{\boldalpha}$.
Fix integers~$q\in \N_1$ and~$\mbfn\eqdef \seq{n_1,\dotsc, n_q}\in\N^q_0$, let~$\msZ^\sym{q}_\mbfn$ be the \emph{pattern inventory}~\eqref{eq:PatternInventory}, $Z^\sym{q}_\mbfn(\mbft)$ be the $q$-chromatic cycle index polynomial~\eqref{eq:ZDefinition}, and recall the relation~\eqref{eq:t:Zn:0} between them.

\begin{thm}[Moments of~$\Dir{\boldalpha}$]\label{t:Moments}
{Fix~$k\in\N_1$,~$\mbfs_1,\dotsc,\mbfs_q\in\R^k$, and~$\boldalpha \in (0,\infty)^k$.} Then
\begin{equation}\label{eq:t:Moments:0}
\int_{\Delta^{k-1}} \prod_{j=1}^q (\mbfs_j\cdot\mbfx)^{n_j} \diff\Dir{\boldalpha}(\mbfx)= \frac{\length{\mbfn}!}{\Poch{\length{\boldalpha}}{\length{\mbfn}}}\, \msZ^\sym{q}_\mbfn[\mbfS;\boldalpha] 
\fstop
\end{equation}
\end{thm}

By `elementary' we mean that our formula is expressed only in terms of elementary functions, and by `closed' that it is both non-recursive and non-iterative.

In order to prove Theorem~\ref{t:Moments}, we start with some straightforward preparatory lemmas.

\begin{lem}
For every~$k\in \N_1$, every~$\mbfv\in \N_0^k$, and every integer~$0\leq m \leq \length{\mbfv}$,
\begin{equation}\label{eq:l:Multinomial:0}
\binom{\length{\mbfv}}{\mbfv}=\sum_{\substack{\mbfw\in \N_0^k\\ \mbfw\leq_\had \mbfv\comma \length{\mbfw}=m}} \binom{\length{\mbfw}}{\mbfw} \binom{\length{\mbfv}-\length{\mbfw}}{\mbfv-\mbfw} \fstop
\end{equation}

\begin{proof}
	See, e.g.,~\cite[Eqn.~(8), p.~1060]{Tau63},~\cite[Eqn.~(10), p.~39]{Car63}.
\end{proof}
\end{lem}

\begin{lem}\label{l:AuxiliaryNu}
The following identity holds
\begin{equation}\label{eq:l:Moments:0}
\mu_\mbfn[\mbfS;\boldalpha] \eqdef \int_{\Delta^{k-1}} \prod_{j=1}^q (\mbfs_j\cdot\mbfx)^{n_j} \diff\Dir{\boldalpha}(\mbfx) = \frac{\mbfn!}{\Poch{\length{\boldalpha}}{\length{\mbfn}}} \sum_{\substack{\mbfM\in \N_0^{k\times q}\\ \car^\transpose \mbfM=\mbfn}} \Poch{\boldalpha}{\mbfM\car} \frac{\mbfS^\mbfM}{\mbfM!}\defeq \nu_\mbfn[\mbfS;\boldalpha] \fstop
\end{equation}

\begin{proof}
	The proof follows by expanding the powers~$(\mbfs_j \cdot \mbfx)^{n_j}$ with the Multinomial Theorem and by the properties of the Dirichlet distribution.
\end{proof}
\end{lem}

\begin{proof}[Proof of Theorem~\ref{t:Moments}]
Our proof follows the same double induction argument as in the proof in the univariate case~\cite[Thm.~3.3]{LzDS19a}.
Let~$\nu_\mbfn[\mbfS,\boldalpha]$ be as in~\eqref{eq:l:Moments:0}, $\zeta_\mbfn[\mbfS;\boldalpha]$ be defined by the right-hand side of~\eqref{eq:t:Moments:0}, and set
\[
\tilde \nu_\mbfn[\mbfS;\boldalpha]\eqdef \frac{\Poch{\length{\boldalpha}}{\length{\mbfn}}}{\mbfn!} \nu_\mbfn[\mbfS;\boldalpha]
\qquad \text{and} \qquad 
\tilde \zeta_\mbfn[\mbfS;\boldalpha]\eqdef \frac{\Poch{\length{\boldalpha}}{\length{\mbfn}}}{\mbfn!} \zeta_\mbfn[\mbfS;\boldalpha] \fstop
\]
By showing that~$\tilde\nu_\mbfn[\mbfS;\boldalpha]=\tilde\zeta_\mbfn[\mbfS;\boldalpha]$, one concludes the assertion by Lemma~\ref{l:AuxiliaryNu}.

\paragraph{Step~1} We claim that
\begin{align}\label{eq:t:Moments:1}
\tilde\nu_{\mbfn-\mbfe_j}[\mbfS;\boldalpha+\mbfe_\ell]
=
\sum_{\mbfe_j\leq_\had\mbfh\leq_\had\mbfn} \mbfS_\ell^{\mbfh-\mbfe_j} \frac{(\length{\mbfh}-1)!}{(\mbfh-\mbfe_j)!} \, \tilde\nu_{\mbfn-\mbfh}[\mbfS;\boldalpha] %
\comma
\end{align}
where, conventionally,
\begin{equation}\label{eq:t:Moments:2}
\tilde\nu_\mbfn[\mbfS;\boldalpha]=0 \quad \text{whenever} \quad \mbfn\not\geq_\had \zero\fstop
\end{equation}
Note that~$\tilde\nu_{\boldsymbol 0}[\mbfS;\boldalpha] = 1$.
We argue by induction on~$\length{\mbfn}$ {simultaneously for every~$\boldalpha \in (0,\infty)^k$. The base step~$\length{\mbfn}=1$ is trivial, because in this case~\eqref{eq:t:Moments:1} reads either~$1=1$ (if~$\mbfn = \mbfe_j$), or~$0=0$ (if~$n_j=0$).}%

\emph{Inductive step.} {If~$n_j = 0$, there is nothing to prove. We thus assume~$n_j \ge 1$.} Let~$\partial_a^b\eqdef \partial_{s_a^b}$, set~$\mbfE_a^b\eqdef [\delta_{ai}\delta_{bj}]_i^j\in \set{0,1}^{k\times q}$, and note that
\begin{align}
\nonumber
\partial_a^b \tilde\nu_\mbfn[\mbfS;\boldalpha]=&\ \sum_{\substack{\mbfM\in \N_0^{k\times q} \\ \car^\transpose\mbfM=\mbfn}} \Poch{\boldalpha}{\mbfM\car} \partial_a^b  \frac{\mbfS^\mbfM}{\mbfM!}= \sum_{\substack{\mbfE_a^b\leq_\had\mbfM\in \N_0^{k\times q} \\ \car^\transpose\mbfM=\mbfn}} \alpha_a \Poch{\boldalpha+\mbfe_a}{(\mbfM-\mbfE_a^b)\car} \frac{\mbfS^{\mbfM-\mbfE_a^b}}{(\mbfM-\mbfE_a^b)!}
\\
\nonumber
=&\ \sum_{\substack{\mbfM\in \N_0^{k\times q} \\ \car^\transpose\mbfM=\mbfn-\mbfe_b}} \alpha_a \Poch{\boldalpha+\mbfe_a}{\mbfM\car} \frac{\mbfS^{\mbfM}}{\mbfM!}
\\
\label{eq:t:Moments:3}
=&\ \alpha_a \tilde\nu_{\mbfn-\mbfe_b}[\mbfS;\boldalpha+\mbfe_a] \fstop
\end{align}

Applying the inductive hypothesis to~$\mbfn-\mbfe_b$ with~$\boldalpha+\mbfe_a$ in place of~$\boldalpha$ {(or trivially if~$n_b=0$),}
\begin{align}
\nonumber
\alpha_a \tilde\nu_{\mbfn-\mbfe_j-\mbfe_b}[\mbfS;\boldalpha+\mbfe_\ell+\mbfe_a]=&\ \alpha_a \sum_{\mbfe_j\leq_\had\mbfh\leq_\had \mbfn-\mbfe_b} \mbfS_\ell^{\mbfh-\mbfe_j} \frac{(\length{\mbfh}-1)!}{(\mbfh-\mbfe_j)!} \, \tilde\nu_{\mbfn-\mbfh-\mbfe_b}[\mbfS;\boldalpha+\mbfe_a]
\\
\label{eq:t:Moments:3.5}
=&\ \alpha_a \sum_{\mbfe_j\leq_\had \mbfh\leq_\had\mbfn} \mbfS_\ell^{\mbfh-\mbfe_j} \frac{(\length{\mbfh}-1)!}{(\mbfh-\mbfe_j)!}\, \tilde\nu_{\mbfn-\mbfh-\mbfe_b}[\mbfS;\boldalpha+\mbfe_a]\comma
\end{align}
where the last equality holds by~\eqref{eq:t:Moments:2}.

Now, let~$\ell\neq a$.
Applying~\eqref{eq:t:Moments:3} with~$\mbfn-\mbfe_j$ in place of~$\mbfn$ and~$\boldalpha+\mbfe_\ell$ in place of~$\boldalpha$,
\begin{align}\label{eq:t:Moments:3.7}
\partial_a^b \tilde\nu_{\mbfn-\mbfe_j}[\mbfS;\boldalpha+\mbfe_\ell] = \alpha_a \tilde\nu_{\mbfn-\mbfe_j-\mbfe_b}[\mbfS;\boldalpha+\mbfe_\ell+\mbfe_a] \fstop
\end{align}
Since~$\ell\neq a$, applying~\eqref{eq:t:Moments:3} with~$\mbfn-\mbfh$ in place of~$\mbfn$ for every~$\mbfh \leq_\had \mbfn$ yields
\begin{equation}\label{eq:t:Moments:3.8}
\begin{aligned}
\partial_a^b& \bracket{\sum_{\mbfe_j\leq_\had\mbfh\leq_\had\mbfn} \mbfS_\ell^{\mbfh-\mbfe_j} \frac{(\length{\mbfh}-1)!}{(\mbfh-\mbfe_j)!} \tilde\nu_{\mbfn-\mbfh}[\mbfS;\boldalpha]}=
\\
&\ =\alpha_a \sum_{\mbfe_j\leq_\had\mbfh\leq_\had\mbfn} \mbfS_\ell^{\mbfh-\mbfe_j} \frac{(\length{\mbfh}-1)!}{(\mbfh-\mbfe_j)!}\, \tilde\nu_{\mbfn-\mbfh-\mbfe_b}[\mbfS;\boldalpha+\mbfe_a] \fstop
\end{aligned}
\end{equation}

Combining~\eqref{eq:t:Moments:3.7},~\eqref{eq:t:Moments:3.8}, and~\eqref{eq:t:Moments:3.5} yields
\begin{align*}
\partial_a^b \paren{\tilde\nu_{\mbfn-\mbfe_j}[\mbfS;\boldalpha+\mbfe_\ell] - \sum_{\mbfe_j\leq_\had\mbfh\leq_\had\mbfn} \mbfS_\ell^{\mbfh-\mbfe_j} \frac{(\length{\mbfh}-1)!}{(\mbfh-\mbfe_j)!} \tilde\nu_{\mbfn-\mbfh}[\mbfS;\boldalpha]}=0 \fstop
\end{align*}

By arbitrariness of~$a,b$, we conclude that the bracketed quantity is a polynomial in the sole variables~$s_\ell^1,\dotsc, s_\ell^q$.
As a consequence, every monomial not in the sole variables~$s_\ell^1,\dotsc, s_\ell^q$ cancels out by arbitrariness of the variables~$\mbfS$.
Thus,
\begin{align*}
\tilde\nu_{\mbfn-\mbfe_j}&[\mbfS;\boldalpha+\mbfe_\ell] - \sum_{\mbfe_j\leq_\had\mbfh\leq_\had\mbfn} \mbfS_\ell^{\mbfh-\mbfe_j} \frac{(\length{\mbfh}-1)!}{(\mbfh-\mbfe_j)!} \tilde\nu_{\mbfn-\mbfh}[\mbfS;\boldalpha]
\\
&= \Poch{\alpha_\ell+1}{\length{\mbfn}-1}\frac{\mbfS_\ell^{\mbfn-\mbfe_j}}{(\mbfn-\mbfe_j)!}-\sum_{\mbfe_j\leq_\had\mbfh\leq_\had\mbfn} \mbfS_\ell^{\mbfh-\mbfe_j} \frac{(\length{\mbfh}-1)!}{(\mbfh-\mbfe_j)!} \Poch{\alpha_\ell}{\length{\mbfn}-\length{\mbfh}}\frac{\mbfS_\ell^{\mbfn-\mbfh}}{(\mbfn-\mbfh)!} \fstop
\end{align*}
The latter quantity is proved to vanish as soon as
\begin{align*}
\frac{\Poch{\alpha_\ell+1}{\length{\mbfn}-1}}{(\mbfn-\mbfe_j)!}=\sum_{\mbfe_j\leq_\had\mbfh\leq_\had\mbfn} \frac{(\length{\mbfh}-1)!}{(\mbfh-\mbfe_j)!} \frac{\Poch{\alpha_\ell}{\length{\mbfn}-\length{\mbfh}}}{(\mbfn-\mbfh)!}\comma
\end{align*}
or, equivalently, by the Chu--Vandermonde identity,
\begin{align*}
\frac{(\length{\mbfn}-1)!}{(\mbfn-\mbfe_j)!}\sum_{i=1}^{\length{\mbfn}} \frac{\Poch{\alpha_\ell}{\length{\mbfn}-i}}{(\length{\mbfn}-i)!} =\sum_{i=1}^{\length{\mbfn}}\sum_{\substack{\mbfe_j\leq_\had\mbfh\leq_\had\mbfn\\ \length{\mbfh}=i}}\frac{(\length{\mbfh}-1)!}{(\mbfh-\mbfe_j)!} \frac{\Poch{\alpha_\ell}{\length{\mbfn}-i}}{(\mbfn-\mbfh)!}\fstop
\end{align*}
The latter is implied by the equality of each of the summands, viz.
\begin{align*}
\frac{(\length{\mbfn}-1)!}{(\mbfn-\mbfe_j)!}\frac{1}{(\length{\mbfn}-i)!} =\sum_{\substack{{\mbfe_j\leq_\had\mbfh\leq_\had\mbfn}\\ \length{\mbfh}=i}} \frac{(\length{\mbfh}-1)!}{(\mbfh-\mbfe_j)!} \frac{1}{(\mbfn-\mbfh)!}\comma
\end{align*}
which is in turn a consequence of %
\eqref{eq:l:Multinomial:0}, after relabeling~$\mbfn$ as~$\mbfn-\mbfe_j$.

\paragraph{Step~2}
We now verify that~$\tilde\nu_\mbfn=\tilde\zeta_\mbfn$.
We argue by strong induction on~$\length{\mbfn}$ with trivial (i.e.~$1=1$) base step~$\length{\mbfn}=0$. 
\emph{Inductive step.} Assume for every~$\boldalpha\in (0,\infty)^k$ %
that~$\tilde\nu_{\mbfn-\mbfh}[\mbfS;\boldalpha]=\tilde\zeta_{\mbfn-\mbfh}[\mbfS;\boldalpha]$ for every~$\mbfh\leq_\had\mbfn$ with~$\mbfh\neq \zero$.
Now, recalling~\eqref{eq:Pattern} {and~\eqref{eq:t:Zn:0}},
\begin{align}
\nonumber
{\mbfn!} \, \partial_p^b \tilde\zeta_\mbfn[\mbfS;\boldalpha]=&\ \sum_{A\vDash \mbfn} \Multi{\mbfn}(A)\, \partial_p^b \prod_{\mbfa\in\supp{A}} \paren{ \tparen{\mbfs_1^{\had a_1}\had\cdots \had \mbfs_q^{\had a_q}}\cdot\boldalpha }^{A(\mbfa)}%
\\
\nonumber
=&\ {\sum_{A\vDash \mbfn} \Multi{\mbfn}(A)\, \sum_{\substack{\mbfh\in\supp{A}\\ \mbfh\geq_\had\mbfe_b}} A(\mbfh)\, h_b\, \alpha_p \,\mbfS_p^{\mbfh-\mbfe_b} \, w[\mbfS;\boldalpha](A - \car_\mbfh)}
\\
\label{eq:t:Moments:4}
=&\ {\alpha_p\sum_{\mbfe_b\leq_\had\mbfh\leq_\had\mbfn} \sum_{\substack{A\vDash \mbfn \\ A \ge \car_{\mbfh}}} \Multi{\mbfn}(A) A(\mbfh)\, h_b \,\mbfS_p^{\mbfh-\mbfe_b} \, w[\mbfS;\boldalpha](A-\car_\mbfh)
 \fstop}
\end{align}
For each~${\mbfe_b \leq_\had} \mbfh\leq_\had\mbfn$ and each~$A\vDash\mbfn$ with~$A\geq \car_\mbfh$, set~$C\eqdef A-\car_\mbfh$.
Note that%
\begin{align}
\nonumber
\Multi{\mbfn}(A)=&\ \mbfn! \prod_{\mbfa\in \supp{A}} \frac{\binom{\length{\mbfa}}{\mbfa}^{A(\mbfa)}}{\length{\mbfa}^{A(\mbfa)}A(\mbfa)!} 
= \frac{\mbfn!\, \length{\mbfh}!}{\length{\mbfh}\, \mbfh!\, A(\mbfh)}\prod_{\mbfa\in \supp{C}} \frac{\binom{\length{\mbfa}}{\mbfa}^{C(\mbfa)}}{\length{\mbfa}^{C(\mbfa)}C(\mbfa)!} 
\\
\label{eq:t:Moments:5}
=&\ \frac{\Multi{\mbfn-\mbfh}(C)}{A(\mbfh)} \frac{(\length{\mbfh}-1)!}{\mbfh!} \frac{\mbfn!}{(\mbfn-\mbfh)!} \fstop
\end{align}
Substituting~\eqref{eq:t:Moments:5} %
 in~\eqref{eq:t:Moments:4} above, and simplifying~$\mbfn!$,
\begin{align*}
\partial_p^b\, \tilde\zeta_\mbfn[\mbfS;\boldalpha]=&\ \alpha_p \sum_{\mbfe_b\leq_\had\mbfh\leq_\had\mbfn} \sum_{C\vDash \mbfn-\mbfh} \frac{\Multi{\mbfn-\mbfh}(C)}{(\mbfn-\mbfh)!} \frac{(\length{\mbfh}-1)!}{(\mbfh-\mbfe_b)!}\ \mbfS_p^{\mbfh-\mbfe_b} \,w[\mbfS;\boldalpha](C)
\\
=&\ \alpha_p \sum_{\mbfe_b\leq_\had\mbfh\leq_\had\mbfn}  \frac{(\length{\mbfh}-1)!}{(\mbfh-\mbfe_b)!}\ \mbfS_p^{\mbfh-\mbfe_b} \,\tilde\zeta_{\mbfn-\mbfh}[\mbfS;\boldalpha] \fstop
\end{align*}

Combining the inductive hypothesis with~\eqref{eq:t:Moments:1} with~$j=b$ and~$\ell=p$, and~\eqref{eq:t:Moments:3} with~$a=p$,
\begin{align*}
\partial_p^b\, \tilde\zeta_\mbfn[\mbfS;\boldalpha]=&\ \alpha_p \sum_{\mbfe_b\leq_\had\mbfh\leq_\had\mbfn}  \frac{(\length{\mbfh}-1)!}{(\mbfh-\mbfe_b)!}\ \mbfS_p^{\mbfh-\mbfe_b}\, \tilde\nu_{\mbfn-\mbfh}[\mbfS;\boldalpha] 
= \partial_p^b\, \tilde\nu_\mbfn[\mbfS;\boldalpha]
\fstop
\end{align*}
By arbitrariness of~$p$ and~$b$ we conclude that~$\tilde\zeta_\mbfn[\mbfS;\boldalpha]-\tilde\nu_\mbfn[\mbfS;\boldalpha]$ is constant as a function of~$\mbfS$, hence vanishing by choosing~$\mbfS=\zero$.
\end{proof}

\subsection{Overview on moments of Dirichlet measures}\label{ss:MomentsOverview}
In order to put Theorem~\ref{t:Moments} into context, we briefly survey previously known results on moments of Dirichlet and related measures.
Moment and Laplace/Fourier-transform methods for~$\Dir{\boldalpha}$ and its infinite-dimensional counterpart, the Dirichlet--Ferguson measure~$\DF{\alpha}$~\cite{Fer73} over a measure space~$(X,\alpha)$ %
had been widely considered unfeasible, %
as we briefly summarize below.

\paragraph{Transforms} It is well-known that the Fourier transform~$\widehat{\Dir{\boldalpha}}$ of the Dirichlet distribution~$\Dir{\boldalpha}$ may be expressed in terms of~${}_k\Phi_2$, the $k$-variate \emph{confluent hypergeometric Lauricella function of type~$D$}~\cite{Lau1893, Ext76}.
The power-series representation of~${}_k\Phi_2$ is \guillemotleft\emph{inconvenient for numerical calculations when~$k>2$}\guillemotright~\cite[p.~4]{Phi88}.
Instead, the complex-contour-integral representation~\cite[Eqn.~(7)]{Erd40} is preferred, but its treatment remains quite involved, see e.g.~\cite{RegGugDiN02}.
In particular, differentiating~$\widehat{\Dir{\boldalpha}}$ in this form does not provide any tractable %
representation for the moments. %

For decades, the Fourier transform~$\widehat{\DF{\alpha}}$ of~$\DF{\alpha}$ was widely considered intractable~\cite{JiaDicKuo04}, which led to the introduction of other characterizing transforms, such as the \emph{Markov--Krein transform}~\cite{KerTsi01} and the \emph{$c$-transform}~\cite{JiaDicKuo04}.
These methods too are unsatisfactory, since there is no counterpart for such transforms of foundational results available for the Fourier transform, such as, for instance, Bochner--Minlos--Sazonov Theorem or L\'evy's Continuity Theorem.
The Fourier transform~$\widehat{\DF{\alpha}}$ was eventually computed in~\cite{LzDS19a} by methods in representation theory.

\paragraph{Moments}
Multivariate moments of~$\Dir{\boldalpha}$ are well-known in the form%
~\eqref{eq:l:Moments:0},
which may be regarded as an extension of the ESF.
Whereas easily computable for small~$k$, this form is unsuitable for letting~$k\to\infty$ and thus provides no insight on multivariate moments of~$\DF{\alpha}$.

Partially in order to overcome this issue, other expressions have been considered:
\begin{enumerate*}[$(a)$]
\item the univariate moment~$\DF{\alpha}(f^n)$ has appeared in~\cite{Reg98} in terms of incomplete Bell polynomials, solely in the case $X\Subset \R$ and~$f(t)\eqdef t$; %

\item more general univariate moments for~$\Dir{\boldalpha}$ have implicitly appeared in~\cite[proof of Prop.~3.3]{LetPic18} in iterative form;

\item univariate moments for both~$\Dir{\boldalpha}$ and~$\DF{\alpha}$ have appeared in full generality in~\cite{LzDS19a} in terms of the cycle index polynomials~$Z_n$~\eqref{eq:CycleIndex} of~$\mfS_n$, which allowed the aforementioned computation of~$\widehat{\DF{\alpha}}$.

As for multivariate moments, they have appeared:
\item in~\cite[Prop.~7.4]{KerTsi01}, in terms of summations over constrained permutations, only in the case~$\length{\boldalpha}=1$;

\item in~\cite[Eqn.~(4.20)]{EthKur94},~\cite[Lem.~5.2]{Fen10}, and~\cite[Cor.~3.5]{LzDSLyt17}, in terms of summations over constrained \emph{set} partitions.
\end{enumerate*}

\paragraph{Other measures}
The measure~$\DF{\alpha}$ is the \emph{simplicial part} of other known measures on the space~$\Mbp$ of non-negative Borel measures on~$X$. 
Among them are: the law~$\GP{\alpha}$ of the \emph{$\gamma$-point process}~\cite{KonDaSStrUs98} with intensity~$\alpha$, and A.M.~Vershik's \emph{multiplicative infinite-dimensional Lebesgue measure}~$\LP_\alpha$~\cite{Ver07,Ver08} with intensity~$\alpha$.
Together with~$\DF{\alpha}$, these measures have a wide range of applications, from the theory of point processes and of measure-valued Markov diffusions, see~\cite[\S1]{LzDS17+} and references therein, to the representation theory of infinite-dimensional Lie groups of currents/multipliers, see~\cite{TsiVerYor01}, or~\cite[\S1]{LzDS19a} for a unified treatment.

\subsection{Dirichlet--Ferguson and Gamma measures}\label{ss:DirichletGamma}
Let~$X$ be a second countable locally compact Hausdorff space.
We let~$\Mbp$ be the space of all finite non-negative non-zero Borel measures on~$X$, endowed with the narrow topology and its Borel $\sigma$-algebra, and by~$\msP$ be the subspace of all Borel probability measures on~$X$. %
Note that there is a homeomorphism
\begin{equation}\label{eq:JIsomorphism}
J\colon \Mbp \longrar \msP\times (0,\infty) \comma \qquad J\colon \mu \longmapsto \tparen{\tfrac{\mu}{\mu X}, \mu X}\fstop
\end{equation}

For any finite Borel measure~$\eta$ on~$X$ and any bounded Borel~$f\colon X\to\R$ we set~$\eta f\eqdef \int f\diff\eta$.
For~$\beta>0$ and~$\sigma\in\msP$, let~$\alpha\eqdef \beta\sigma$ be the measure with total mass~$\beta$ and \emph{simplicial part}~$\sigma$. 

\paragraph{Dirichlet--Ferguson measures}
The \emph{Dirichlet--Ferguson measure~$\DF{\alpha}$ with intensity} (\emph{measure})~$\alpha$ is the unique Borel probability measure on~$\msP$ with Fourier transform~\cite[Thm.~3.10]{LzDS19a}
\begin{align*}
\widehat{\DF{\alpha}}(f)\eqdef \int_\msP e^{\imu \, \eta f} \diff\DF{\alpha}(\eta)= \sum_{n=0}^\infty \frac{\imu^n}{\Poch{\beta}{n}} Z_n\tparen{\alpha f, \alpha f^2, \dotsc, \alpha f^n} \comma \qquad f\in\Cb 
\fstop
\end{align*}

For bounded {Borel}~$f_1,\dotsc, f_q\colon X\to\R$, set~$%
\Omega_\mbfn[f_1,\dotsc, f_q;\alpha]\eqdef \seq{\alpha\tparen{f_1^{a_1}\cdots f_q^{a_q}}}_{\mbfa\leq_\had\mbfn} %
$.%

By a straightforward adaptation of the proof for the univariate case~\cite[Thm.~3.10]{LzDS19a}, as a corollary of Theorem~\ref{t:Moments} we obtain an explicit expression for the moments of~$\DF{\alpha}$.

\begin{cor}[Multivariate moments of~$\DF{\alpha}$] \label{c:MomentsDF} We have 
\[
\int_{\msP} \prod_{j=1}^q (\eta f_j)^{n_j}\diff \DF{\alpha}(\eta) = \frac{\length{\mbfn}!}{\Poch{\beta}{\length{\mbfn}}}Z^\sym{q}_\mbfn\tparen{\Omega_\mbfn[f_1,\dotsc, f_q;\alpha]} \comma \qquad f_j\in\mcB_b\comma j\in [q]\fstop
\]
\end{cor}

We recover Theorem~\ref{t:Moments} by choosing a Borel partition~$\seq{X_i}_{i\leq k}$ of~$X$ with~$\alpha_i\eqdef \alpha X_i>0$, and simple functions~$f_j$, $j\in [q]$, constantly equal to, respectively,~$\mbfS_i^j$ on each~$X_i$ for all~$i\in [k]$.

\paragraph{Gamma measures} %
Denote by~$\diff G_{\beta,1}(t)= \Gamma(\beta)^{-1}\, e^{-t}t^{\beta-1}\car_{(0,\infty)}(t)\diff t$ the Gamma distribution on~$(0,\infty)$ with shape parameter~$\beta>0$ and scale parameter~$1$.
The \emph{gamma measure~$\GP{\alpha}$ with intensity} (\emph{measure})~$\alpha$ is the (unique) Borel probability on~$\Mbp$ satisfying, cf.~\eqref{eq:JIsomorphism},
\[
J_\pfwd \GP{\alpha}= \DF{\alpha}\otimes G_{\beta,1}\fstop
\]

Similarly to the case of~$\DF{\alpha}$, multivariate moments of~$\GP{\alpha}$ can be expressed in terms of~$Z^\sym{q}_\mbfn$.

\begin{cor}[Multivariate moments of~$\GP{\alpha}$] 
We have
\begin{equation}\label{eq:c:GammaMoments:0}
 \int_{\Mbp} \prod_{j=1}^q (\eta f_j)^{n_j}\diff \GP{\alpha}(\eta) = \length{\mbfn}! \, Z^\sym{q}_\mbfn\tparen{\Omega_\mbfn[f_1,\dotsc, f_q;\alpha]} \comma \qquad  f_j\in\mcB_b\comma j\in [q] \fstop
\end{equation}
\end{cor}

\subsection{Generating functions}
We conclude with computations of generating functions.
\begin{cor}
For every choice of~$\mbff\eqdef\seq{f_j}_{j\in[q]}$, with~$f_j\in\mcB_b$, and every~$\alpha\in\Mbp$,
\begin{align}\label{eq:c:EGFGamma:0.1}
\sum_{\mbfn\in\N_0^q} \frac{\length{\mbfn}!}{\mbfn!} Z^\sym{q}_\mbfn\tparen{\Omega_\mbfn[f_1,\dotsc, f_q;\alpha]} \mbfz^\mbfn &= \exp\bracket{-\int_X \log(1-\mbfz\cdot\mbff) \diff \alpha}\comma 
\end{align}
for every~$\mbfz\in \R^q$ with~${\alpha\text{-}\esssup \sum_j \abs{z_j f_j}} <1$.
Analogously, for every choice of~$k\in\N_1$, of~$\mbfs_1,\dotsc,\mbfs_q \in \R^k$, and of~$\boldalpha\in (0,\infty)^k$,
\begin{gather}\label{eq:c:EGFGamma:0.2}
\sum_{\mbfn\in\N_0^q} \frac{\length{\mbfn}!}{\mbfn!} \msZ^\sym{q}_\mbfn\bracket{\mbfS;\boldalpha} \mbfz^\mbfn = \prod_{i=1}^k \paren{1-\mbfz \cdot \mbfS_i}^{-\alpha_i} \comma \qquad \mbfz \in \R^q \, : \, \max_i \sum_j \tabs{z_j \mbfS_i^j} < 1 \fstop
\end{gather}
Finally, for every uniformly bounded~$\mbft\eqdef \seq{t_\mbfa}_{\mbfa\in \N^q_*}$, letting~$\mbft\restr_\mbfn\eqdef \seq{t_\mbfa}_{\mbfa\leq_\had\mbfn}$,
\begin{equation}\label{eq:c:EGFGamma:0.15}
\sum_{\mbfn\in\N_0^q} \frac{\length{\mbfn}!}{\mbfn!}Z^\sym{q}_\mbfn(\mbft\restr_\mbfn) \mbfz^\mbfn = \exp \bracket{\sum_{\mbfa\in\N^q_*} \frac{1}{\length{\mbfa}}\binom{\length{\mbfa}}{\mbfa} t_\mbfa \mbfz^\mbfa} \comma \qquad  \mbfz \in \R^q \, : \, \length{\mbfz} < 1 \fstop
\end{equation}
\begin{proof}
Computing the exponential generating function of the moments on the left-hand side of~\eqref{eq:c:GammaMoments:0} gives the known expression for the Laplace transform of~$\GP{\alpha}$, e.g.~\cite[Eqn.~(7)]{TsiVerYor01}, which shows~\eqref{eq:c:EGFGamma:0.1}. 
{(The condition~$\alpha\text{-}\esssup \sum_j \abs{z_j f_j} < 1$ ensures absolute convergence of the left-hand side of~\eqref{eq:c:EGFGamma:0.1}.)}
The identity~\eqref{eq:c:EGFGamma:0.2} follows from~\eqref{eq:c:EGFGamma:0.1} by specializing the latter to simple functions and applying Proposition~\ref{p:Z=Z}.
In order to show~\eqref{eq:c:EGFGamma:0.15}, set~$\xi_\mbfa\eqdef \frac{1}{\length{\mbfa}}\binom{\length{\mbfa}}{\mbfa} t_\mbfa \mbfz^\mbfa$ for simplicity of notation. 
{Note that~$\mbfz^\mbfn = \prod_{\mbfa} \mbfz^{\mbfa A(\mbfa)}$ whenever~$A \vDash \mbfn$.
By~\eqref{eq:ZDefinition} and~\eqref{eq:MultinomialMultiYoung},} %
\begin{equation}\label{eq:c:EGFGamma:1}
\sum_{\mbfn\in\N_0^q}\frac{\length{\mbfn}!}{\mbfn!} Z^\sym{q}_\mbfn(\mbft) \mbfz^\mbfn = \sum_{n\in \N_0} \sum_{A\in\Part{q}{n}}  \prod_{\mbfa\in\N^q_*} \frac{\xi_\mbfa^{A(\mbfa)}}{A(\mbfa)!} = \sum_{A\in \Part{q}{\cup}}  \prod_{\mbfa\in\N^q_*} \frac{\xi_\mbfa^{A(\mbfa)}}{A(\mbfa)!} \fstop
\end{equation}
Reindexing the last term over the labeling~$m\eqdef A(\mbfa)$,
\begin{equation}\label{eq:c:EGFGamma:2}
\sum_{A\in \Part{q}{\cup}}  \prod_{\mbfa\in\N^q_*} \frac{\xi_\mbfa^{A(\mbfa)}}{A(\mbfa)!} = \prod_{\mbfa\in\N^q_*} \sum_{m=0}^\infty \frac{\xi_\mbfa^m}{m!} = \prod_{\mbfa\in\N^q_*} e^{\xi_\mbfa} = \exp\bracket{\sum_{\mbfa\in\N^q_*} \xi_\mbfa}\fstop
\end{equation}
Combining~\eqref{eq:c:EGFGamma:1} and~\eqref{eq:c:EGFGamma:2} proves the assertion.
\end{proof}
\end{cor}

\begin{rem}
Alternative expressions for the multivariate moments of the Gamma measure may be obtained by differentiating its characteristic functional, see~\cite[p.~5]{LzDSLyt17}.
Such expressions are however not informative on the algebraic and combinatorial meaning of the moments in connection with~$Z^\sym{q}_\mbfn$, as they rather rely on the multivariate multi-factor Leibniz rule.
A similar approach does not apply to the Dirichlet--Ferguson measure, due to the convoluted form of its characteristic functional.
\end{rem}

\appendix

\section{An interpretation in population genetics%
}\label{ss:PopGen}

In this section, we discuss an interpretation of the pESF in population genetics in terms of \emph{in situ} epigenetic modifications.
We refer the interested reader to~\cite{Fel14} for a more complete account of epigenetic modifications and, in particular, of \emph{DNA methylation} accessible to a broad readership.
The latter is the attachment of a methyl group to the 5-carbon atom of a cytosine nucleotide base.
In the following, we will focus on DNA methylation since, among other epigenetic modifications, it happens \emph{in situ}, that is, at a given genomic location, without disassembly, relocation, or replacement of the underlying DNA segment.

\paragraph{Deletion uniformly at random}
If we set, conventionally,~$\probP_n[A-\car_{\zero}]=\probP_n[A]$ for each $q$-chromatic partition~$A\in\Part{q}{n}$, %
formula~\eqref{eq:ConsistencySpecial} further simplifies as
\begin{equation}\label{eq:ConsistencyReduced}
\probP_{n-1}[B] = \sum_{j=1}^q \sum_{\mbfb\in\supp B\cup\set{\zero}} \probP_n [B-\car_\mbfb+\car_{\mbfb+\mbfe_j}] \frac{(b_j+1) \tparen{B(\mbfb+\mbfe_j)+1}}{n} \fstop
\end{equation}
{This is in line with the original interpretation by Kingman in terms of allele sampling.}

In the polychromatic setting, we treat each color $j \in [q]$ as an \emph{in situ} epigenetic modification, or \emph{marker}.
In this setting, we use the term `\emph{in situ}' to indicate that {the} marker specification is subordinate to the allele; thus, at a fixed genetic locus, an individual is identified by both its allele type and the specific marker it displays for that allele.
As an example, a gene consisting of~$k$ base-pairs {displays one out of}~$4^k$ alleles. Each allele may have up to~$k$ cytosine bases, and each such base can be either methylated or non-methylated, so that an allele with~$k_{\rm{cyt}}$ cytosine bases ($k_{\rm{cyt}}\leq k$) may have up to~$q=2^{k_{\rm{cyt}}}$ marked versions.

When a sample of~$n$ individuals is partitioned by their combined allele-marker types, the result is a polychromatic partition~$A$.
In this framework, two individuals belong to the same block (with coloring $\mbfa$) if they share the same allele.
Within such a block, the component~$a_j$ is the number of individuals displaying the same marker type~$j$ for that specific allele.

If we remove one individual uniformly at random from our sample, the resulting polychromatic partition~$B$ is formed in two possible ways:
\begin{itemize*}
\item the removed individual is one of~$A(\mbfe_j)$ individuals, each with unique allele type and marker~$j$, that is, it forms a singleton in~$A$ with block type~$\mbfe_j$;
\item the removed individual belongs to one of~$A(\mbfa)$ groups of individuals sharing (within the same group) the same {allele type}, and of~$a_j$ individuals within that group sharing the same marker~$j$, that is, it belongs to one out of the~$A(\mbfa)$ blocks of the same type~$\mbfa$ for some~$\mbfa$ with~$\length{\mbfa}\geq 2$, and is one of the $a_j$~elements in that block with color~$j$.
\end{itemize*}
In this second case, we always have~$\mbfa=\mbfb+\mbfe_j$ for some~$\mbfb\in\supp B$, the vector~$\mbfe_j$ accounting for the individual we have removed.
Thus, formula~\eqref{eq:ConsistencyReduced} states that~$\probP_\bullet$ is a partition structure precisely when it is consistent under removal of one individual chosen uniformly at random from our sample of~$n$ individuals in the population.

In particular, for the pESF with~$\nu=\delta_\mbfp$, the color of each element is independent of the block it belongs to, that is, each allele can carry any marker, and markers are independent of their carrier alleles.
That each allele carries a marker is a harmless assumption, in that we can always interpret a fixed marker (say, the first color) as the unmodified version of the allele.

\paragraph{Coefficients}
Many combinatorial coefficients appearing throughout have a simple explanation in light of the above population-genetics interpretation.
For example, in the combinatorial coefficient in~\eqref{eq:DefPhi}, the factor~$\mbfa!^{A(\mbfa)}$ accounts for the fact that~$a_j$ individuals,~$j\in [q]$, carrying the same allele with the same epigenetic marker are indistinguishable within that allele-marker group:
permuting them (i.e.\ their labels) does not change which allele and marker they share.
The factor~$A(\mbfa)!$ accounts for the fact that allele-marker pairs are unlabeled in the polychromatic partition~$A$.

\paragraph{pESF vs.\ rESF}
The distinction between rESF and pESF is apparent from their interpretation in population genetics:
\begin{itemize*}%
\item In the pESF (with~$\nu=\delta_\mbfp$), the probability weight of a monochromatic block of color~$j$ and size~$a$ is proportional to $p_j^a$, because~$a$ individuals are each assigned color~$j$.
In this model, the color is like an epigenetic modification that can be acquired independently of the genetic lineage (the block).
{A lineage is `polychromatic' because blocks can acquire new colors.}

\item In the rESF, the probability weight of the same block is proportional to $\theta_j$, because the color was chosen only once when the allele was formed.
{In this model, the color is (an epigenetic modification of) the allele itself.}
The mutation occurs at the beginning of the lineage.
Every member of the lineage \emph{inherits} this color by descent, making the lineage (the block) monochromatic. 
\end{itemize*}

\subsection*{AI declaration}
The authors acknowledge the use of ChatGPT~5.6~Sol and Claude~Opus~5 as complementary tools for proofreading, literature research, and the simplification and preliminary verification of human-generated proofs.
All AI-generated comments incorporated into the manuscript have been independently checked and adapted by the authors.

{\small{
\bibliographystyle{abbrvurl}

}}

\end{document}